\documentclass[11pt,english]{article}
\usepackage[sc]{mathpazo}
\usepackage{geometry}
\usepackage[noblocks]{authblk}
\usepackage{graphicx}
\usepackage{booktabs}
\usepackage{threeparttable}
\usepackage{multicol}
\usepackage{multirow}
\usepackage{color}
\usepackage{algorithm}
\usepackage{tabularx}
\usepackage{bm}
\usepackage{setspace}
\usepackage{graphicx}
\usepackage{subfigure}
\usepackage{makecell}
\usepackage[]{natbib}
\usepackage{algorithmic}
\usepackage{amsmath, amssymb}
\usepackage[thmmarks, amsmath, thref]{ntheorem}
\usepackage[justification=centering]{caption}
\allowdisplaybreaks[4]
\usepackage[colorlinks=true, allcolors=blue]{hyperref}
\renewcommand{\algorithmicrequire}{\textbf{Input:}}

\newtheorem{proposition}{PROPOSITION}
\newtheorem{lemma}{LEMMA}
\newtheorem{definition}{DEFINITION}
\newtheorem{insight}{INSIGHT}
\theoremstyle{plain}
\theoremheaderfont{\bfseries}
\theorembodyfont{}
\theoremseparator{.}
\theoremsymbol{\ensuremath{\blacksquare}}

\usepackage{indentfirst} 
\usepackage{threeparttable}
\usepackage{longtable}

\usepackage{lmodern}

\begin{document}\sloppy
	
	\title{Robotic Sorting Systems: Robot Management and Layout Design Optimization}
	\author[a]{Tong Zhao} 
	\affil[a]{\small\emph{Department of Industrial Engineering, Tsinghua University, 100084 Beijing, People’s Republic of China}\normalsize}
	\author[b]{Xi Lin}
	\affil[b]{\small\emph{Department of Civil and Environmental Engineering, University of Michigan, Ann Arbor, Michigan 48109}\normalsize}
	\author[a]{Fang He\footnote{Corresponding author. E-mail address: \textcolor{blue}{fanghe@tsinghua.edu.cn}.}}
	\author[a]{Hanwen Dai}
	\maketitle
	
	\begin{abstract}
		\noindent In the contemporary logistics industry, automation plays a pivotal role in enhancing production efficiency and expanding industrial scale. In particular, autonomous mobile robots have become integral to modernization efforts in warehouses. One noteworthy application in robotic warehousing is the robotic sorting system (RSS), which is distinguished by its cost-effectiveness, simplicity, scalability, and adaptable throughput control. Previous research on RSS efficiency often assumed an ideal robot management system, ignoring potential traffic delays and assuming constant travel times. To address this gap, we introduce a novel robot traffic management method, named Rhythmic Control for the Sorting Scenario (RC-S), for RSS operations, along with an analytical estimation formula that establishes the quantitative relationship between system performance and configurations. Simulations validate that RC-S reduces average service time by 10.3\% compared to the classical cooperative A* algorithm, while also improving throughput and runtime. Based on the performance analysis of RC-S, we develop a layout optimization model that considers system configurations, desired throughput, and costs to minimize expenses and determine the optimal layout. Numerical studies show that facility costs dominate at lower throughput levels, while labor costs prevail at higher throughput levels. Additionally, due to traffic efficiency limitations, an RSS is well-suited for small-scale operations like end-of-supply-chain distribution centers.
		\hfill\break
		
		
	\end{abstract}
	\section{Introduction} \label{sec:intro}
Improved logistics and delivery services have fueled the rapid growth of e-commerce in the 21st century. The development of stable supply chains, same-day or next-day delivery options, and hassle-free return policies increases consumer confidence in online shopping \citep{ReginaLoBiondo}. To handle the increasing volume of online orders, a new generation of warehouses specifically catering to individual customers has become a hot topic for logistics companies. This type of warehouse efficiently meets the demand for small orders with tight delivery schedules through the implementation of automated equipment \citep{BOYSEN2019396}. For instance, the robotic mobile fulfillment system (RMFS), which employs a rack-moving mechanism, is widely used in intra-warehouse logistics, such as Amazon’s KIVA systems \citep{wurman2008coordinating}.

The main activities performed in a warehouse include: (1) receiving, (2) transfer and put-away, (3) order picking/selection, (4) accumulation/sorting, (5) cross-docking, and (6) shipping \citep{de2007design}. In this paper, we focus on the sorting process. Sorting involves categorizing and consolidating parcels according to their order information and shipping destinations. Conventional sorting systems commonly employ conveyor-based sorters , where the actuators move along with the conveyor belt, sequentially passing through each outbound station and releasing the loaded parcels at the appropriate locations \citep{BOYSEN2019796}. These systems are highly appreciated for their efficiency and stability, while suffering from the inflexibility and significant space occupation \citep{BOYSEN2023582}. 

\begin{figure}[hbt!]
	\centering
	\subfigure[A dual-layer RSS solution] {\includegraphics[width=7cm]{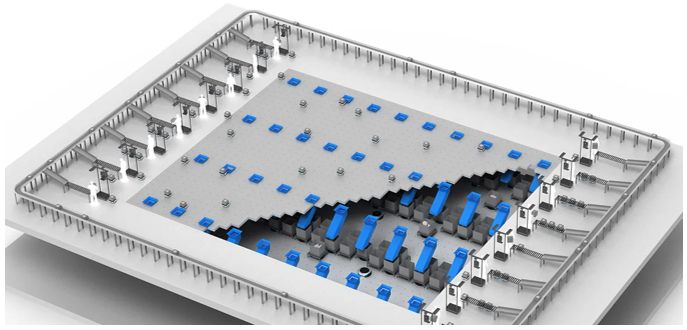}} 
	\subfigure[Conveyor-belt sorting robot] {\includegraphics[width=6cm]{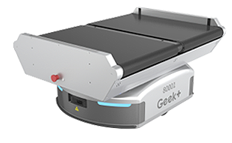}} 
	\caption{The application of RSS (source: www.geekplus.com)} 
	\label{fig:101}  
\end{figure}

A novel sorting system utilizing mobile robots attracts increasing attention in recent years, namely the robotic sorting system (RSS) \citep{zou2021robotic}. Figure~\ref{fig:101}(a) presents a dual-layer solution of RSS. The site is arranged with numerous delivery ports, each associated with a specific sorting category, gathering all parcels of that category. A delivery port can be a bin, a cage cart, or a chute in a dual-layer structure. Robots deliver parcels to delivery ports, serving as the combination of conveyor belts and actuators in conventional sorting systems. Loading stations are located on the periphery of the site, where workers load parcels onto the robots. A top-mounted tray or conveyor actuator enables the robot to load or unload parcels. Figure~\ref{fig:101}(b) shows a specific model of sorting robot developed by Geekplus, a global technology company specialized in smart logistics and robots. The process for a robot to execute a sorting task is as follows:
\begin{enumerate}
	\item The robot receives a parcel at a loading station, along with destination information;
	\item The robot moves along aisles to the designated outlet;
	\item Upon reaching the target outlet, the robot releases the loaded parcel;
	\item The robot returns to the waiting zone behind a loading station and becomes idle.
\end{enumerate}

E-commerce orders fluctuate during special selling seasons; therefore, warehouse throughput must be designed to accommodate these variations. Unlike conveyor-based sorting systems, the independent and modular nature of robotic operations contributes to the flexibility and scalability of RSS, making it well-suited to the dynamic demands of e-commerce \citep{pub99983}. During off-hours, warehouse managers can easily configure the control software to change the status of loading stations, the correspondences between outlets and destinations, and the aisle network topology \citep{XU2022102808}. To leverage this flexibility, robotics companies have introduced the innovative Robots-as-a-Service (RaaS) business model, allowing logistics clients to rent robots as needed. RaaS eliminates the high costs associated with purchasing and maintaining robots and peripheral equipment, enabling warehouse managers to adjust the number of rented robots based on demand fluctuations.

Despite these advantages, a research gap remains in investigating the flexibility of RSS and minimizing system costs, as accurately modeling system efficiency is challenging. In robotic systems, robots rely on path planning and coordination algorithms to complete tasks. The problem of finding conflict-free trajectories for all robots is known as the multi-agent path finding (MAPF) problem. The MAPF problem has been proven to be NP-hard \citep{Yu_LaValle_2013}, which means that finding a relatively good solution can be time-consuming, especially when the scale of the robot system increases, resulting in instability of RSS efficiency. \cite{poms13626} analyzed a dataset from the China Post sortation center and found that the congestion effect significantly hampers robot efficiency in RSS, particularly when dealing with a substantial parcel flow. \cite{zou2021robotic} proposed several closed queueing network (CQN) models to quantitatively analyze the impact of congestion and designed an algorithm to estimate throughput. Validated through numerical experiments, these CQN models accurately reflect the performance of a real case from Deppon Express. However, the models primarily focus on the process of agents transferring between queues rather than on traffic flow and do not fully consider traffic management methods to mitigate conflicts. This makes them less suitable for analyzing the efficiency of MAPF solvers, such as the method proposed in this paper. On the other hand, other existing studies often overlooked issues related to robot coordination, relying on the idealized assumption that there are no conflicts or deadlocks. Consequently, their results may become distorted as the actual travel time of robots deviates from the free-flow travel time. Furthermore, these studies do not adequately explore the flexibility of RSS, particularly the ability to adjust the number of robots.

Given the existing gap, this paper aims to address critical challenges in designing an efficient RSS, focusing on three key areas. First, we propose a novel robot traffic management method to enhance system throughput, ensuring stable and efficient operations in multi-robot collaboration settings while overcoming the complexities of robot interactions. Second, we develop a precise model to quantify system efficiency in RSS, enabling accurate throughput estimation to reveal the relationship between system configurations and sorting capacity. Third, we highlight the need for management strategies that leverage RSS flexibility to provide cost-effective responses to fluctuating demands. These contributions aim to improve the performance and scalability of RSS implementations.

To address the aforementioned challenges, we first explore the potential of centralized control strategy for sorting robots. To fully utilize global information while considering robot operational patterns and safety distance regulations, we adapt a network-level control strategy for coordinating autonomous vehicles and propose a collision-free spatio-temporal path planning method for RSS. An efficient heuristic algorithm is introduced to enhance the scalability. Second, building upon this control strategy, we derive an estimation formula for throughput in RSS given system configurations, which serves as an efficiency constraint in the planning problem. Finally, we propose an optimization model that aims to minimize initial investment and average operational costs under fluctuating demands. The solution includes layout design and resource allocation recommendations, specifically on how to adjust the number of robots and workers. Experiments provide insights into the cost structure of RSS under varying unit price conditions. To the best of our knowledge, this paper is among the first to incorporate robot traffic issues during operations into RSS configuration design. Our intention is to bridge operational and strategic planning in multi-robot coordination contexts, thereby contributing a comprehensive theoretical framework to support decision-making in such complex and flexible systems.

The remainder of this paper is organized as follows. In Section \ref{sec:LR}, we provide an extensive review of relevant research on RSS systems and conventional sorting systems, highlighting the contributions of this research. Section \ref{sec:description} offers a detailed problem description. In Section \ref{sec:RC-S}, we introduce an innovative traffic management framework to coordinate multiple robots in a warehouse setting. Based on the proposed framework, Section \ref{sec:analysis} presents a system throughput estimation formula, which quantifies the impact of layout configuration on system efficiency. Section \ref{sec:validation} demonstrates the superiority of our proposed traffic management framework compared with benchmark methods and validates the accuracy of the throughput estimation formula. To emphasize the impact of traffic issues on RSS, experiments are conducted to show that queueing network model has estimation biases for throughput in certain scenarios. Moving on to Section \ref{sec:layout}, we propose a layout optimization model to minimize total costs and introduce an efficient solution approach. Section \ref{sec:numerical} presents a sensitivity analysis to investigate the optimal layout design and the corresponding cost structure under different unit costs, and insights distilled from the results. Finally, section \ref{sec:conclusion} provides a summary of the entire paper.
	\section{Literature Review} \label{sec:LR}
The cost savings in labor through automation are particularly crucial in the picking process, and improved sorting efficiency enables parcel companies to offer convenient and low-cost same-day delivery services \citep{dekhne2019automation}. However, the operations research community has scarcely explored the realm of robotized sorting systems, which represent the latest advances in warehouse automation. Among the existing literature, studies focus mainly on the efficiency of small parcel sorting systems in distribution centers. This body of work delves into the comprehensive assessment of transportation efficiency within the robot fleet under various system configurations. In the majority of these studies, robot sorting tasks are modeled as services in queuing network models, and total throughput is the key performance evaluation metric  \citep{9532245,9216812,9102081,zou2021robotic,XU2022102808,fang2025dynamic}. \cite{9216812,9102081} both assumed that each robot moves independently. \cite{9532245} introduced a semi-open queuing network (SOQN) model and subsequently improved its theoretical aspects in their follow-up research \citep{zou2021robotic}, simplifying it to the closed queueing model (CQN) by assuming a sufficiently high parcel arrival rate. Their queueing model took into account the queuing behavior of robots at loading stations and outlets, partially addressing the issue of robot occupancy on the aisle. Furthermore, their second study compared how the minimum total costs of the system varied under different network topologies and different unit costs of system components, given the target sorting efficiency. \cite{fang2025dynamic} further emphasized the importance of dynamic robot routing and destination assignment policies in handling severe heterogeneity and fluctuations in sorting demand. They proposed an SOQN model to analyze network congestion and an MIP model to reassign the destination-outlet mapping in response to drastic changes in sorting demand. \cite{XU2022102808} investigated an RSS with a parcel-to-loading-station assignment mechanism, which could be regarded as a pre-sorting strategy. This study removed the constraint on the number of robots, resulting in an open queueing network (OQN) model. Experimental results demonstrated that the introduction of pre-sorting effectively reduces the average travel distance of the robots but might simultaneously lead to additional upstream congestion.

\cite{liu2019multi} and \cite{TAN2021101279} focused on the task assignment problem with the objective of minimizing the sorting makespan. Both adopt a travel time model based on the closed-form travel time expression proposed by \cite{pub99983}, and the latter formulates a mixed-integer programming model. Compared to queuing network models, travel time models exhibit significant deviation when the number of robots is high, as they allow multiple robots to simultaneously occupy a facility. \cite{BOYSEN2023582} conducted research on an order sorting system, which involves additional tasks including Piece-to-Order assignment and Order-to-Collection Point assignment. Therefore, the primary focus of their study is not solely on scheduling and control of robots, and the travel time of robots to a specific outlet is assumed to be constant. \cite{poms13626}, on the other hand, investigates a human-robot hybrid sorting system, in which the scale of robot sorting is fixed, while the manual capacity is adjustable according to the demand. Regrettably, in this context, the considerable flexibility inherent in the robotic system has not been fully exploited.

\begin{table}[!hbt]
	\centering
	\caption{Literature summary on RSS investigation}
	\label{tab:table11}
	\begin{threeparttable}
		\begin{tabular}{ccccc}
			\toprule
			Reference &System & Modeling & Objective & Congestion \& Deadlock\\
			\midrule
			\cite{liu2019multi}      & parcel sortation & TT, MIP & MS, C & $\times$\\
			\cite{9532245}      & parcel sortation & SOQN & T, C & $\circledcirc$\\
			\cite{9216812}      & parcel sortation & Q & T & $\times$\\
			\cite{9102081}      & parcel sortation & Q, MIP & T, C & $\times$\\
			\cite{zou2021robotic}      & parcel sortation & CQN, MIP & T, C & $\circledcirc$\\
			\cite{poms13626} & parcel sortation & TT, MIP & T & $\circledcirc$\\
			\cite{TAN2021101279}      & parcel sortation & TT, MIP & MS & $\times$\\
			\cite{XU2022102808}      & parcel sortation  & OQN & T & $\times$\\
			\cite{BOYSEN2023582}      & order sortation  & TT, MIP & T & $\times$\\
                \cite{fang2025dynamic}    & parcel sortation & SOQN, MIP & T, C & $\circledcirc$\\
			This paper & parcel sortation & TT, MIP & T, TC & $\checkmark$\\
			\bottomrule
		\end{tabular}
		\begin{tablenotes}
			\footnotesize
			\item[] \emph{Note.} TT, travel time model; MIP, mixed integer programming; Q, queueing model; CQN, closed queueing network model; OQN, open queueing network model; SOQN, semi-open queueing network model; MS, makespan; C, part of system costs; T, throughput; TC, total system costs; $\circledcirc$,  partially addressing the issue of robot occupancy on the aisle.
		\end{tablenotes}
	\end{threeparttable}
\end{table}

Throughout the existing research, a limited number of studies consider the floor space cost while others relax the constraints on aisle resources. Furthermore, finding collision-free paths for multiple vehicles has been proved to be an NP-hard problem \citep{Surynek_2010}, and some state-of-the-art algorithms also struggle to completely eliminate the issue of aisle congestion resulting from an increasing number of robots. However, except for the study by \cite{poms13626}, which employs a fitting method to predict the positive correlation between traffic flow and congestion, most of the studies neglect the mutual influence among multiple robots in aisle occupation for a compact system, leading to a misestimate of traffic capacity. Table~\ref{tab:table11} presents a comprehensive overview of the literature reviewed on the RSS, highlighting its relevance to our paper. 

Some researchers recognized the limitations of queuing models in capturing the impact of aisle layout design on multi-robot systems and chose to address the traffic issues by modeling them as MAPF problems. \cite{wagner6095022} proposed a novel M* algorithm that allowed for the planning of a larger number of robot trajectories. However, they faced challenges due to significant memory requirements and time complexity, and could not guarantee a feasible solution. Some researchers designed priority-based algorithms and integrated them with task assignment to reduce waiting time \citep{nguyen2019generalized,Liu2019path}, sacrificing optimality for improved stability and computational efficiency. \cite{Li_Tinka2021} introduced an RHCR algorithm framework to tackle the windowed MAPF problem, providing an efficient solution for managing newly arrived orders in warehouse operations. However, despite their effectiveness, MAPF-based studies have faced challenges in deriving a closed-form expression for throughput estimation, which is crucial for creating reliable and efficient system configurations. Consequently, few studies on multi-robot systems comprehensively address both operational robot management and strategic layout design while considering their interplay.

In summary, existing research on RSS systems lacks sufficient focus on robot traffic management and overall system costs. This study addresses these gaps by first developing a dedicated robot traffic management method for RSS at the operational level. Building upon this method, we establish an optimal layout and configuration design model aimed at minimizing overall system costs in the long run, thus enhancing the applicability of the findings to real-world logistics operations.
	\section{Problem Description} \label{sec:description}
\begin{figure}[!htp]
	\centering
	\includegraphics[width=0.45\textwidth]{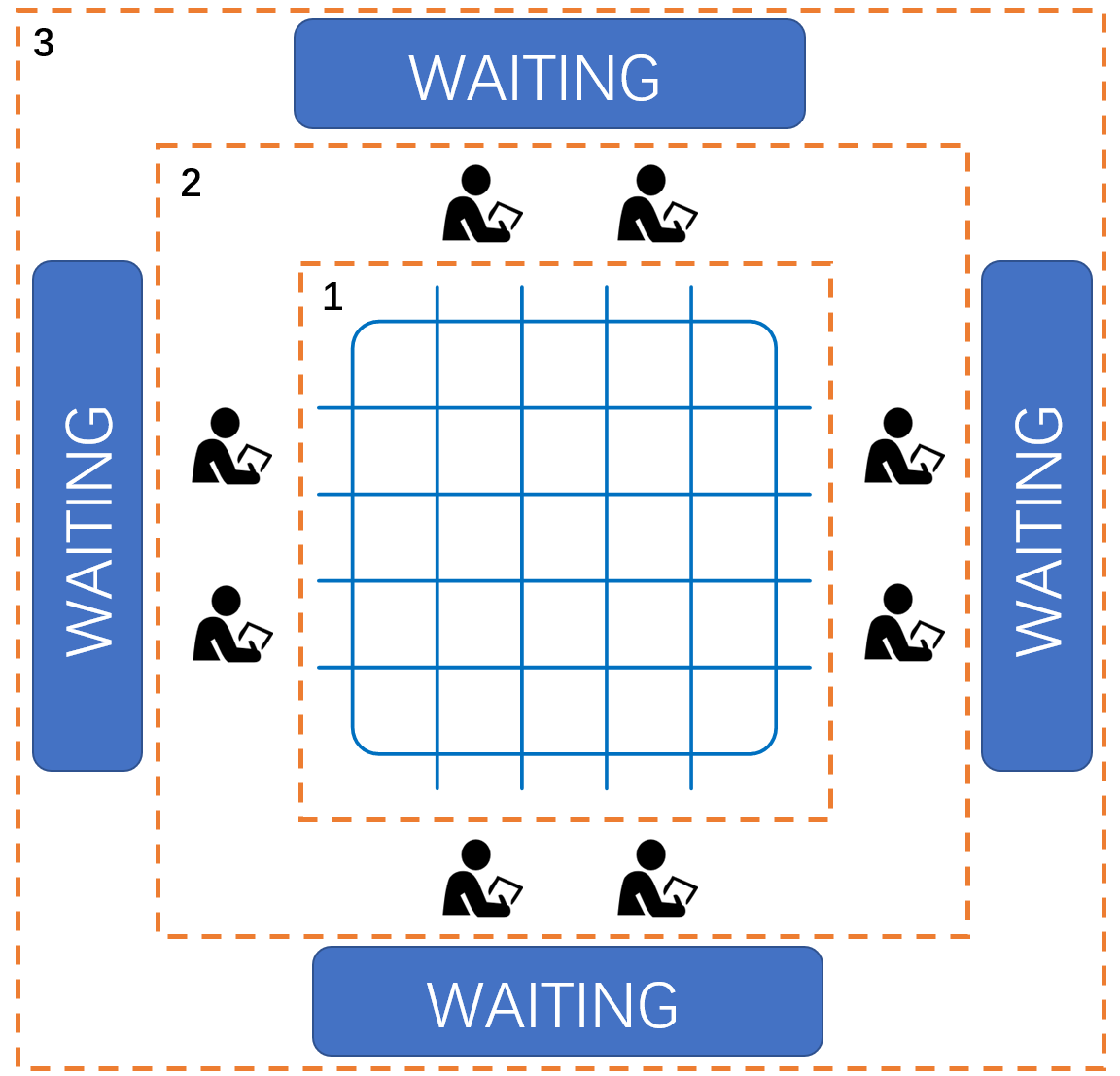}
	\caption{\label{fig:component} Layout of RSS}
\end{figure}

In an RSS, three zones are deployed from center to the periphery of the field \citep{10.1007/978-981-16-8174-5_10}: 1) sorting zone; 2) loading zone; 3) waiting zone, as shown in Figure~\ref{fig:component}. Parcels are loaded to robots by workers in the loading zone. Robots drop off each parcel at its target outlet in the sorting zone to fulfill a sorting task. The waiting zone is the space where robots idle and queue behind loading stations. The area of every zone is determined during the design phase and remains unchanged regardless of sorting demands or the numbers of robots and workers. Let $W_l$ and $W_w$ denote the widths of the loading zone and the waiting zone, respectively, both considered constant in this study. The decision of site planning mainly concerns the configuration of the aisle network in the sorting zone.

Without loss of generality, we make the following assumptions in this study:
\begin{itemize}
    \item The parcels arrive in batches. When a parcel is loaded onto a robot, the next parcel will arrive at the loading station immediately. That is, the loading time is negligible and is implicitly captured in the loading rate of the workers. 
    \item The arrival rate of parcels destined for each outlet at each loading station is identical. (When sorting demand is unevenly distributed, this can be achieved by adjusting the number of outlets corresponding to each destination. Additionally, outlets with extremely low or high arrival rates may be grouped together to balance the flow.)
    \item Robots release parcels without slowing down or stopping, and drop-off time is considered negligible. (To satisfy this assumption, robots must release parcels in advance based on their speed before reaching the target outlet, similar to the behavior of conveyor-based sorters.); 
    \item Drop-off operations around an outlet are independent, i.e., each outlet allows at most four robots releasing parcels simultaneously.
    \item Each robot can carry only one parcel during one delivery process;
    \item Once a robot exits the network through an exit node, it immediately joins the queue of the associated loading station without repositioning between loading stations. Travel time within the waiting and loading zones is neglected.
    \item The transportation of parcels after they are collected at the outlet is not considered.
    \item Battery charging requirements of robots are not considered.
\end{itemize}

\begin{figure}[hbt!]
	\centering
	\subfigure[Sample gird layout of a compact $4+4$ aisle network] {\includegraphics[width=7.5cm]{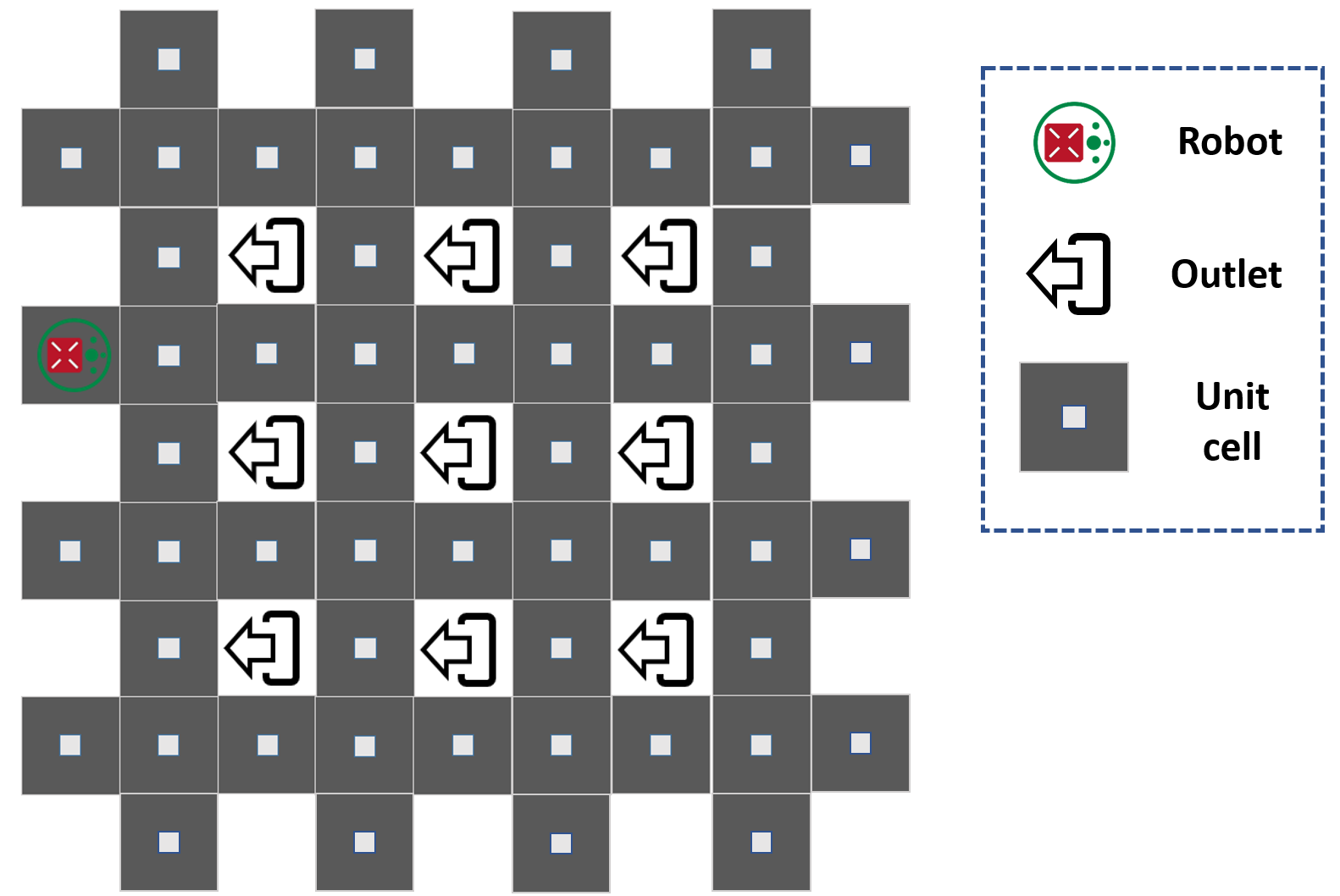}} 
	\subfigure[Graph representation of sample layout] {\includegraphics[width=7.5cm]{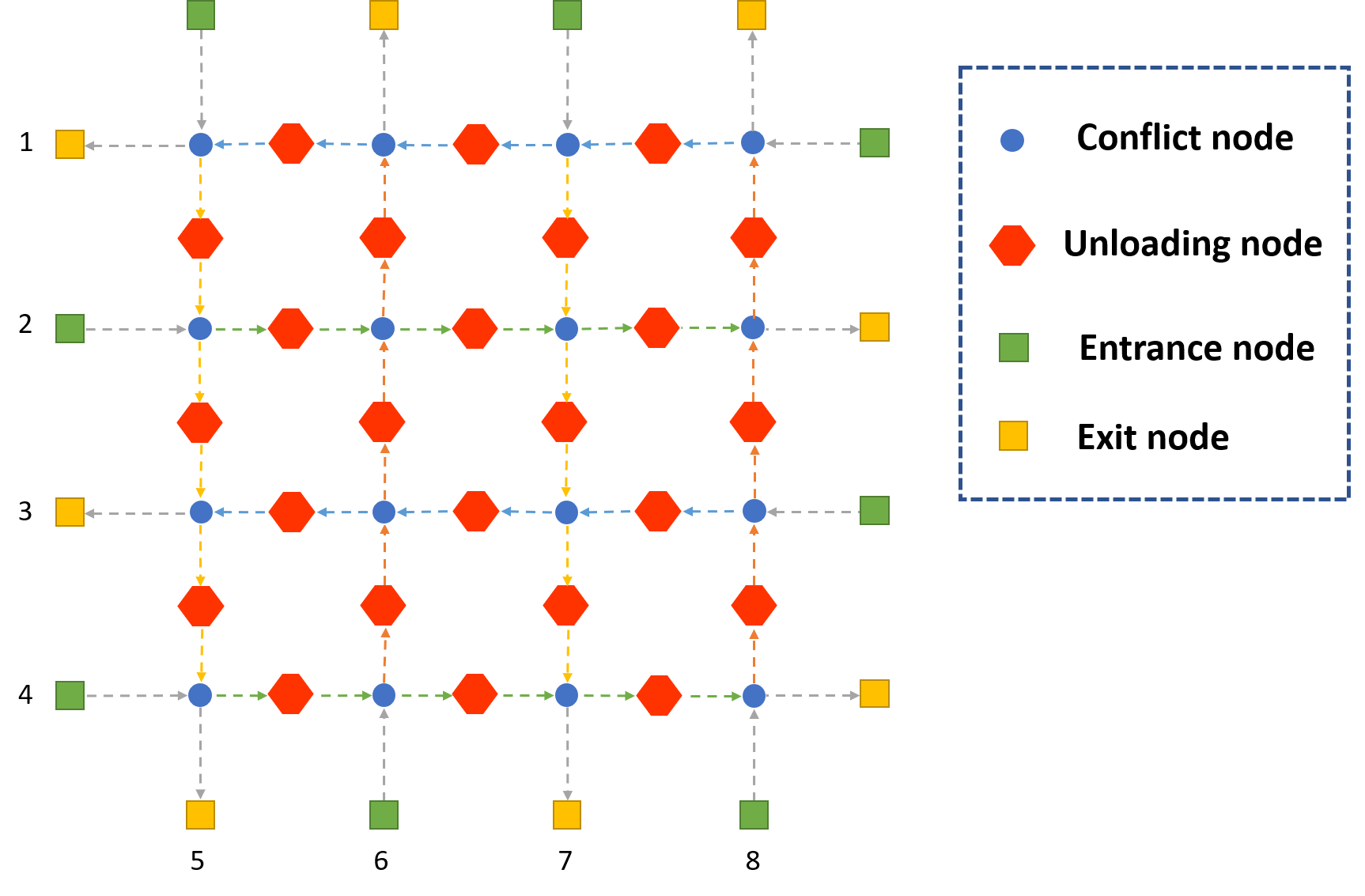}} 
	\caption{Illustration of aisle network in the sorting zone} 
	\label{fig:layout}  
\end{figure}

To better evaluate the performance of this system with high loads and complicated sorting demands, the microscopic traffic management of robots in the sorting zone needs careful investigation. Notably, we use the term “aisle” instead of “road” in the following content to clarify that each road in the network is single-lane. A sample grid-based map in RSS is depicted in Figure~\ref{fig:layout}(a). It consists of $n_{v}$ vertical aisles and $n_{h}$ horizontal aisles. Each cell in the map has a side length of $D$. Robots locate themselves using QR codes at the center of each cell and can only move between adjacent cells. We assume that the length of sorting robots is less than $D$, allowing them to rotate within a cell. We use the term “outlet” to refer to the special cells designated for collecting parcels delivered by robots. Outlets are evenly distributed throughout the sorting zone, seperated by aisles. This layout allows robots to deliver parcels to outlets on both sides of the aisle and ensures that the outlets are arranged with sufficient density. 

Each loading station in the loading zone is connected to a pair of adjacent entrance and exit nodes. Therefore, the number of loading stations cannot exceed the number of aisles. A dedicated worker is required at each loading station to perform parcel loading operations. A loading station is defined as \textit{active} if a worker is assigned and engaged in tasks at that station. Figure~\ref{fig:layout}(b) shows the graph representation of the sample layout in Figure~\ref{fig:layout}(a). In the graph, each cell that makes up the aisles is modeled as either a conflict node (blue dot), an unloading node (red hexagon), or an entrance/exit node (green/yellow square), with edges of length $D$ connecting adjacent nodes. Conflict nodes are at the intersections of aisles, unloading nodes are where robots deliver parcels, and entrance/exit nodes are at aisle ends. 

\begin{figure}[!ht]
	\centering
    \includegraphics[width=0.9\textwidth]{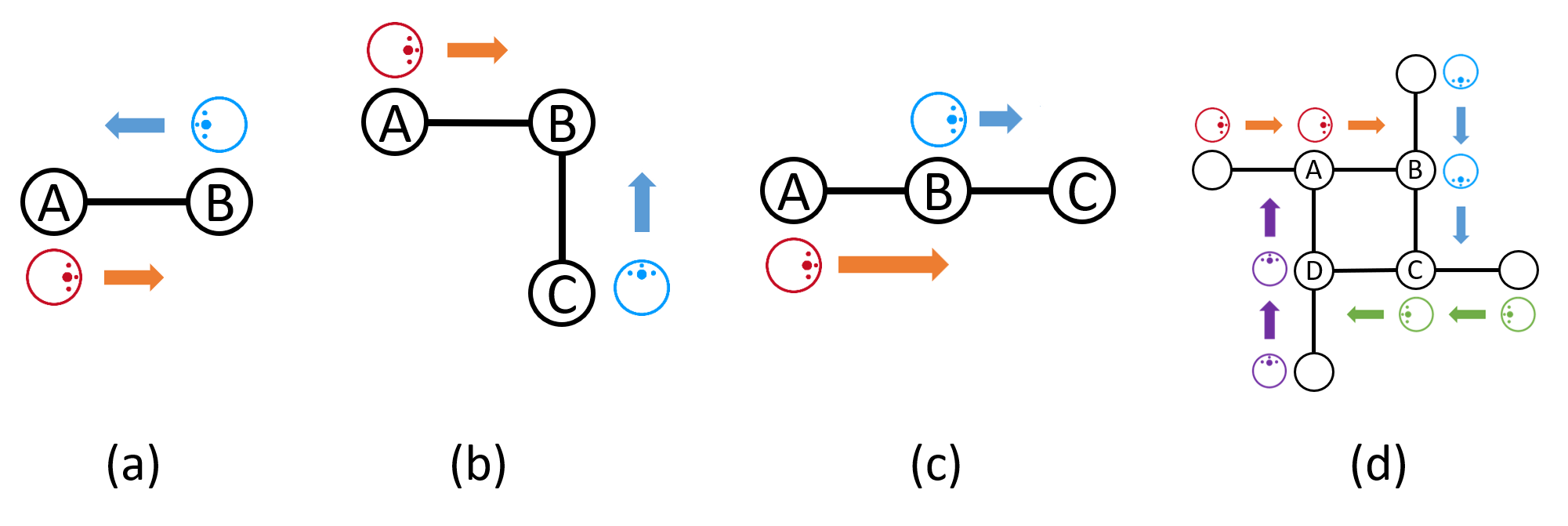}
	\caption{\label{fig:conflicts} Common types of conflicts: (a) a swapping conflict, (b) a vertex conflict, (c) a following conflict, and (d) a deadlock conflict}
\end{figure}

The sorting robots share the same kinematic parameters, with a maximum velocity of $v_{max}$ and maximum acceleration/deceleration of $c_{max}$. Robots possess the capability of in-place rotation, with angular velocity denoted as $\omega_{r}$. To ensure collision-free movement of robots, the robot management system needs to avoid conflicts in path planning. Common conflicts, as illustrated in Figure~\ref{fig:conflicts}, include the following four types \citep{stern2019multi}: (a) Swapping conflict occurs when two robots are planned to swap locations at the same time. It does not exist in a one-way network; (b) Vertex conflict occurs when more than two robots are planned to occupy the same vertex; (c) Following conflict occurs when one robot is planned to occupy the vertex currently occupied by another robot, while the latter is either stationary or moving at a lower speed; (d) Deadlock conflict involves multiple robots, each trying to maintain the minimum headway with another, resulting in mutual blockage. Among these conflicts, both swapping conflicts and deadlock conflicts necessitate a re-planning of paths. The occurrence of either can result in significant losses in system efficiency.

In this study, the RSS layout planning problem is structured as a two-stage cost-minimization model subject to throughput constraints, integrating long-term layout design with operational-level traffic control of robots:
\begin{enumerate}
    \item Site planning stage (strategic level):
    The objective is to minimize the total cost, including facility cost $C_f$ and expected operations cost $C_o$ across different operating periods $\sigma \in\mathcal{S}$:
    \begin{align}
	&\min\ C_{d} = C_f + C_o
    \end{align}
    Each period $\sigma$ is associated with a required throughput $T^\sigma$, representing changes in sorting demand. Decision variables include the number of horizontal aisles $n_h$, vertical aisles $n_v$, loading stations $n_l$, and the number of workers $n_w^\sigma$ and robots $n_r^\sigma$ allocated in each period $\sigma$. Layout and resource configurations must meet the throughput requirements:
     \begin{align}
	&\tilde{T}_{O}(n_h,n_v,n_w^\sigma,n_r^\sigma) \geq T^\sigma \qquad \forall \sigma \in \mathcal{S}
    \end{align}
    where $\tilde{T}_O$ is the estimated system throughput. A trade-off exists between $C_f$ and $C_o$: larger networks increase facility costs and travel distances, while smaller layouts may cause congestion and reduce network capacity.
    \item Operations stage (operational level):
    Given a fixed layout and resource configurations, the goal is to coordinate robot traffic efficiently to achieve stable and predictable system throughput. The robot traffic management method regulates robot movements and trajectory assignments, directly determining the system throughput. The estimation of system throughput $\tilde{T}_O$ serves as a critical link between operational performance and upper-level planning.
\end{enumerate}

\begin{figure}[!htp]
	\centering
    \includegraphics[width=0.9\textwidth]{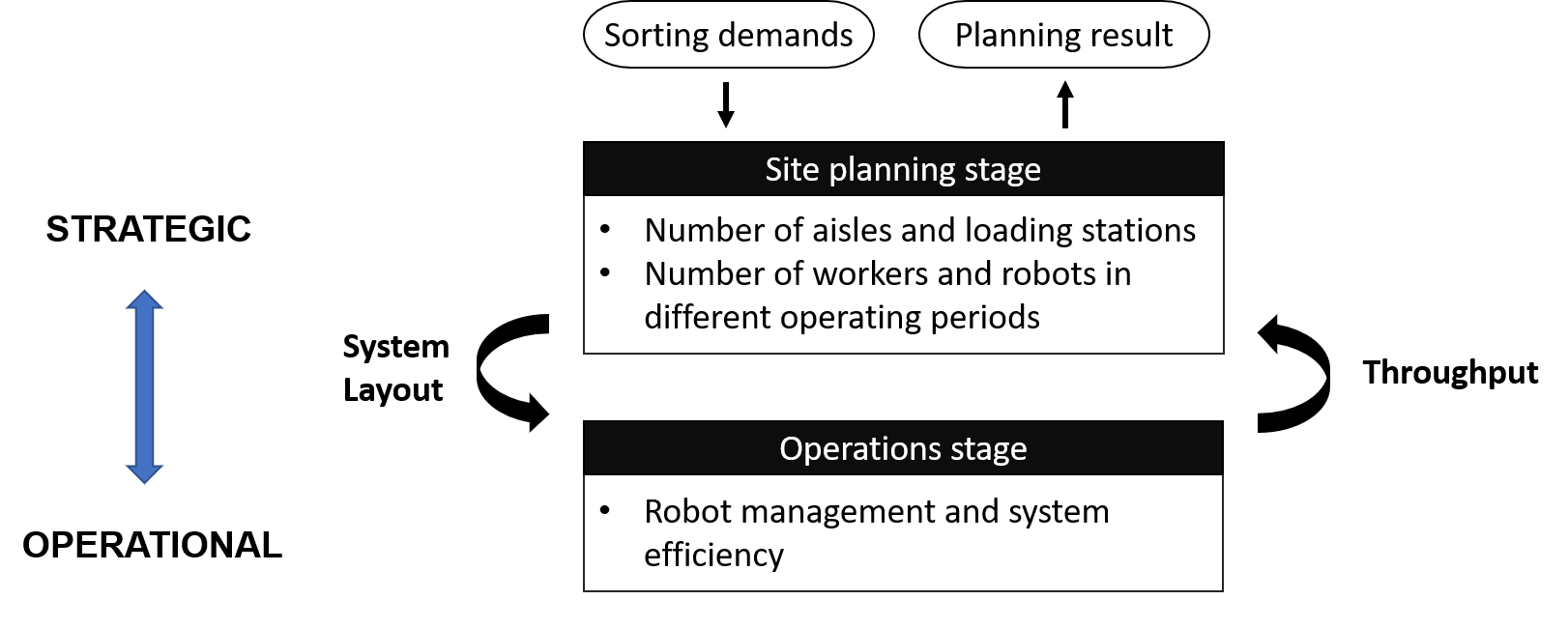}
	\caption{\label{fig:methodology} Overview of research methodology}
\end{figure}

Figure~\ref{fig:methodology} provides an overview of this integrated framework. In Sections \ref{sec:RC-S}–\ref{sec:validation}, we introduce the robot traffic management approach, analyze the system efficiency, and validate the performance and accuracy of the proposed throughput estimation method. In Section \ref{sec:layout}, we formulate the layout optimization model to minimize the total cost of RSS under given sorting demand.

	\section{Management of High-density Robot Traffic in RSS} \label{sec:RC-S}

Managing a large number of robots within a dense network poses significant challenges, especially when balancing system throughput with the computational load of algorithms under unpredictable traffic conditions. To address the multi-robot path finding problem and eliminate traffic gridlock within an acceptable computational time, this study adopts an innovative autonomous vehicle management scheme called Rhythmic Control (RC) \citep{RCII}, which enables uninterrupted scheduling of robot fleets. In this section, we begin with the key concepts inherited from the original RC: the virtual platoon and the cycle. By leveraging the reservation mechanism of RC, we propose a new centralized framework called Rhythmic Control for the Sorting Scenario (RC-S). RC-S exhibits a high degree of orderliness and serves as the foundation for the theoretical analysis of system efficiency, which will be discussed in detail in Section \ref{sec:analysis}. 

\subsection{Virtual platoon and cycle}
\begin{figure}[t!]
	\centering
	\subfigure[Virtual platoons in a four-lane network] {\includegraphics[width=7.5cm]{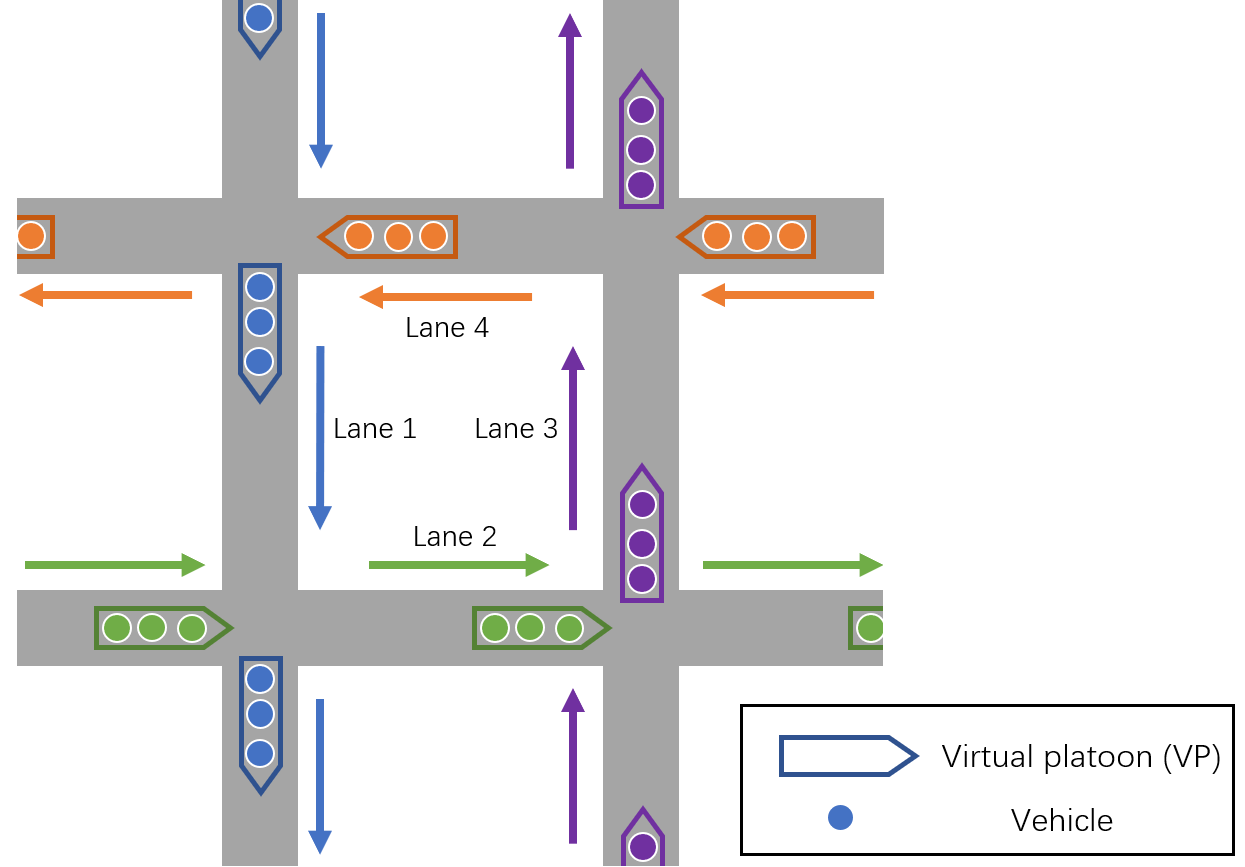}\label{fig:RC_1}} 
	\subfigure[Trajectory of virtual platoons in different directions. \citep{RCII}] {\includegraphics[width=6cm]{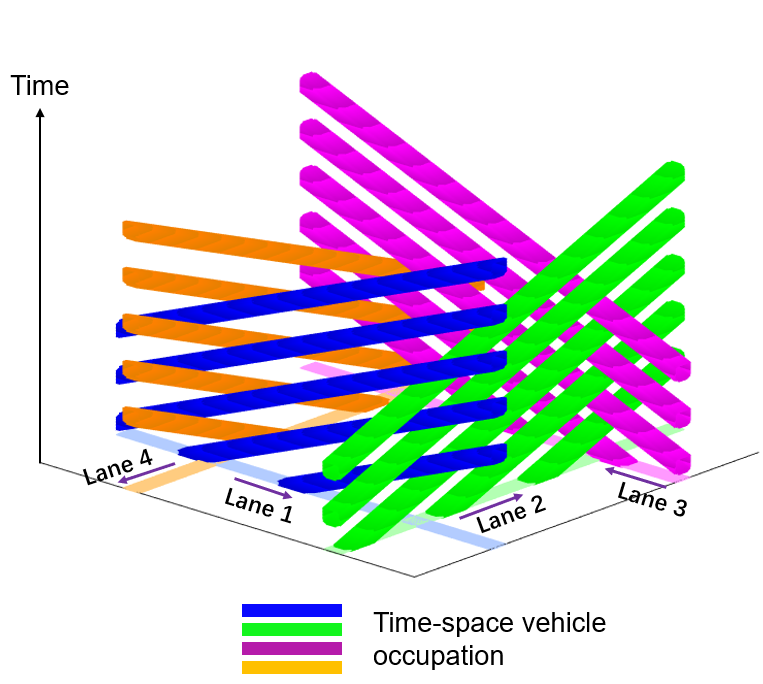}\label{fig:RC_2}} 
	\caption{Concept of Rhythmic Control} 
	\label{fig:RC}  
\end{figure}

\cite{RCII} investigated a method of incorporating rhythm into traffic management at the macro-traffic level and introduced the concept of \textit{virtual platoon} (VP). A VP is a reserved spatio-temporal slot along a lane, characterized by a uniform linear trajectory and generated at fixed intervals. Each VP could be identified by its corresponding lane and the interval in which it is released. As shown in Figure~\ref{fig:RC}(a), VPs move along a designated spatio-temporal trajectory and maintain a constant safe distance between each other. Figure~\ref{fig:RC}(b) shows the spatio-temporal trajectory of VPs in a four-lane network, without any collision or deceleration. We discretize time based on the generation time of VPs at the entrance, and each time interval is referred to as a \textit{cycle}. The duration of each time interval is the cycle length. By staggering the entry times of VPs at different entrances within the same cycle, VPs pass a intersection alternately with relatively short headway and the intersection capacity could be maximized.

Whenever a vehicle needs to pass through a lane, it is required to follow the movement of a VP in that lane until it reaches its destination and exits the lane. Before entering the network, the vehicle must wait for an unoccupied VP at the entrance. It first sends its destination to the control system, which then plans a conflict-free route and reserves one unit of capacity in the corresponding VPs in different segments of the route. The capacity is released when the vehicle leaves each VP. Through the reservation mechanism, RC ensures that the number of vehicles within each VP does not exceed its predetermined capacity. Consequently, the trajectories of vehicles are conflict-free.

RC is essentially a traffic control method, similar to a combination of globally coordinated traffic signal control and road network flow control. Within this framework, it is essential to clarify the movement patterns of each entity so that RC can be integrated into a sorting-specific robot management system. For detailed information on the RC framework, we refer the reader to \cite{RCII}. The following points provide a foundational overview that is critical to our study on the implementation of the RC-S scheme:
\begin{itemize}
    \item A virtual platoon (VP) represents a spatio-temporal slot that can be reserved and occupied by robots and is generated in a rhythmic manner at each entrance of the network.
    \item VPs move in a straight line at a constant speed while maintaining a fixed spacing, ensuring that the total number of VPs present in the aisle network remains constant.
    \item After entering the network, a robot must continuously follow a specific VP or transition to another VP until it exits.
    \item The capacity of each VP is one, which means at most one robot can occupy it at any time.
    \item A robot may enter the network only when the required VP capacity is available; otherwise, it must wait for a later cycle.
\end{itemize}

\subsection{Concept of RC-S} \label{sec:RC-S details}
\begin{figure}[!h]
	\centering
	\subfigure[] {\includegraphics[width=0.4\linewidth]{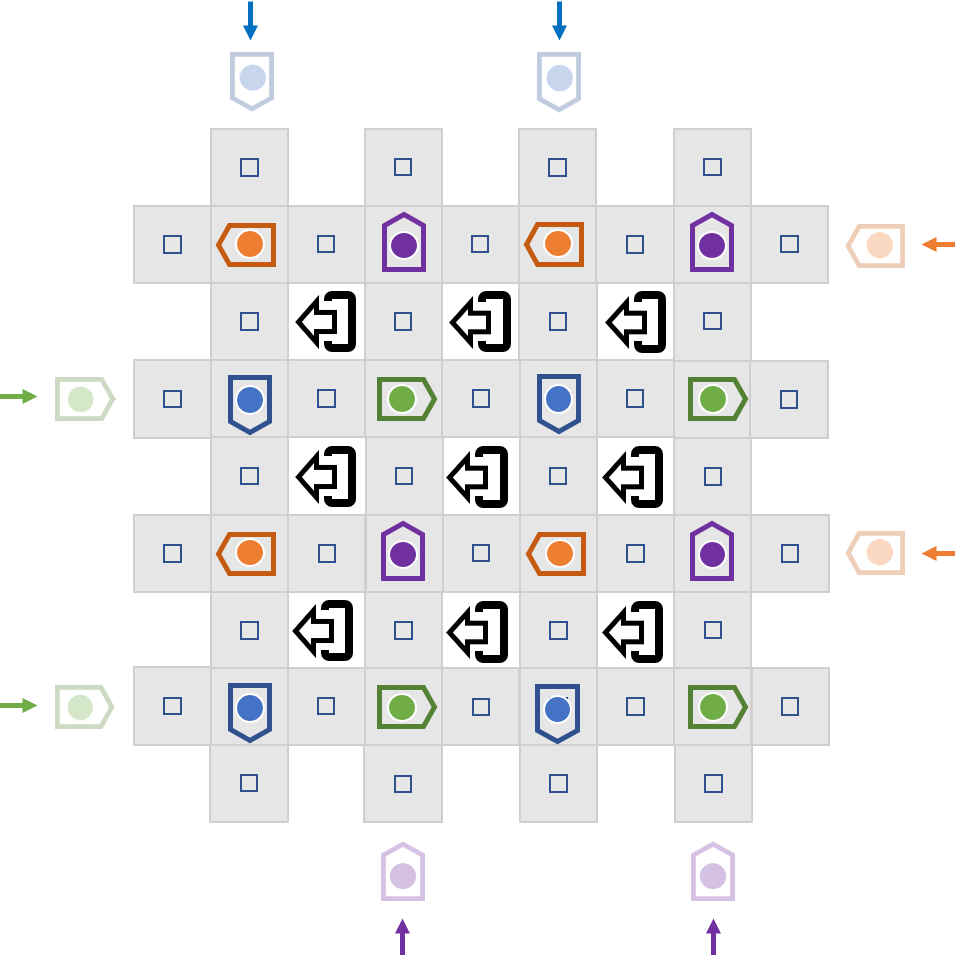}} 
	\subfigure[] {\includegraphics[width=0.4\linewidth]{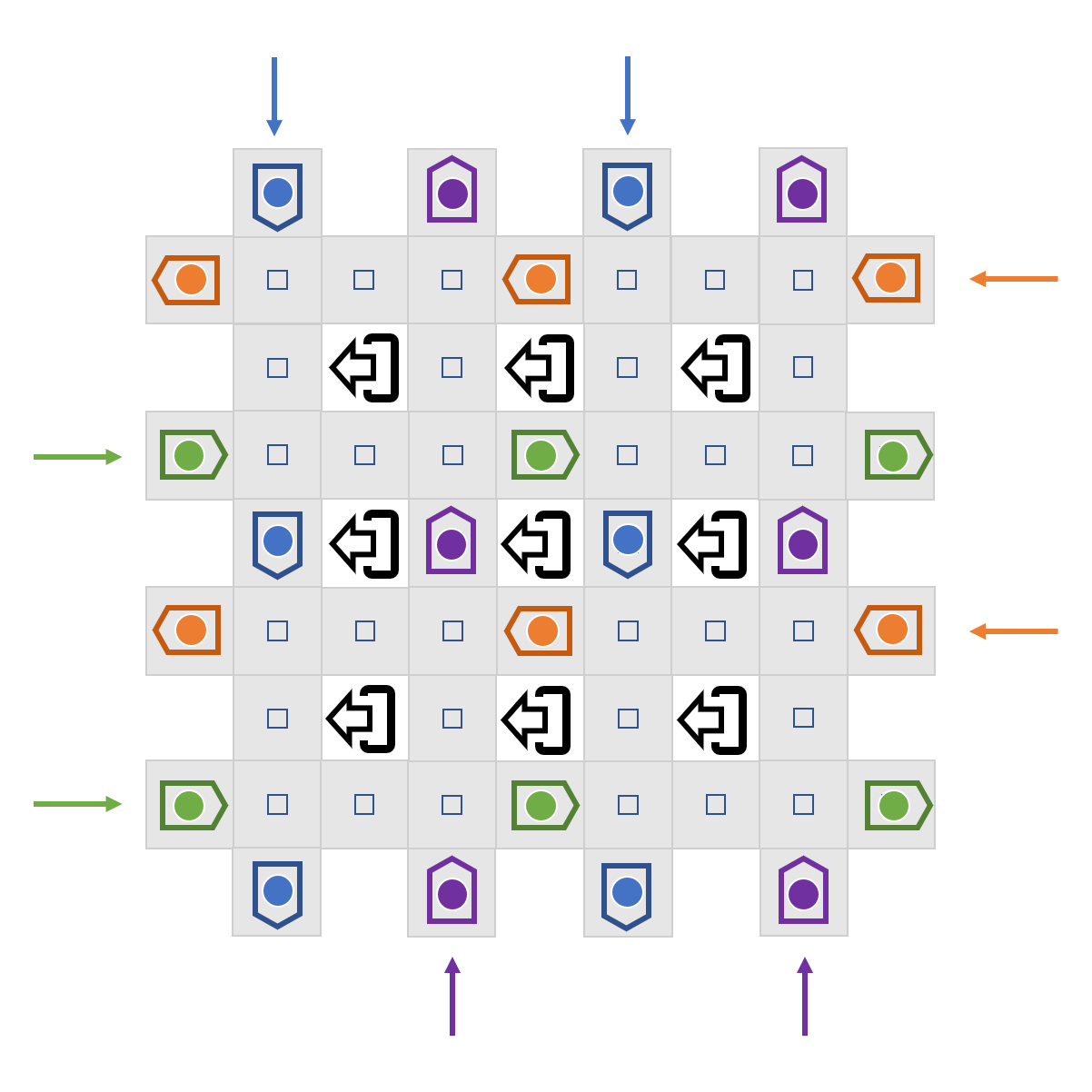}} 
	\qquad
	\subfigure[] {\includegraphics[width=0.4\linewidth]{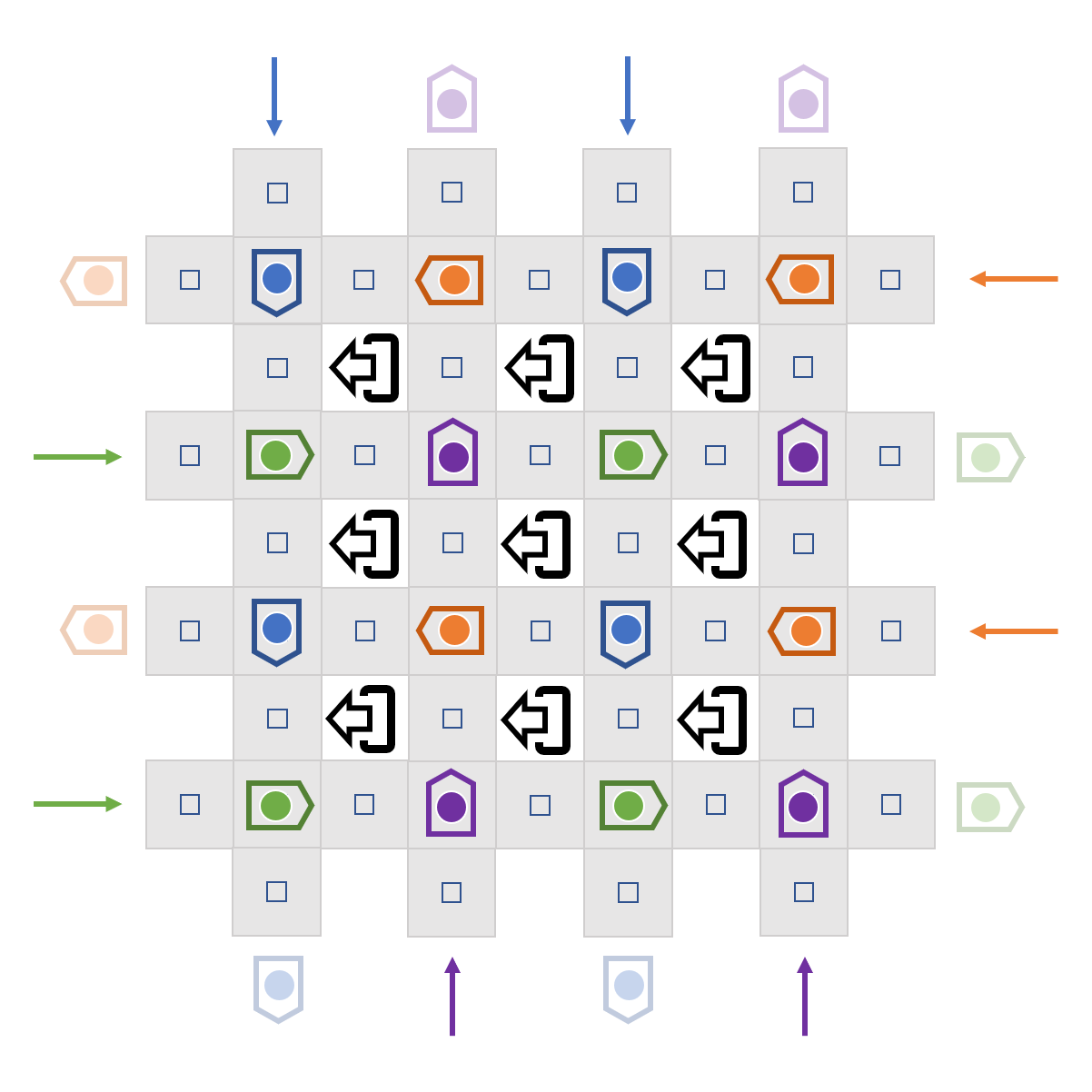}} 
	\subfigure[] {\includegraphics[width=0.4\linewidth]{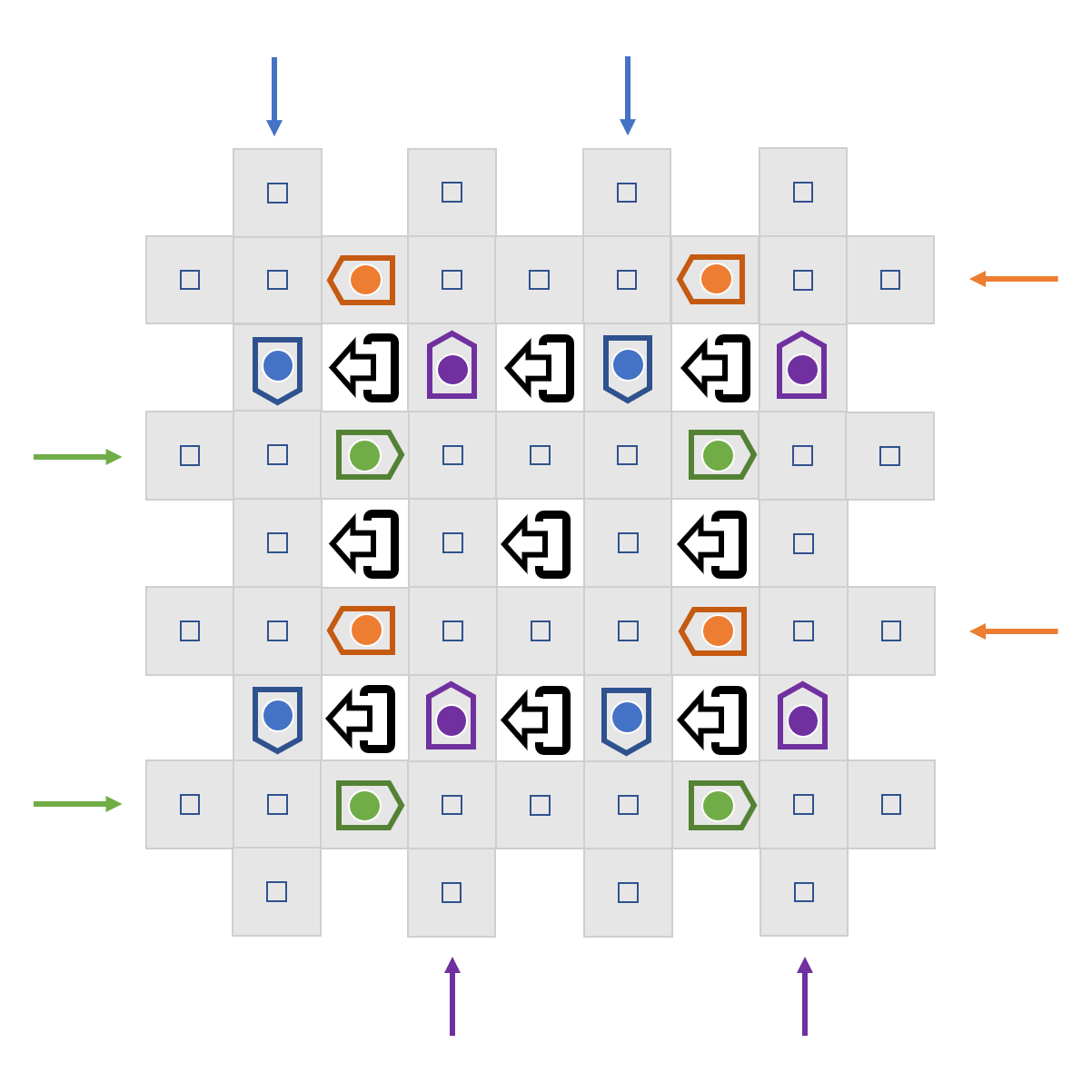}} 
	\caption{Illustration of RC-S. (a)-(d) are the four phases in one cycle.}
	\label{fig:RCS}
\end{figure}

In order to embody RC in the RSS, we regulate the behavior of robots, and propose the RC-S scheme, which can be regarded as a centralized MAPF solver. Under RC-S, the behaviors of robots are limited to three actions: moving straight, dropping off parcel, and turning. We discretize the map into grid-cells, defining each grid-cell where robots can pass through as a node. Given $(n_v,n_h)$, the network in the sorting zone is represented by a directed graph $\mathcal{G} = (\mathcal{V}=\mathcal{V}_c\cup \mathcal{V}_u\cup \mathcal{V}_{en}\cup \mathcal{V}_{ex},\mathcal{E})$, where $\mathcal{V}_c$, $\mathcal{V}_u$, $\mathcal{V}_{en}$, $\mathcal{V}_{ex}$ and $\mathcal{E}$ stand for the set of conflict nodes, unloading nodes, entrance nodes, exit nodes and edges, respectively. To reduce potential conflicts and increase the density of outlets, we adopt a one-way network, where aisles from different directions are arranged in an interleaved pattern, as shown in Figure~\ref{fig:RCS}. Based on the network structure, each VP can accommodate at most one robot. Let $\tau_e$ and  $\tau_c$ denote the fixed travel time of VPs on each grid-cell and the cycle length of rhythmic control, respectively. That is, at each entrance node, VPs enter the network at fixed intervals of $\tau_c$.

In each cycle, the distribution of VPs could be divided into four phases (Figure~\ref{fig:RCS}(a)-(d)), and the duration of each phase is $\tau_e$. Each VP moves forward one grid-cell along its assigned direction after $\tau_e$, transitioning the system from one distribution phase to the next, in the sequence: (a) $\rightarrow$ (b) $\rightarrow$ (c) $\rightarrow$ (d) $\rightarrow$ (a) $\rightarrow \dots$. We could obtain:
\begin{align}
    \tau_c = 4\tau_e \label{eq:tau_c_1}
\end{align}
$\tau_e$ and $\tau_c$ are restricted by the maximum velocity of robot $v_{max}$ and the maximum loading rate of each loading station $r_{l}$, respectively:
\begin{gather}
    \tau_e \geq \frac{D}{v_{max}} \label{eq:tau_e_1}\\
    \tau_c \geq \frac{1}{r_{l}} \label{eq:tau_c_2}
\end{gather}
where $D$ is the length of each grid-cell. At phase(a) of a cycle, new VPs are released at the entrance node of each aisle; at phase(c), VPs leave the sorting zone through the exit node of each aisle. Thus, if we do not count the VPs at the exit nodes, then the number of VPs in the network remains the same throughout the entire sorting process. Each intersection is occupied by two different VPs at phases (a) and (c) within a cycle, respectively. Each outlet could be served by four VPs in four different directions at phases(b) or (d).  The fixed number of VPs and the regularity in their phases allow us to manage and analyze this system from a global perspective. The speed of each VP could be derived by:
\begin{align}
    v_{VP} &= \frac{D}{\tau_e} \label{eq:speed_of_VP_1}
\end{align}

\begin{figure}[hbt!]
	\centering
	\includegraphics[width=0.9\textwidth]{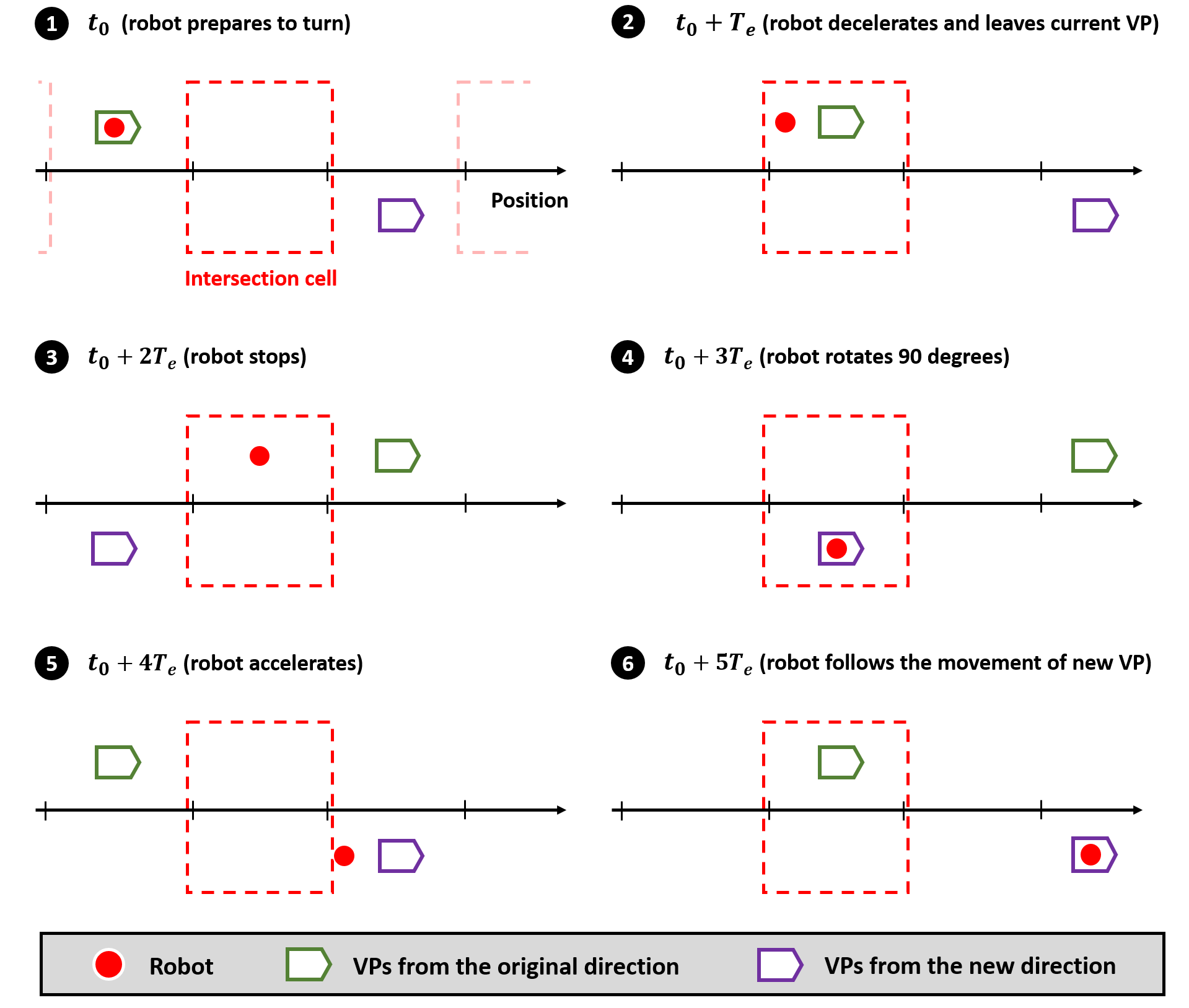}
	\caption{\label{fig:Turning} Illustration of turning and transfer process for a robot}
\end{figure}
\begin{figure}[h!]
	\centering
	\subfigure[] {\includegraphics[width=8cm]{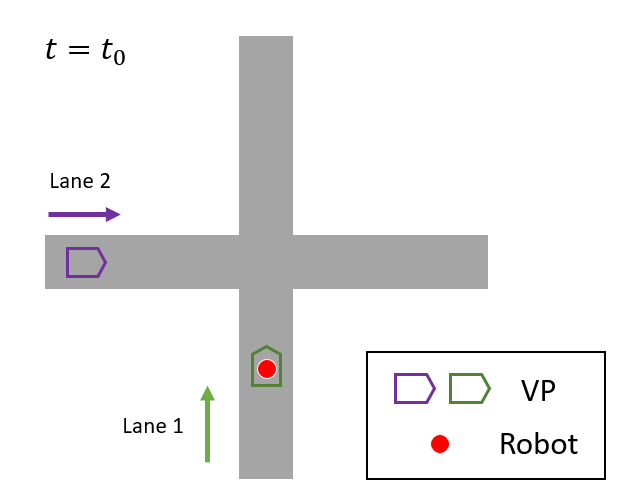}\label{fig:Turning2_1}} 
	\subfigure[] {\includegraphics[width=8cm]{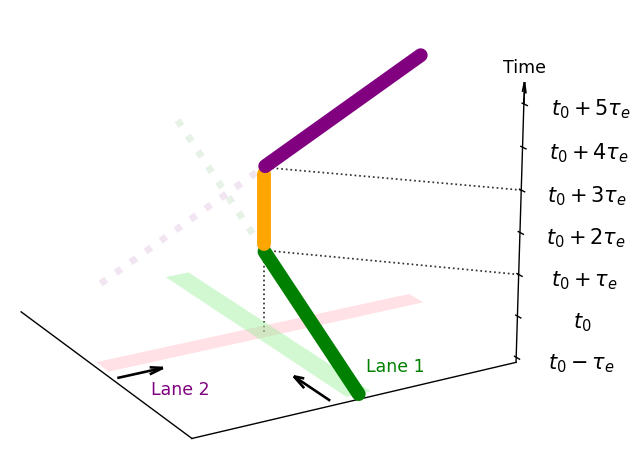}\label{fig:Turning2_2}} 
	\caption{(a) Intersection (b) Trajectory of a turning robot} 
	\label{fig:Turning2}  
\end{figure}

RSS does not require auxiliary turning lanes, since robots are able to rotate within the intersection grid-cell. The process of robot turning and transfer between VPs is demonstrated in Figure~\ref{fig:Turning}. As shown in the figure, the turning robot first decelerates to zero, stops following its original VP, rotates 90 degrees at the intersection, joins a new VP, and then follows its movement. The turning action occupies one cycle in both the original and the subsequent VPs, leading to a reduction in traffic capacity.

During this process, the parameters of RC-S must satisfy the following conditions:
\begin{gather}
   d(\tau_e,v_{VP},v_{max},c_{max}) \geq 2D  \label{eq:speed_of_VP_3}\\
   \tau_e \geq \frac{\pi}{2\omega_{r}} \label{eq:speed_of_VP_4}
\end{gather}
Constraint (\ref{eq:speed_of_VP_3}) imposes a restriction on $\tau_e$ based on the robot's maximum acceleration/deceleration rate $c_{max}$. The function $d(\tau_e, v_{VP}, v_{max}, c_{max})$ represents the maximum distance that a robot can cover during acceleration and deceleration, and its derivation is provided in Appendix B. Constraint (\ref{eq:speed_of_VP_4}) sets a minimum requirement on $\tau_e$ to ensure that the robot can complete a 90-degree turn within one phase, given its rotation speed $\omega_r$. As a result, using Expressions (\ref{eq:tau_c_1})\textendash(\ref{eq:speed_of_VP_4}) together with the kinematic parameters of the robots, we could determine the parameters for RC-S. 

The trajectory of a robot for a single delivery task starts at a loading station, passes the target outlet corresponding to the parcel's destination, and ends at another loading station.  We define a \textit{feasible path} under the RC-S scheme as follows:

\begin{definition}
A feasible path in RC-S is a sequence of VP–cycle pairs, denoted by $\{ (p_1, l_1), (p_2, l_2), \ldots, (p_n, l_n) \}$, where $(p_i, l_i)$ indicates that the robot occupies VP $p_i$ during cycle $l_i$. The sequence must satisfy the following conditions:
\begin{enumerate}
    \item \textbf{Start and end points}: The first VP $p_1$ originates from an entrance node $\nu_{en} \in \mathcal{V}_{en}$ in cycle $l_1$, and the last VP $p_n$ ends at an exit node $\nu_{ex} \in \mathcal{V}_{ex}$ in cycle $l_n$.
    
    \item \textbf{Temporal consistency}: For every consecutive pair $(p_i, l_i)$ and $(p_{i+1}, l_{i+1})$:
    \begin{itemize}
        \item If $p_{i+1} = p_i$, then $l_{i+1} = l_i + 1$ (the robot remains in the same VP).
        \item Otherwise, $l_{i+1} = l_i$ (the robot transitions between VPs).
    \end{itemize}
    
    \item \textbf{Spatial connectivity}: If $p_i \ne p_{i+1}$, then $p_i$ and $p_{i+1}$ must intersect at a common conflict node $\nu_c$ from different directions, and the corresponding entry times $t_i$ and $t_{i+1}$ must satisfy:
    \[
        \tau_c \cdot (l_i - 1) \leq t_i < t_{i+1} < \tau_c \cdot \left(l_i + \frac{1}{2}\right), \quad t_{i+1} = t_i + \frac{\tau_c}{2}
    \]
    where $\tau_c$ denotes the cycle length.
\end{enumerate}
\label{def:feasible_path}
\end{definition}

\begin{proposition}
    Given $\mathcal{G} = (\mathcal{V},\mathcal{E})$, each feasible path under the RC-S control scheme can be represented as a sequence of node–cycle pairs $\{(\nu_i, l_i)\}$.
\label{proposition:feasible_path}
\end{proposition}

The proof of Proposition~\ref{proposition:feasible_path} is provided in Appendix B.
Definition~\ref{def:feasible_path} specifies the candidate spatio-temporal trajectory patterns for robots under RC-S. Condition 1 ensures that the trajectory connects valid start and end nodes on the network boundary. Conditions 2 and 3 guarantee continuity and compliance with RC-S by restricting robot movements in accordance with the RC-S rules. Moreover, since robots are not permitted to stop except when turning, feasible paths cannot contain transfers between VPs in the same direction. Proposition~\ref{proposition:feasible_path} enables the transformation of VP reservations into node–cycle pairs, forming the basis for the optimization model presented later.

When a robot performs a sorting task, it moves from one loading station $i$ to another loading station $j$, with $i,j \in \mathcal{L}$, where $\mathcal{L}$ denotes the set of all loading stations. Let $\mathcal{O}$ represent the set of outlets. We denote by $\mathcal{R}$ the set of all feasible paths, that is, $\mathcal{R} = \bigcup_{i\in \mathcal{L}, j\in \mathcal{L}, k\in \mathcal{O}} \mathcal{R}(i,k,j)$, where $\mathcal{R}(i,k,j)$ is the subset containing feasible paths starting from loading station $i$, passing outlet $k$, and ending at loading station $j$. To simplify the notation, we use $\mathcal{R}_{i,k}$ to represent the set of all feasible paths originating from station $i$ and passing outlet $k$, and $\hat{\mathcal{R}}_{i}$ to represent the set of all feasible paths terminating at station $i$, that is, $\mathcal{R}_{i,k} = \cup_{j\in \mathcal{L}} \mathcal{R}(i,k,j)$ and $\Hat{\mathcal{R}}_{i} = \cup_{j\in \mathcal{L}, k\in \mathcal{O}} \mathcal{R}(j,k,i)$. 

The number of paths in $\mathcal{R}(i,k,j)$ is theoretically infinite, since without additional constraints, each feasible path could include an arbitrary number of turns. We limit the total number of feasible paths by allowing at most three turns and prohibiting robots from stopping within the network except for turning. These restrictions ensure that $\mathcal{R}$ is finite and that each feasible path uses no more than four distinct VPs. The rationale behind this choice will be explained in the analysis of the average travel distance in the next section.

To expand the search space, feasible paths starting at future cycles should also be considered. We further define the \textit{entry cycle}:
\begin{definition}
	The entry cycle of a specific feasible path is the earliest cycle in which all VPs in the path are unoccupied.
\label{def:entry cycle}
\end{definition}
A feasible path corresponds to a specific spatio-temporal trajectory pattern, represented by the sequence of VP reservations it occupies. The entry cycle denotes the earliest possible release cycle for the trajectory pattern, given the current VP reservation status.

Before the beginning of each cycle, the RC-S scheme performs path planning for all loaded robots that are waiting at entrance nodes and do not yet have assigned paths. These robots execute their sorting tasks based on the spatio-temporal trajectories specified by their assigned feasible paths. The corresponding feasible path assignment (FPA) problem is formulated as an integer programming model with binary decision variables. The objective is to minimize the total completion time of all unassigned tasks in the current cycle, together with penalties for any task that cannot be assigned a feasible path.

To reduce computational complexity, VP reservations from previously assigned paths are treated as fixed input. Consequently, the current planning process is based on the real-time availability of VP resources, leading to a greedy insertion approach. The main notation used in the FPA model is listed in Table~\ref{Notataion_fpa}.
\begin{longtable}{p{.20\textwidth} p{.80\textwidth}}
    \caption{Notations and explanation} \\
    \hline
    \textbf{Notations} & \textbf{Explanation} \\ \hline
        \multicolumn{2}{l}{\textit{Sets}} \\
            $\mathcal{V}$ & set of all nodes\\
		$\mathcal{L}$ & set of all loading stations\\
		$\mathcal{O}$ & set of all outlets\\
            $\mathcal{C}$ & set of all cycles starting from the current cycle\\
            $\mathcal{R}_{i,k}$ & set of feasible path patterns starting from loading station $i\in \mathcal{L}$ and passing by outlet $k \in \mathcal{O}$, independent of starting cycle\\
            $\mathcal{R}_i^{(0)}$ & Set of feasible paths available to loading station $i\in \mathcal{L}$ at the beginning of current cycle, based on the current destination outlet of the pending sorting task. If no task is pending, $\mathcal{R}_i^{(0)} = \emptyset$ \\
            $\hat{\mathcal{R}}_i$ & set of feasible paths ending at loading station $i \in \mathcal{L}$\\
        \hline
	\multicolumn{2}{l}{\textit{Parameters}} \\
            $c_{i}^r$ & completion time of feasible path $r \in \mathcal{R}_i^{(0)}$, measured in number of cycles\\
            $\hat{c}_{i}$ & penalty for delaying the assignment of a feasible path to a sorting task from loading station $i\in \mathcal{L}$ until the next cycle\\
            $\delta_{i}^{r,t,\nu,l}$ & incidence between feasible path $r \in \mathcal{R}_i^{(0)}$ starting at cycle $t \in \mathcal{C}$, node $\nu  \in \mathcal{V}$ and cycle $l \in \mathcal{C}$\\
            $N_\nu^l$ & remaining capacity in node $\nu \in \mathcal{V}$ in cycle $l \in \mathcal{C}$\\
            $d_{i}$ & sorting demand from loading station $i\in \mathcal{L}$, $d_{i} \in \{0,1\}$\\
	\hline	
	\multicolumn{2}{l}{\textit{Variables}} \\
            $x_{i}^{r,t}$ & decision on whether feasible path $r \in \mathcal{R}_i^{(0)}$ starting at cycle $t \in \mathcal{C}$ is reserved\\
            $\hat{x}_{i}$ & decision on whether the sorting task from loading station $i\in \mathcal{L}$ should be postponed to the next cycle\\
	\hline
    \label{Notataion_fpa}
\end{longtable}
\begin{align}
	\textbf{(FPA)}& &\nonumber\\
	&\underset{x_{i}^{r,t}, \hat{x}_{i}}{\min}\ \sum_{i\in \mathcal{L}} [\sum_{r\in \mathcal{R}_i^{(0)}} \sum_{t\in \mathcal{C}} (c_{i}^r + t) x_{i}^{r,t} + \hat{c}_{i} \hat{x}_{i}]&\nonumber \\
	s.t.\quad & \sum_{i\in \mathcal{L}}\sum_{r\in \mathcal{R}_i^{(0)}} \sum_{t\in \mathcal{C}} \delta_{i}^{r,t,\nu,l} x_{i}^{r,t} \leq N_\nu^l & \forall \nu \in \mathcal{V}, l \in \mathcal{C} \label{const:grid capacity} \\
	& \hat{x}_{i} + \sum_{r\in \mathcal{R}_i^{(0)}} \sum_{t\in \mathcal{C}} x_{i}^{r,t} = d_{i} & \forall i\in \mathcal{L} \label{const:maxflow} \\
        & x_{i}^{r,t},\hat{x}_{i} \in \{0,1\} & \forall i\in \mathcal{L}, r\in \mathcal{R}_i^{(0)}, t\in \mathcal{C} \label{const:inter}
\end{align}
We introduce two sets of binary decision variables. Variable $x_{i}^{r,t}$ equals one if the feasible path $r\in \mathcal{R}_i^{(0)}$ starting at cycle $t$ is assigned to the robot at loading station $i$, and zero otherwise. Variable $\hat{x}_{i}$ equals one if this robot is deferred to the next cycle, and zero otherwise. $\mathcal{R}_i^{(0)}$ denotes the set of all feasible paths available to loading station $i$ in the current cycle $t_0$. If a sorting task is pending with destination outlet $k$, then $\mathcal{R}_i^{(0)} = \mathcal{R}_{i,k}$; otherwise, $\mathcal{R}_i^{(0)} = \emptyset$. In the objective function, $c_{i}^r$ denotes the completion time of path $r$, and $\hat{c}_{i}$ denotes the penalty for delaying the assignment of a feasible path to a loaded robot until the next cycle. This penalty is introduced to prevent robots from being assigned paths with excessively long completion times. Constraints~(\ref{const:grid capacity}) enforce the spatio-temporal occupancy limits of each node $\nu$, where $N_\nu^l$ denotes the remaining capacity of node $\nu$ in cycle $l$. $N_\nu^l$ is dynamically updated based on the FPA decisions in previous cycles. Conflict nodes are allowed to accommodate up to two robots per cycle; all other nodes are restricted to one. The parameter $\delta_{i}^{r,t,\nu,l}$ equals one if feasible path $r$ starting at cycle $t$ will occupy node $\nu$ in cycle $l$, and zero otherwise. Constraints~(\ref{const:maxflow}) ensure that exactly one assignment is made for each sorting task, where $d_{i}$ equals one if there is a sorting task pending at loading station $i$. Finally, Constraints~(\ref{const:inter}) impose integrality conditions.

It is essential to solve the FPA problem efficiently within each cycle. In large-scale networks, the number of feasible paths grows rapidly, making it challenging to obtain a solution within a reasonable computation time. We propose a heuristic algorithm to generate a feasible solution in a relatively short time, as shown in Algorithm~\ref{alg3.0}. 

\begin{algorithm}
    \caption{Heuristic Algorithm for Implementing RC-S}
    \label{alg3.0}
    \small
    \begin{algorithmic}
	\renewcommand{\algorithmicrequire}{\textbf{Initialization:}}
	\REQUIRE Initialize the reservation table that records the status of VP reservation for the upcoming $N_{record}$ cycles, labeled as unoccupied. Set the waiting time of each robot that is waiting to be assigned a feasible path as 0. Set the maximum search range in each cycle as $N_{search}$. Pre-compute the set of feasible path patterns with less than three turns as $\mathcal{R} = \bigcup_{i \in \mathcal{L},\ k \in \mathcal{O}} \mathcal{R}_{i,k}$. Initialize the candidate robot list as $\mathcal{S}_r = \emptyset$. Set cycle index $l = 1$.
	\renewcommand{\algorithmicrequire}{ \textbf{Step 1.}}
	\REQUIRE At the beginning of cycle $l$, for each loading station, if the first robot in the queue has not been assigned a feasible path, add it to the candidate list $\mathcal{S}_r$.
        \renewcommand{\algorithmicrequire}{ \textbf{Step 2.}}
	\REQUIRE If $\mathcal{S}_r$ is empty, go to Step 5. Otherwise, select the robot $a$ with the longest waiting time from $\mathcal{S}_r$. Obtain its origin loading station $i$ and target outlet $k$. Search the shortest feasible path with the earliest entry cycle for robot $a$:\\
        \hspace*{1em} \textbf{for} future cycle $t \in [l, l+N_{\text{search}}]$: \\
        \hspace*{3em} \textbf{for} feasible path $r \in \mathcal{R}_{i,k}$  (sorted by path duration in ascending order): \\
        \hspace*{5em} Check if reserving path $r$ starting at cycle $t$ is feasible under current reservation table.\\
         \hspace*{5em} If feasible, record $(r, t)$ as the assigned path for robot $a$, and go to Step 3.\\
        \hspace*{3em} \textbf{end for}\\
        \hspace*{1em} \textbf{end for} 
        \renewcommand{\algorithmicrequire}{ \textbf{Step 3.}}
        \REQUIRE If a feasible pair $(r, t)$ was recorded in Step 2, then the VP-cycle occupation of this assignment is given by $\{(p_i, l_i + t)\}$. Update the reservation table with $\{(p_i, l_i + t)\}$, and assign the path to robot $a$. If no such pair is found, increase the waiting time of robot $a$ by $\tau_c$.
	\renewcommand{\algorithmicrequire}{ \textbf{Step 4.}}
	\REQUIRE Remove robot $a$ from $\mathcal{S}_r$, and return to Step 2.
	\renewcommand{\algorithmicrequire}{ \textbf{Step 5.}}
	\REQUIRE Add the newly generated VPs in the next cycle into the reservation table, labeled as unoccupied. Remove the VPs that have exited the network. 
        \renewcommand{\algorithmicrequire}{ \textbf{Step 6.}}
	\REQUIRE Set $l = l + 1$. Update the status of all robots and assign parcels to empty robots queueing in loading zone. Return to Step 1.
    \end{algorithmic}
\end{algorithm}

In Algorithm~\ref{alg3.0}, we introduce a heuristic search method that iteratively identifies the shortest feasible path with the earliest entry cycle, as a substitute for solving the original FPA model. During initialization, the reservation table is constructed with dimensions $n_h \times n_v \times N_{record}$, since the aisle network consistently contains $n_h \times n_v$ VPs at any given time. The value of $N_{record}$ is required to be greater than the maximum duration of feasible paths to fully capture the impact of each assignment on future spatio-temporal resource availability.

For each feasible path, we retain only the earliest admissible release cycle, as defined in Definition 2. By introducing a maximum search range $N_{search}$ in Step 3, we ensure that the number of pairs of feasible paths and entry cycles is finite, with an upper bound equal to the number of trajectory patterns $\lvert \mathcal{R} \rvert$ multiplied by $N_{search}$. 

In each cycle, VPs move according to the RC-S pattern, accompanied by the generation of new VPs and the departure of others. The reservation table must be updated accordingly, as shown in Step 5. In Step 6, the path-assignment results for the current cycle are executed by the robots, which then perform movement, parcel drop-off, queueing, and loading within the cycle. The algorithm then proceeds to the planning phase of the next cycle.

Compared with classic MAPF methods, the proposed approach eliminates potential path conflicts by design and greatly reduces the search range while maximizing intersection capacity. The optimal strategy of RC-S requires selecting a combination of a feasible path and an entry time for each robot, and checking the occupancy status of VPs. Algorithm~\ref{alg3.0} has a complexity of $O((n_h + n_v) \cdot N_{search} \cdot \lvert \mathcal{R} \rvert)$ in each cycle. Moreover, when the algorithm is terminated early, it can still produce a feasible solution in which robots that are not served simply wait until the next cycle. The performance of Algorithm~\ref{alg3.0} will be validated in Section~\ref{sec:validation}.
	\section{Performance Analysis of RC-S} \label{sec:analysis}
This section elaborates on the network capacity by analyzing the characteristics of RC-S and investigates how the system configurations affect system throughput. A throughput estimation formula is proposed, serving as a constraint for the planning model in Section \ref{sec:layout}.

Two important metrics in traffic analysis \citep{DAGANZO2010434} are the total vehicular distance traveled per second of operation, $\bar{m} (veh\cdot m/s)$, and the average travel distance, $\bar{l}$. The relationship between system throughput $\tilde{T}_O$ and these two metrics is shown in Equation~(\ref{eq:throughput_1}): 
\begin{align}
	& \tilde{T}_O = \frac{\bar{m}}{\bar{l}} \label{eq:throughput_1}
\end{align}
Under the mechanism of RC-S, robots are required to follow VPs for movement, so their average speed within the sorting zone is approximately equal to $v_{VP}$. Although acceleration and deceleration inevitably occur during turning, we limit the maximum number of turns in each feasible path and mitigate their impact on the average robot speed through the design of the robot turning and VP transfer process, as described in Section~\ref{sec:RC-S}. Therefore, the first metric $\bar{m}$ can be estimated as the product of the number of occupied VPs, denoted by $n^{\text{occupied}}_{VP}$, and the speed of each VP, denoted by $v_{VP}$:
\begin{align}
	& \bar{m} \approx n^{\text{occupied}}_{VP} \cdot v_{VP} \label{m_bar}
\end{align}
The speed of each VP can be obtained from Expressions (\ref{eq:speed_of_VP_1})\textendash(\ref{eq:speed_of_VP_3}), given the parameters of RC-S. For the number of occupied VPs, an intuitive assumption is that it depends on the number of workers engaged in loading parcels in the loading zone and on the spatial utilization under RC-S. We therefore introduce two factors: the workforce factor $\kappa$, representing the impact of the number of workers on the quantity of VPs that can be occupied; and the attenuation factor $\beta$, representing the utilization of VPs by the feasible paths generated by RC-S. In addition, the number of available robots $n_r$ also limits the maximum number of occupied VPs. Combining these effects, the expression for calculating $n^{\text{occupied}}_{VP}$ can be formulated as: 
\begin{align}
	& n^{\text{occupied}}_{VP}=\min\left\{\kappa\cdot \beta \cdot n_{VP},n_r\right\} \label{VP_occupied}
\end{align}
where $n_{VP}$ denotes the total number of VPs in the network, which is constant given the scale of the network. According to the phase analysis in Section~\ref{sec:RC-S details}, the RC-S scheme stipulates that the number of VPs is equal to the number of intersection nodes at any given time. Therefore, we obtain:
\begin{align}
	& n_{VP}=n_h\cdot n_v \label{number_of_VP}
\end{align}
We now investigate the forms of the two factors $\kappa$ and $\beta$.
\newline

\noindent{\large \textbf{Workforce factor}}

\begin{figure}[hbt!]
	\centering
	\includegraphics[width=0.8\textwidth]{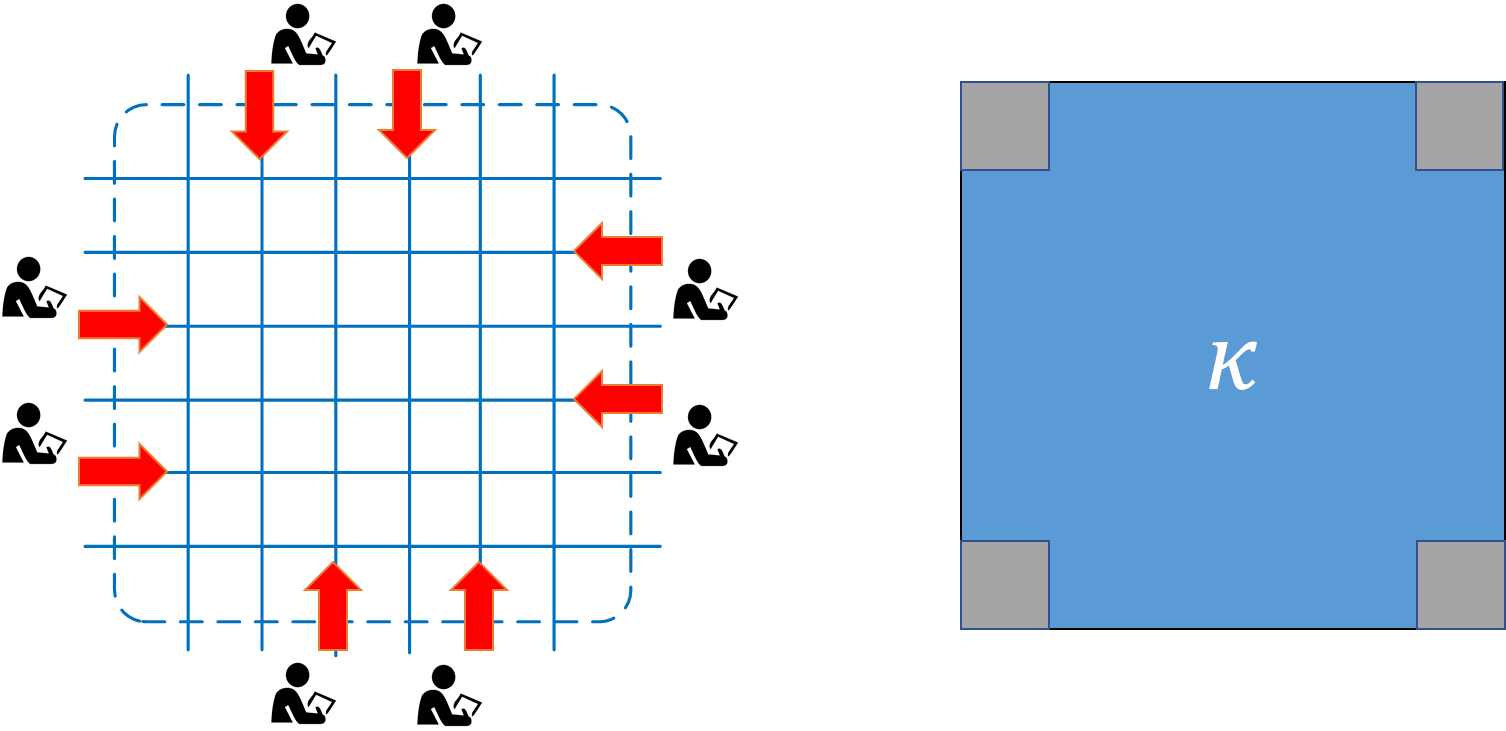}
	\caption{\label{fig:workforce} The impact of the number of workers on the VPs that can be occupied}
\end{figure}

The workforce factor $\kappa$ represents the impact of the number of workers on the number of VPs that can be occupied. It is determined by the network size $(n_h, n_v)$ and the number of workers $n_w$. For each aisle, the newly generated VPs at the entrance can be occupied by loaded robots when the corresponding loading station is active and the worker at that station is engaged in tasks. The number of such VPs is a critical determinant of the maximum throughput. Consider, for instance, a scenario where eight workers are stationed in an RSS with a $6+6$ aisle network. To ensure balanced traffic flow and mitigate congestion at the periphery of the aisle network, the workers are evenly divided into four groups and positioned along the four edges of the site, as shown in Figure~\ref{fig:workforce}. Four horizontal aisles and four vertical aisles are served by these workers. We merge the areas covered by these aisles and obtain the blue polygon.

For ease of discussion, we use $\alpha$ to denote the ratio of the number of workers to the network's maximum worker capacity.
\begin{align}
	& \alpha = \frac{n_w}{n_h+n_v} \label{eq:alpha}
\end{align}
The proportion of the blue area, referred to as the workforce factor $\kappa$, is derived by:
\begin{align}
	& \kappa = 1-{(1-\alpha)}^2 \label{eq:kappa}
\end{align}
We propose the following proposition regarding the relationship between the workforce factor and the number of occupied VPs:
\begin{proposition}
    If either (\romannumeral1) the candidate feasible paths in RC-S include no more than two turns, or (\romannumeral2) the average travel distance of a task is less than $(2-\alpha) \min\{n_v, n_h\}$, then an upper bound on the proportion of occupied VPs is $\kappa$. 
\label{proposition:workforce}
\end{proposition}
The proof is presented in Online Appendix B. Proposition~\ref{proposition:workforce} provides theoretical support for estimating the maximum network throughput. Next, we analyze the utilization rate of VPs under the RC-S scheme.
\newline

\noindent{\large \textbf{Attenuation factor}}

The attenuation factor $\beta$ represents the utilization of VPs by the feasible paths generated by RC-S. We begin by examining the path reservation strategy within the framework of RC-S. On a first-come-first-served basis, the system allocates an available robot to each newly requested sorting task, assigning it the currently shortest viable route. This robot then proceeds through the network along a predefined spatio-temporal path, ensuring collision- and queue-free delivery within the aisle network. However, this spatio-temporal reservation mechanism obligates the robot to adhere strictly to its assigned spatio-temporal trajectory, precluding any alterations such as deceleration or temporary halts. Such inflexibility can lead to underutilization of traffic capacity, as VPs that are momentarily unoccupied cannot always be sequentially linked to construct a continuous path if they do not intersect at the same point in time. To quantify the impact of this reduction, we introduce the attenuation factor $\beta$, representing the spatial utilization.

\begin{figure}[hbt!]
	\centering
	\includegraphics[width=0.8\textwidth]{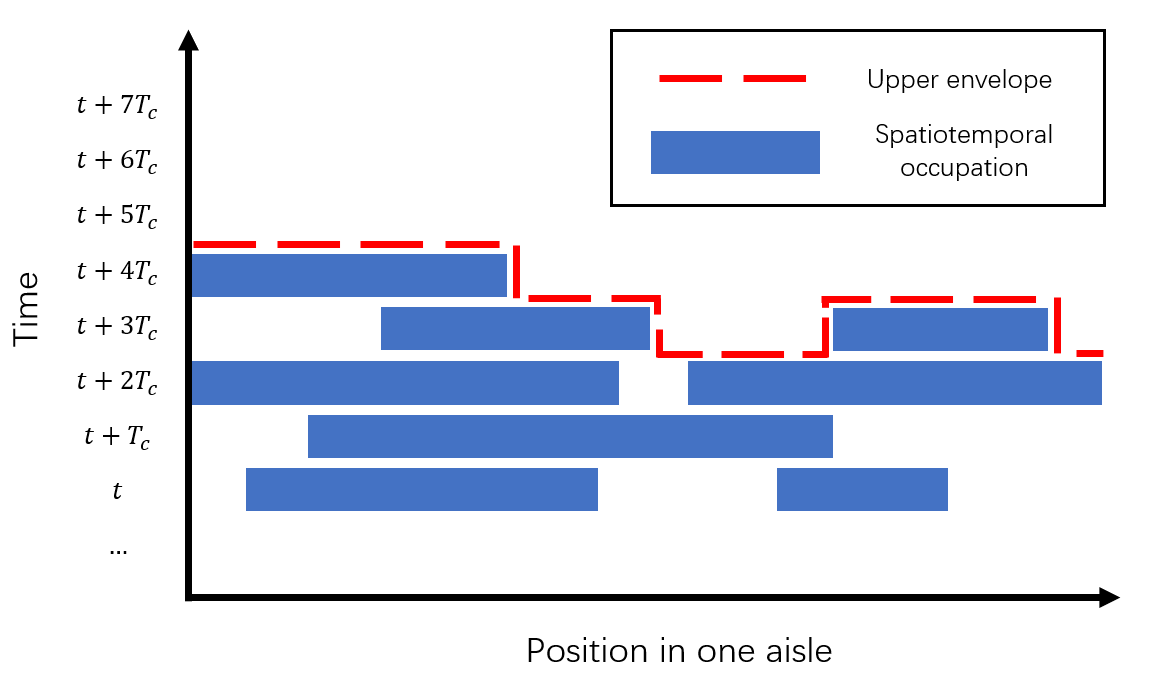}
	\caption{\label{fig:pileup} The "pile-up" process of spatio-temporal occupation under RC}
\end{figure}

\begin{figure}[hbt!]
	\centering
	\includegraphics[width=0.9\textwidth]{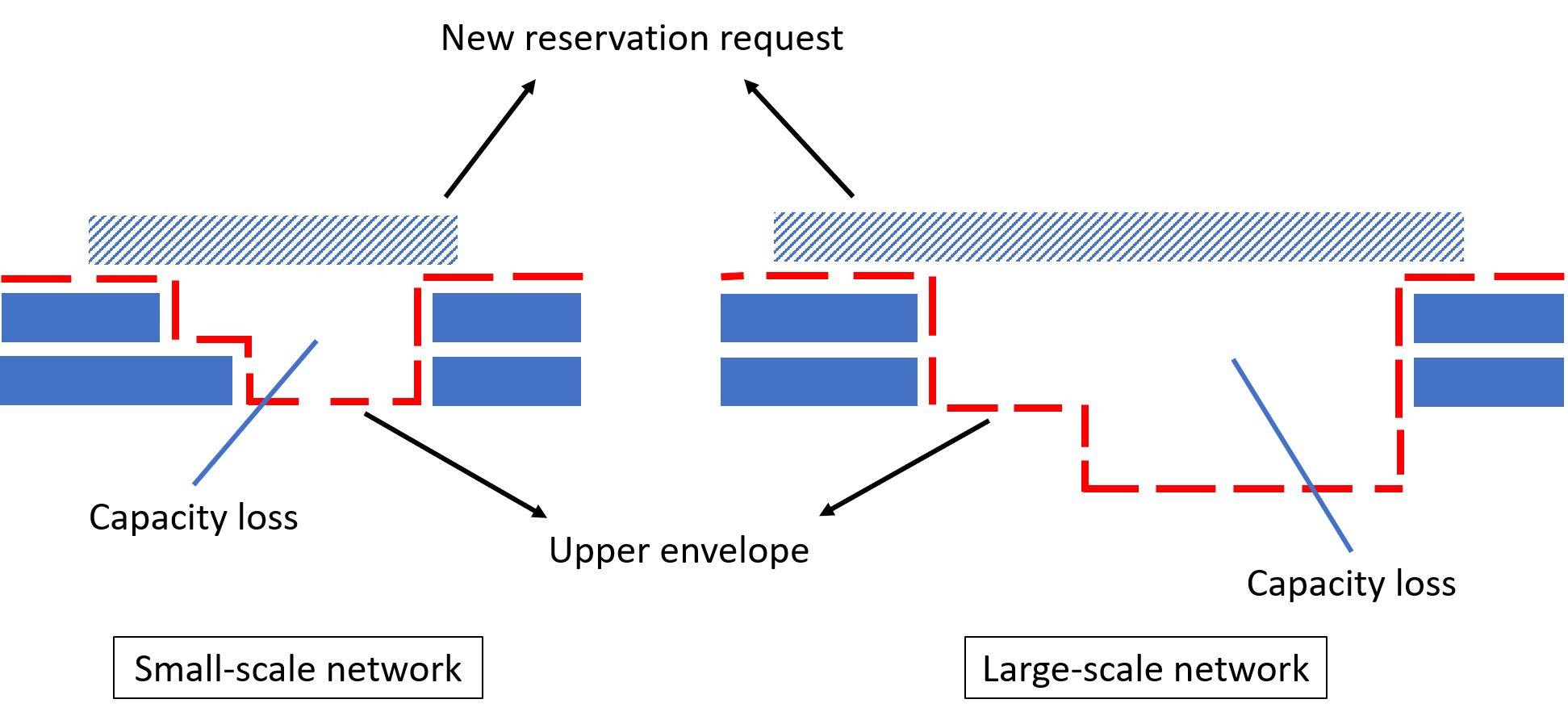}
	\caption{\label{fig:RCloss} The loss of capacity in different networks}
\end{figure}
The attenuation factor $\beta$ increases with $n_h+n_v$. Under the RC-S protocol, the spatio-temporal occupancy of an aisle exhibits a "pile-up" effect, as depicted in Figure~\ref{fig:pileup}. The figure's horizontal axis marks the spatial position within an aisle, while the vertical axis tracks discretized time intervals. Each robot's transit through this path necessitates occupying a VP, visualized as a blue bar in the illustration. The collective occupancy by different robots generates the "pile-up" phenomenon. The capacity loss at any moment is represented by the unoccupied area beneath the upper boundary of this accumulation.

When comparing networks of varying scales, i.e., those with shorter versus longer aisle lengths, we adopt a continuous approach: we assume that new sorting tasks emerge with starting points evenly spaced along the boundary and target destinations uniformly spread across the network. In stable operation, the flow distribution of paths under consistent control logic remains analogous across different network sizes, with variations only in the total path lengths. Accordingly, the pattern of unoccupied spaces beneath the upper boundary remains consistent across scales. Figure~\ref{fig:RCloss} illustrates this by comparing capacity losses across aisles of different lengths. Longer aisles tend to exhibit extended periods of spatio-temporal occupancy, potentially exacerbating capacity loss--the unoccupied intervals shown in Figure~\ref{fig:RCloss}. Consequently, we can reasonably infer that the absolute amount of lost capacity increases quadratically with the number of aisles $(n_h+n_v)$, and that the proportion of the lost traffic capacity grows linearly (because the number of available VPs increases linearly with the length of the aisles). In summary, the expression of the attenuation factor $\beta$ is as follows:

\begin{align}
	\beta(n_h,n_v) = \frac{1}{a+b\cdot (n_h+n_v)} \label{eq:beta}
\end{align}

where $a$ and $b$ are two parameters. Parameter $a$ is slightly greater than 1 and $b$ is relatively small, indicating that for a small-scale network, the attenuation factor should be close to 1 because there are fewer conflict points in a feasible path. We use simulations to obtain the actual performance of RC-S and apply linear regression to calibrate the values of the parameters in the attenuation factor. Experiments show that appropriate values are $a = 1.4$ and $b= 0.012$. Note that, due to the continuity assumption, there is a certain bias in $\beta$ when the number of aisles is small. Specifically, in smaller networks, this factor may overestimate the system throughput, especially in scenarios where the activation rate of loading stations is low. This is because some path shapes may not exist in such networks. Experimental validation of this effect will be presented in the subsequent section. To this end, we derive the final expression for the number of occupied VPs:
\begin{align}
	& n^{\text{occupied}}_{VP}(\alpha,n_r) = \min \left\{\frac{1-{(1-\alpha)}^2}{a+b(n_h+n_v)}\cdot n_h n_v,\ n_r\right\} \label{eq:N_VP_final}
\end{align}

\begin{lemma}
	$n^{\text{occupied}}_{VP}(\alpha, n_r)$ is non-decreasing and concave in both $\alpha$ and $n_r$ over its domain.
\label{lemma:1}
\end{lemma}

The proof of Lemma~\ref{lemma:1} can be established by verifying that the partial derivatives of the two terms in the minimum function are non-negative and that the Hessian matrices are both negative semi-definite. Lemma~\ref{lemma:1} indicates that opening more loading stations will increase the capacity of VPs within the network, thereby raising the upper limit of the total throughput. However, the marginal returns of constructing new loading stations will gradually diminish.
\newline

\noindent{\large \textbf{Average travel distance}}

\begin{figure}[hbt!]
	\centering
	\includegraphics[width=0.95\textwidth]{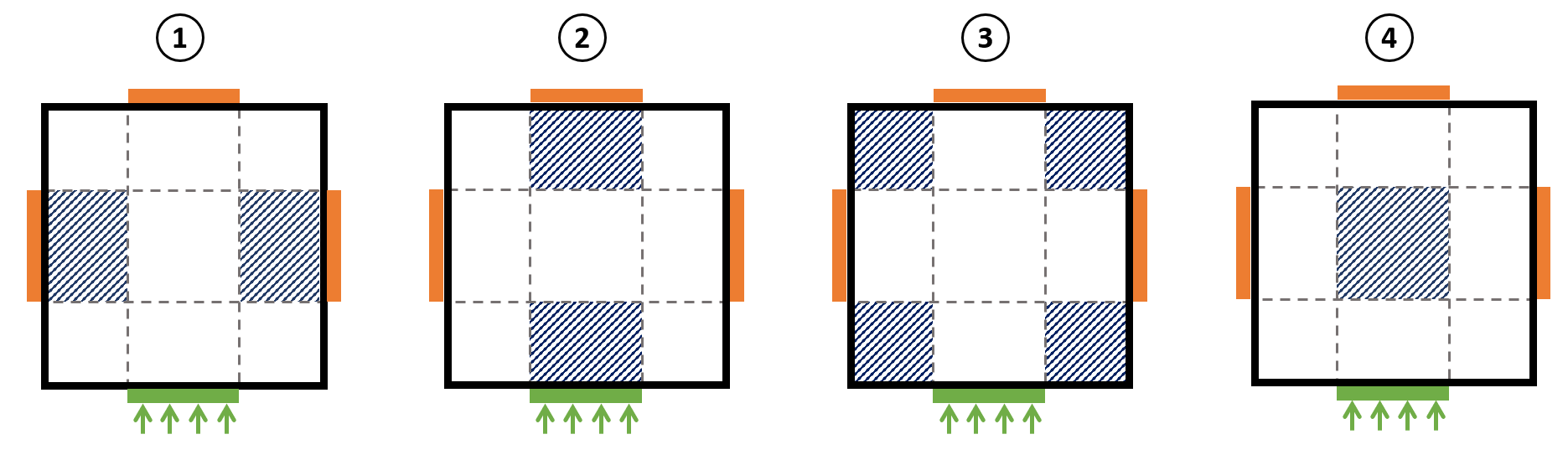}
	\caption{\label{fig:division}Division of areas according to the layout}
\end{figure}

We now derive the average travel distance $\bar{l}$ of a sorting task as the remaining component of the throughput estimation formula. Under the mechanism of RC-S, appropriate path selection follows the principles of evenly distributed sorting demand, balanced worker workload, and minimized travel distance. The workload-balancing principle ensures that loading stations do not experience starvation, thus preserving the stability of the system. It also aligns well with real-world scenarios. To maintain a balanced flow of robots between loading stations, we arrange the workers at the centers of each side, covering a length equal to $\alpha$ times the side length. Figure~\ref{fig:division} illustrates this layout, where the sorting zone is divided into four types of areas based on the spatial locations of the outlets. Sorting demand is assumed to be uniformly and continuously distributed across the entire area. As mentioned earlier, each time a robot turns, it occupies additional space in the VP fleet, so it is necessary to restrict the maximum permissible number of turns in one delivery. However, from the diagram, it can be observed that to ensure the outlets in Area~3 are reachable from each loading station, paths with at least three turns are required. According to our simulation experiments, introducing paths with four or more turns does not effectively improve the sorting throughput. Therefore, we prioritize paths with fewer turns and do not consider paths that require more than three turns. We denote the path length serving Area~$i$ by $l_i$, with the proportion of this area represented by $p_i$. 
We next present the rationale for their calculation and the corresponding formulas. For a detailed derivation, please refer to Appendix~B.

\begin{figure}[hbt!]
	\centering
	\includegraphics[width=0.95\textwidth]{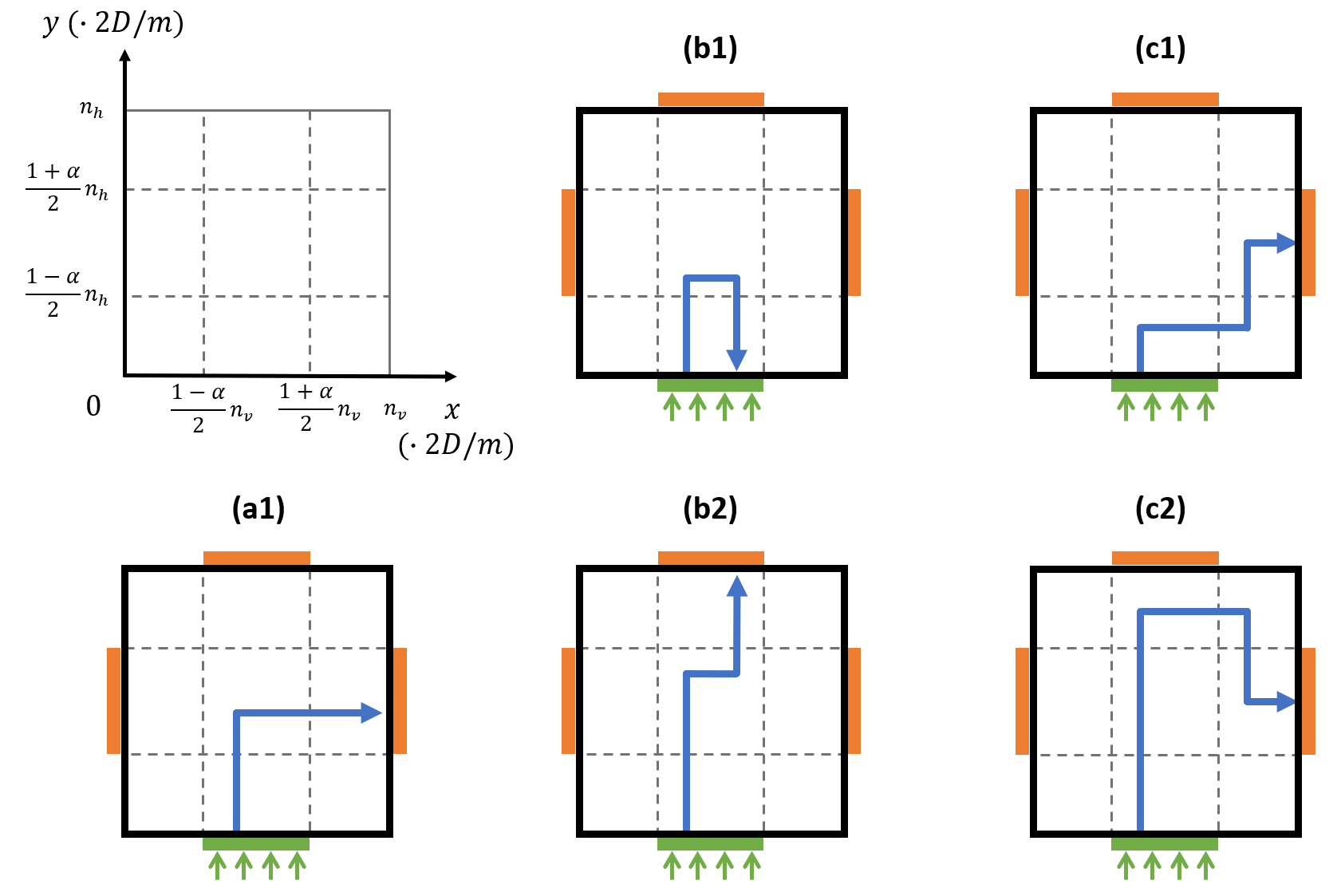}
	\caption{\label{fig:different_path}Service paths for different categories}
\end{figure}

First, for Area~1, servicing the outlets requires one change of direction, as shown in Figure~\ref{fig:different_path}(a1). Assuming balanced workload, the turning points should have equal probability of occurring within the central area. We then perform a weighted accumulation based on the length and take the average. Consequently, we obtain the following expressions for the area proportion and the expected path length:
\begin{align}
    p_1 &= \alpha(1-\alpha) \label{eq:p1} \\
    E[l_1] &= 2D\cdot \left[\frac{\frac{9+\alpha^2}{6}(n_h^2+n_v^2)-n_h n_v-\frac{1}{3}}{n_h+n_v}+1\right] \label{eq:l1}
\end{align}
For paths that serve Area~2, their endpoints are located at either the lower or upper loading stations, as shown in Figure~\ref{fig:different_path}(b1) and (b2), respectively. Let $(x_1,y_1)$ and $(x_2, y_2)$ denote the coordinates of the two turning points in a path. $x_1$ and $x_2$ are uniformly distributed over the length covered by the workers, while $y_1$ and $y_2$ are identical and uniformly distributed across the area. Notably, when $x_1 = x_2$, the path actually does not involve any turns. We assume that the robot chooses the nearest top or bottom exit point when making the second turn. The expressions for $p_2$ and $E[l_2]$ are as follows:
\begin{align}
    p_2 &= \alpha(1-\alpha) \label{eq:p2} \\
    E[l_2] &= 2D\cdot \left[\frac{\alpha(n_h^2+n_v^2)}{3(n_h+n_v)}-\frac{2}{3\alpha(n_h+n_v)}+\frac{3n_h n_v}{2(n_h + n_v)}\right] \label{eq:l2}
\end{align}
For Area~3, as previously mentioned, only paths involving three turns can serve it. Area~3 includes four individual small areas. The two small areas at the bottom, whose service paths are shown in Figure~\ref{fig:different_path}(c1), have a length equal to that of the paths with one turn in Figure~\ref{fig:different_path}(a1); the two areas at the top, with service paths depicted in Figure~\ref{fig:different_path}(c2), involve two extra segments of path in the $y$-direction (or in the $x$-direction when the starting point is on the left or right side). As a result, the expressions for $p_3$ and $E[l_3]$ are as follows:
\begin{align}
    p_3 &= {(1-\alpha)}^2 \label{eq:p3} \\
    E[l_3] &= 2D\cdot \left(\frac{n_h+n_v}{2} + \frac{(1+\alpha)}{4}\cdot \frac{n_h n_v}{n_h + n_v}\right) \label{eq:l3}
\end{align}
Finally, for Area~4, which is in the center, we assume that it is equally served by the paths of Areas~1 and 2.
\begin{align}
    p_4 &= {\alpha}^2 \label{eq:p4} \\
    E[l_4] &= \frac{1}{2}E[l_1] + \frac{1}{2}E[l_2] \label{eq:l4}
\end{align}

In summary, the average travel distance of a sorting task can be obtained by the following weighted summation:
\begin{align}
	\bar{l}(n_h, n_v, \alpha) = \sum_{i=1}^4 p_iE[l_i] \label{eq:travel distance}
\end{align}

We now further investigate the properties of Equation~(\ref{eq:travel distance}).
\begin{lemma}
	For uniformly distributed sorting demands, if the length-to-width ratio is less than 2 and $\min\{n_h,n_v\} \geq 4$, then (\romannumeral1) as $\alpha$ increases, $\bar{l}(n_h, n_v, \alpha)$ first decreases and then increases;\quad (\romannumeral2) when the site size is fixed ($n_h\cdot n_v = \text{constant}$), a square-shaped site has a smaller average travel distance;\quad (\romannumeral3) when the site is square-shaped, the range of average travel distance is given by $2D\cdot n_h < \bar{l}(n_h, n_v, \alpha) < 2D\cdot\frac{9}{8}n_h$.
\label{lemma:2}
\end{lemma} 

The proof of Lemma~\ref{lemma:2} is provided in Appendix~B. The conditions in the lemma indicate that its conclusions is applicable to an RSS where the length-to-width ratio does not exceed two. This is a reasonable setup in practice, as excessively elongated sites accommodate fewer outlets for the same area. Lemma~\ref{lemma:2}(\romannumeral1) states that the sorting efficiency of robots exhibits a trend of first increasing and then decreasing with the increase as the number of loading stations increases. This is because, under the worker load-balancing criterion, the average distance from the corner loading stations to each outlet is greater than that from loading stations located in the middle of the four sides. Lemma~\ref{lemma:2}(\romannumeral2) suggests that, under the permissible conditions, designing the aisle network in a square shape can reduce the demand for robots. As mentioned in Section~\ref{sec:description}, the outlet density remains constant, thus the distribution of outlets in RSS does not need to be altered. Lastly, Lemma~\ref{lemma:2}(\romannumeral3) provides the range of the average distance to be covered in a single delivery within RC-S. This serves as a foundation for the efficiency analysis in the following sections.

Combining Equations~(\ref{eq:alpha})\textendash(\ref{eq:travel distance}), the estimation formula for the sorting throughput under a given configuration $(n_h,n_v,n_w,n_{r})$ is as follows:
\begin{align}
	&\tilde{T}_{O}(n_h,n_v,n_w,n_{r}) = \frac{D\cdot n^{\text{occupied}}_{VP}(\frac{n_w}{n_h+n_v},n_r)}{\tau_e \cdot \bar{l}(n_h, n_v, \frac{n_w}{n_h+n_v})} \label{eq:T_O_final}
\end{align}

If we fix the values of $n_h$ and $n_v$, the upper bound of the system throughput can be obtained by setting $n_r$ to infinity and assigning workers to all loading stations: 
\begin{align}
	&\tilde{T}_{O}(n_h,n_v,n_w,n_{r}) \leq \frac{D\cdot n^{\text{occupied}}_{VP}(\frac{n_{l}}{n_h+n_v},\infty)}{\tau_e \cdot \bar{l}(n_h, n_v, \frac{n_{l}}{n_h+n_v})} \label{ineq:T_O_upper_bound}
\end{align}

We use $\tilde{T}_{M}(n_h,n_v,n_{l})$ to denote this upper bound on the right-hand side of Inequality~(\ref{ineq:T_O_upper_bound}). We now investigate the properties of this bound.
\begin{proposition}
	$\frac{\partial \tilde{T}_{M}(\cdot)}{\partial n_{l}}$ has a unique zero point within the valid range of $n_{l}$. Denote this value by $n_{l}^{(0)}$, then $\tilde{T}_{M}(\cdot)$ is monotonically increasing on the interval $(0, n_{l}^{(0)}]$ and monotonically decreasing in the interval $(n_{l}^{(0)}, n_h+n_v]$.
\label{proposition:TM}
\end{proposition}

The proof of Proposition~\ref{proposition:TM} can be established by combining Lemma~\ref{lemma:1} and Lemma~\ref{lemma:2}, and is presented in detail in Appendix~B. Proposition~\ref{proposition:TM} reveals that under the RC-S robot control scheme, the maximum throughput has a critical point, and that this critical point is not attained when all loading stations are activated, due to workload balancing and uniformly distributed demand. When the proportion of active loading stations approaches 1, the growth rate of $n^{\text{occupied}}_{VP}$ becomes slower than the growth rate of the service path length, resulting in a decrease in throughput. We will verify this formula through experiments in the next section.
        \section{Numerical Validation} \label{sec:validation}

\subsection{Comparative analysis: RC-S versus benchmark traffic management framework} \label{sec:vali1}
In this section, we compare our proposed RC-S with cooperative A* (CA*) \citep{Silver_2021} for traffic management in a warehouse, in order to illustrate the superiority of RC-S in terms of robot travel distance and system throughput. CA* is a framework based on a simple prioritized-planning scheme: each agent is first assigned a unique priority, and, in order of these priorities, the algorithm finds the shortest path for each agent that avoids conflicts with agents of higher priority. We use SIPP \citep{phillips2011sipp}, an efficient variant of location-time A*, as the lower-level solver of CA*. Moreover, other well-known MAPF frameworks, such as PBS \citep{ma2019searching} and RHCR \citep{Li_Tinka2021}, are found to be less efficient than CA* in large-scale settings ($n_r > 100$) of our experiments, primarily due to the increased likelihood of cycle conflicts within the RSS framework, and are thus excluded from this comparative analysis.

\begin{table}[!htbp]
	\centering
	\caption{RC-S parameters used in the simulations}
	\begin{tabular}{cccc}
		\toprule
		$D$ (m) & $\tau_e$ (s) & $\tau_c$ (s) & $v_{VP}$ (m/s) \\
		\midrule
		1      & 0.5 & 2 & 2\\
		\bottomrule
	\end{tabular}
	\label{tab:table51}
\end{table}
We compare the performance in two scenarios: 12+12 aisles and 20+20 aisles. To ensure uniformity in the experimental setup, we use directed maps in all tests to prevent swapping conflicts. All potential loading-station locations are activated in the experiments, and parcel destinations are evenly distributed among outlets. The RC-S parameters, detailed in Table~\ref{tab:table51}, are based on sorting robot data from the Geekplus company. For a fair comparison of algorithm performance, robots move at a constant speed of 2~m/s, matching the speed of VPs in RC-S. The simulations do not consider acceleration and deceleration phases. We conduct 50 independent experiments for each scenario. In each experiment, we set a warm-up period of 30 minutes and then record the system performance for 60 minutes of continuous operation. Our algorithms and simulations are implemented in Python~3.11, and all experiments are conducted on a personal computer running Windows 11 with an Intel i7-9700F CPU and 16~GB of RAM. All processes are run single-threaded.

\begin{table}[!htbp]
	\centering
	\caption{Results on 12+12 aisle network}
	\begin{tabular}{cccp{0.2cm}ccp{0.2cm}cc}
		\toprule
		$n_r$ & \multicolumn{2}{l}{Throughput ($\times10^3$/h)} &&  \multicolumn{2}{l}{Service time (s)} &&  \multicolumn{2}{l}{Run-time per cycle(ms)} \\ \cline{2-3} \cline{5-6} \cline{8-9}
		& RC-S & CA* && RC-S & CA* && RC-S & CA* \\
		\midrule
        40&8.83&4.72&&11.73&14.17&&5.67&5.99\\
        80&15.12&9.03&&12.92&14.57&&6.80&16.85\\
        120&17.45&12.40&&13.89&15.24&&8.50&35.37\\
        160&18.06&14.32&&14.26&15.79&&9.31&55.31\\
        200&18.33&14.94&&14.51&16.03&&9.63&63.36\\
		\bottomrule
	\end{tabular}
	\label{tab:table61}
\end{table} 
\begin{table}[!htbp] 
	\centering
	\caption{Results on 20+20 aisle network}
	\begin{tabular}{cccp{0.2cm}ccp{0.2cm}cc}
		\toprule
		$n_r$ & \multicolumn{2}{l}{Throughput ($\times10^3$/h)} &&  \multicolumn{2}{l}{Service time (s)} &&  \multicolumn{2}{l}{Run-time per cycle (ms)} \\ \cline{2-3} \cline{5-6} \cline{8-9}
		& RC-S & CA* && RC-S & CA* && RC-S & CA* \\
		\midrule
        50&7.68&3.94&&18.85&21.75&&15.53&11.61\\
        100&14.91&7.79&&19.33&21.95&&15.79&29.94\\
        200&26.82&15.07&&20.77&22.54&&19.76&99.61\\
        300&32.32&22.58&&21.48&23.33&&29.27&226.76\\
        400&33.66&24.37&&23.41&23.92&&37.40&363.45\\
		\bottomrule
	\end{tabular}
	\label{tab:table62}
\end{table}
Tables~\ref{tab:table61} and \ref{tab:table62} report the system throughput, average service time, and average run-time. The comparative results indicate that the run-time of RC-S is consistently lower than that of CA*, and in large-scale cases, it can be less than one tenth of that of CA*. Moreover, RC-S always assigns shorter service paths to robots, resulting in an average service time reduction of 10.3\% compared to CA*. In most experiments, RC-S demonstrates higher throughput than CA*, with the only exception observed in the 20+20\textendash aisle scenario with $n_r=400$. This observation underscores the suitability of the RC-S algorithm for RSS, which is characterized by a compact network and numerous outlets acting as obstacles. The difference in algorithm performance is primarily due to the following factors: (1) RC-S optimizes by selecting from candidate spatio-temporal paths, whereas CA* and similar cell-based search algorithms consider variable waiting times at each cell and allow robots to return to previously visited cells; and (2) RC-S imposes limits on both the total number and the maximum density of robots and remains conflict-free, whereas the high robot density around outlets in the benchmark algorithm leads to significant resource expenditure on conflict resolution.

\subsection{Throughput estimation formula validation} \label{sec:vali2}
This section validates the throughput estimation formula proposed in Section~\ref{sec:analysis} for different active ratios of loading stations. We consider five scenarios with different network sizes: $n_h = n_v = 12,14,16,18,20$, respectively. The number of robots is set to be equal to the number of activated VPs in the aisle network, plus five robots per loading station to ensure stable operation. Other experimental parameters are specified in detail in Section~\ref{sec:vali1}. 

\begin{figure}[hbt!]
	\centering
	\subfigure[] {\includegraphics[width=0.42\linewidth]{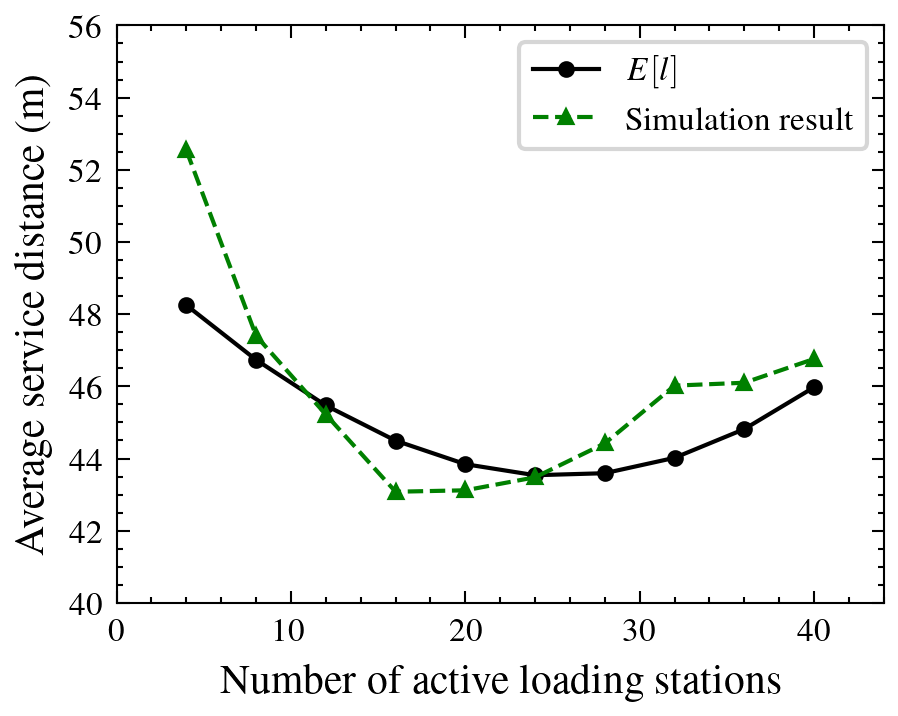}\label{fig:validation1}} 
	\quad
	\subfigure[] {\includegraphics[width=0.42\linewidth]{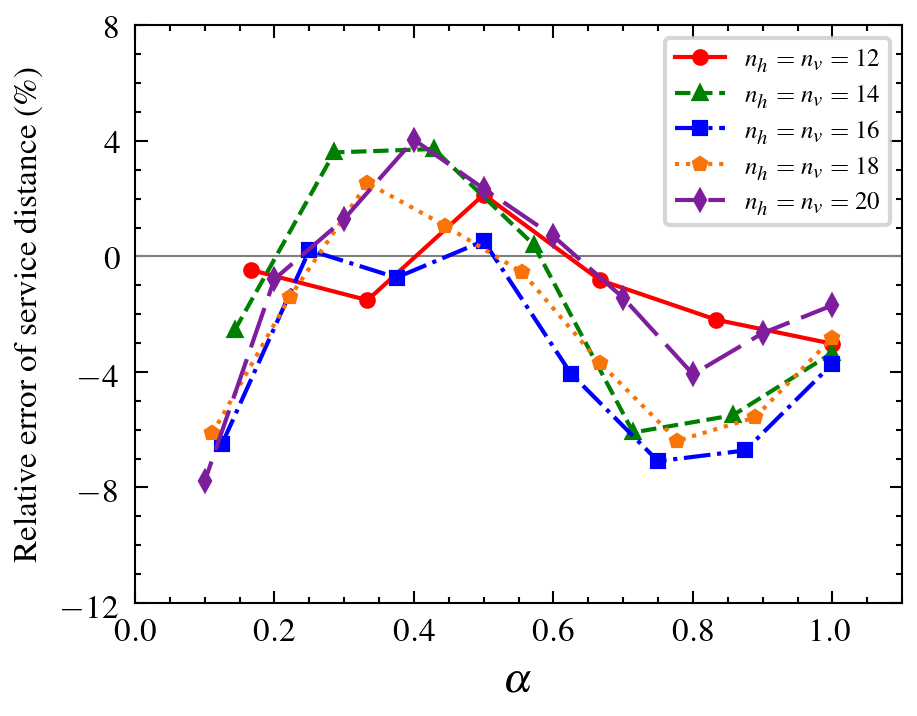}\label{fig:validation2}} 
	\qquad
	\subfigure[] {\includegraphics[width=0.42\linewidth]{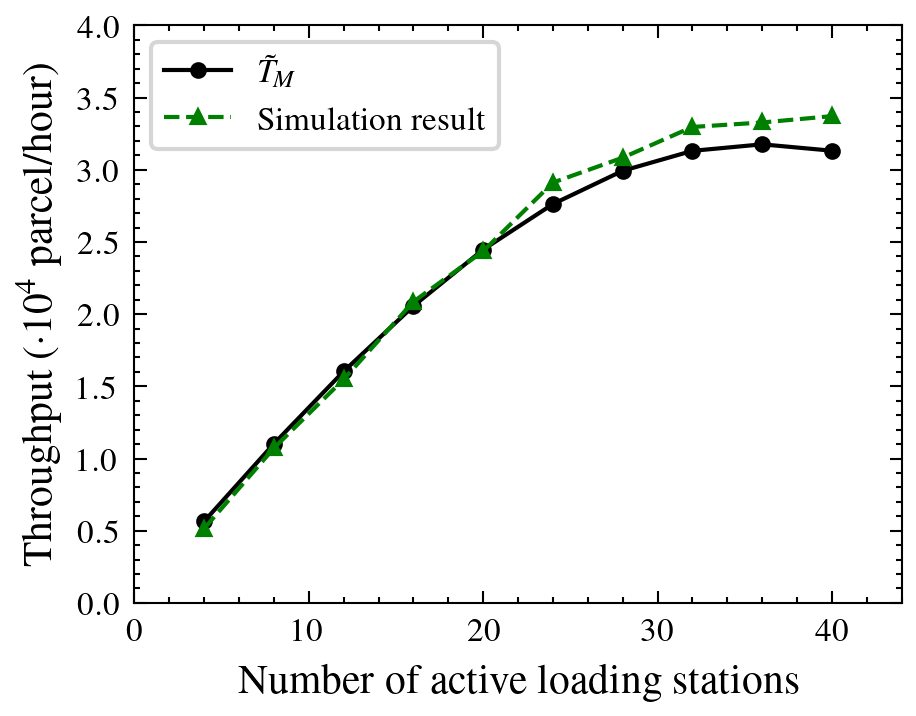}\label{fig:validation3}}
	\quad
	\subfigure[] {\includegraphics[width=0.42\linewidth]{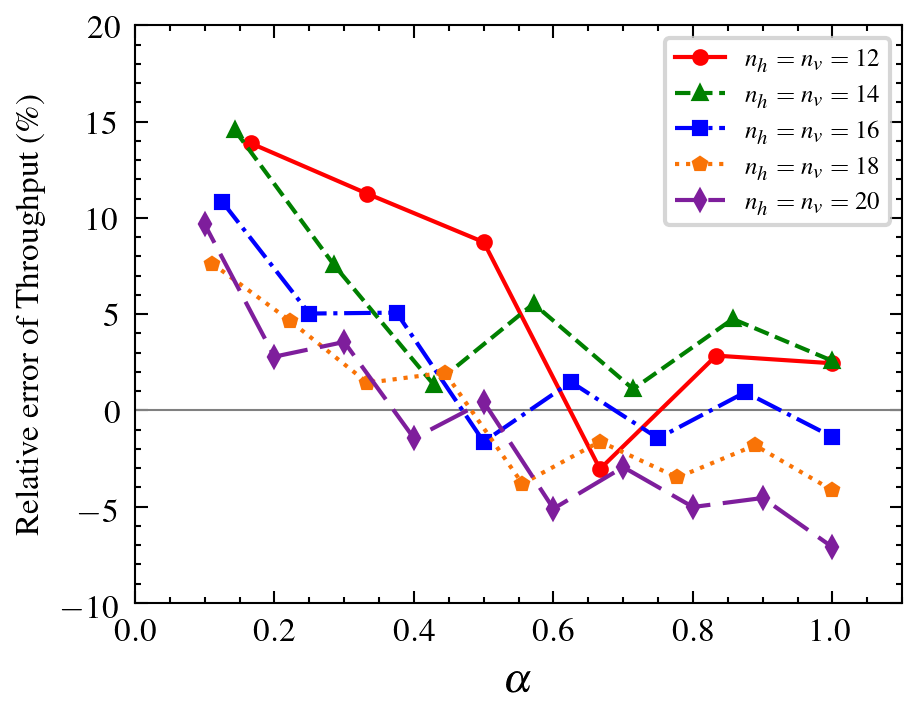}\label{fig:validation4}} 
	\caption{Analytical formula compared with simulation results.}
	\label{fig:validation}
\end{figure}

\begin{figure}[hbt!]
	\centering
	\includegraphics[width=0.54\textwidth]{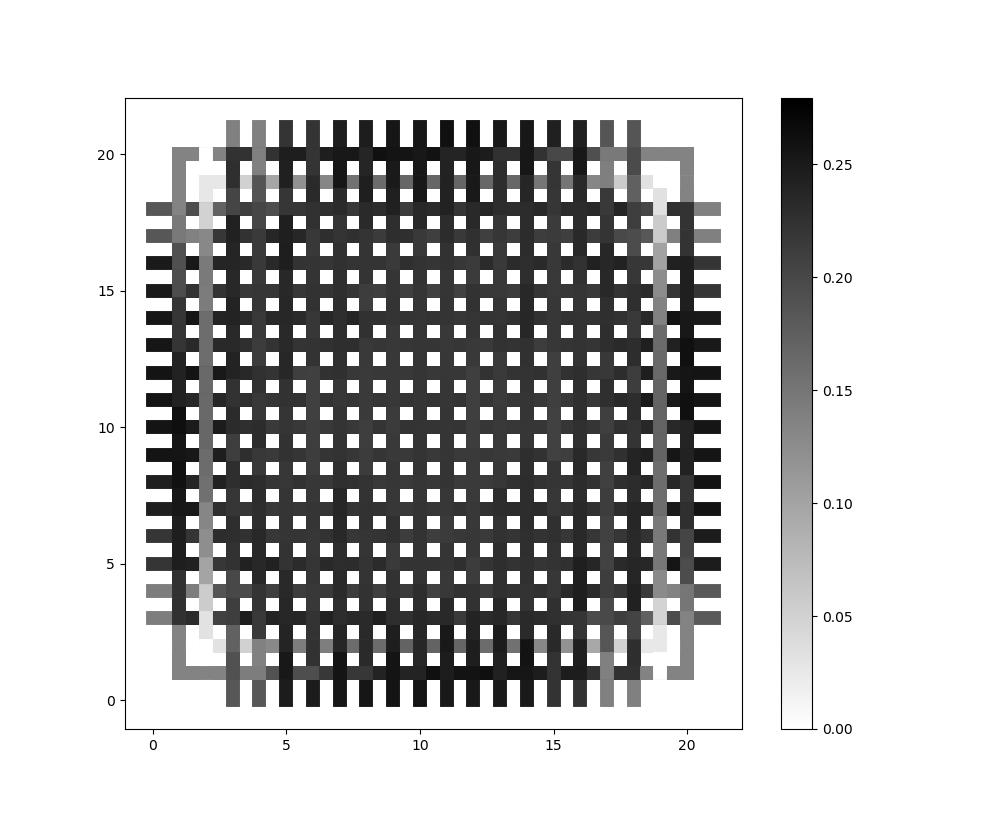}
	\caption{\label{fig:load} Robot flow distribution in simulation, $n_h=n_v=20, n_w=40$}
\end{figure}

Figure~\ref{fig:validation1} shows the trend of the average service distance as a function of the number of active loading stations when $n_h=n_v=20$. It indicates that the formula and the simulation results exhibit a similar trend of first decreasing and then increasing. This suggests that the service frequency of each robot is highest when about half of the loading stations are activated. Figure~\ref{fig:validation2} shows the relative error between the formula and the simulation results under five scenarios. It can be seen that the error is always within $\pm 8\%$, indicating a good fit.

Figure~\ref{fig:validation3} shows the trend of the system throughput as a function of the number of activated loading stations when $n_h=n_v=20$. The trend demonstrates a clear marginal effect of the number of activated loading stations: when the system is operating near its maximum capacity, there is a certain degree of resource waste. Figure~\ref{fig:validation4} shows the relative error between the formula and the simulation results for five cases, which remains within $\pm 15\%$. Figure~\ref{fig:load} depicts the average traffic flow distribution during stable operation of the system in simulations with the maximum number of workers. The darker the color of the aisle segments, the higher the traffic load. We observe that the corner loading stations have almost no incoming or outgoing traffic flow. This indicates that Algorithm~\ref{alg3.0}, in practical execution, cannot effectively balance the distribution of robots over loading stations in this extreme situation.

\subsection{Underestimation of efficiency in queueing models} \label{sec:vali3}
Queueing models have been widely used in the analysis of flexible machining systems, where flexibility is defined as the ability to arbitrarily specify the processing sequence of different workpieces \citep{buzacott1993stochastic}. In robotic systems, robots function similarly to pallets in traditional manufacturing systems, transferring workpieces from one machine to another for processing. However, we find that queueing models could not accurately capture robot efficiency when the number of robots has not reached the network's capacity limit. We apply the CQN model established in \cite{zou2021robotic}, in which a two-tier RSS with a closest Robot-to-Loading-Station assignment rule is modified to meet the assumptions in Section~\ref{sec:description}. Eight scenarios of large-scale networks with few robots are tested, using the same experimental setup as in Section~\ref{sec:vali2}. The results are shown in Table~\ref{tab:table60}, where Eq.~(\ref{eq:T_O_final}) represents the proposed estimation formula. The CQN model significantly underestimates system throughput in scenarios with few robots. It lacks the capability to accurately depict traffic capacity under low flow conditions with RC-S control and overestimates the impact of conflicts under low-demand conditions, thereby underestimating individual robot efficiency. Although the CQN model provides accurate estimates when the network approaches capacity, it leads to unnecessary costs due to surplus robots during off-peak periods.
\begin{table}[!htbp]
	\centering
	\caption{Error in throughput estimation under RC-S control}
	\begin{tabular}{cc cccc p{0.1cm} cccc}
		\toprule
		& $n_h + n_v$ & \multicolumn{4}{c}{$20+20$} &&  \multicolumn{4}{c}{$24+24$} \\ \cline{3-6} \cline{8-11}
	& $n_r$ & 50 & 100 & 200 & 300 && 50 & 100 & 200 & 300 \\
		\midrule
    \multirowcell{3}{Throughput\\($\times 10^3$/h)}&Simulation & 8.43	& 16.25	& 28.60	& 32.90 && 7.27 & 14.30 & 27.01 & 36.49
\\
    &CQN & 5.98 & 11.87 & 22.98 & 31.47 && 4.96 & 9.88 & 19.52 & 28.61\\
    &Eq.(\ref{eq:T_O_final}) & 7.83 & 15.66 & 31.32 & 33.86 && 6.55 & 13.10 & 26.19 & 39.08\\
    \midrule
    \multirowcell{2}{Relative error\\(\%)} & CQN & 29.03 & 26.95 & 19.65 & 4.34 && 31.82 & 30.92 & 27.73 & 21.61\\
    &Eq.(\ref{eq:T_O_final}) & 7.12 & 3.65 & -9.52 & -2.92 && 9.92 & 8.41 & 3.04 & -7.10\\ 
		\bottomrule
	\end{tabular}
	\label{tab:table60}
\end{table} 

We also conduct an experiment on a real case of Deppon Express described in \cite{zou2021robotic}, where the network scale is $n_h+n_v=18+6$, with 108 outlets, 6 loading stations, and 170 robots. The aisle width $D$, maximum robot speed $v_{max}$, and worker loading rate $r_l$ are set according to the operational data provided by Deppon Express. We compare the performance of RC-S with the results from \cite{zou2021robotic}, as shown in Table~\ref{tab:table202}, where $T_{D}$, $T_{CQN}$, and $T_{RCS}$ represent the throughput of the Deppon system, the CQN model and the simulation under the RC-S scheme, respectively. Together with the experimental results in Section~\ref{sec:vali1}, this real-case study shows that the RC-S scheme significantly improves the operational efficiency of robots, thereby enhancing the overall system throughput. Meanwhile, the CQN model effectively captures congestion in the RSS system implemented by Deppon Express but underestimates the traffic capacity of the network when the RC-S scheme is applied.

\begin{table}[!htbp]
	\centering
	\caption{Real case validation}
	\begin{tabular}{cccccccc}
		\toprule
		$D$ (m) & $r_w$ (/s) & $v_{max}$ (m/s) & $\tau_c$ (s) & $v_{VP}$ (m/s) & $T_{D}$ (/h) & $T_{CQN}$ (/h) & $T_{RCS}$ (/h) \\
		\midrule
		0.6 & 0.42 & 2 & 1.2 & 2 & 7,163 & 6,907 & 8,758\\
		\bottomrule
	\end{tabular}
	\label{tab:table202}
\end{table}

\subsection{Evaluating RC-S across diverse scenarios} \label{sec:vali4}
In real sorting centers, due to cost constraints and demand arrival rates, loading stations are often not distributed across the entire perimeter of the RSS, resulting in only partial availability of network entrances and exits. The placement of outlets is also restricted by site-specific factors, such as load-bearing columns or the presence of certain equipment. Additionally, demand arrival rates are heterogeneous and are not uniformly distributed across outlets. In this section, we evaluate the performance of RC-S across various application scenarios and test the impact of different system configurations on efficiency by simulation. For each set of experiments, we evaluate overall performance across 50 independent runs. Each run includes a 30-minute warm-up period, followed by an observation period during which average performance is recorded over one hour.

\begin{figure}[ht!]
	\centering
	\includegraphics[width = 10cm]{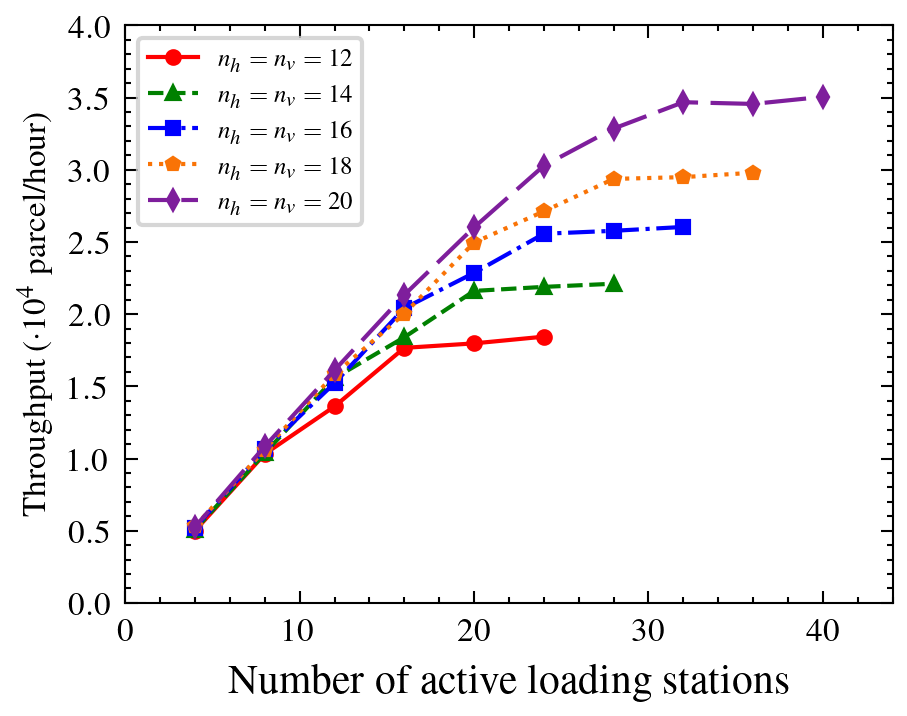}
	\caption{Effect of different numbers of active loading stations on throughput}
	\label{fig:num_of_station}
\end{figure}

We first compare the impact of different numbers of active loading stations on average throughput. We consider five scenarios with varying network sizes: $n_h = n_v = 12,14,16,18,20$, respectively. The number of robots is set to be sufficient in each scenario, ensuring that it does not become a bottleneck for the sorting process. Other RC-S parameters are specified in detail in Table~\ref{tab:table51}. Figure~\ref{fig:num_of_station} shows the trend of throughput across the five scenarios. It can be observed that large-scale networks not only have a higher upper bound on throughput but also exhibit steeper slopes, representing higher average efficiency per loading station. This is because, under RC-S, larger networks provide more space to distribute incoming traffic across different aisles, thereby reducing congestion at each entrance.

\begin{figure}[!ht]
\centering
\hspace{-0.5 cm}
\subfigure[Division of outlet areas in the network\label{fig:outlet_1}] {\includegraphics[width= 5.5 cm]{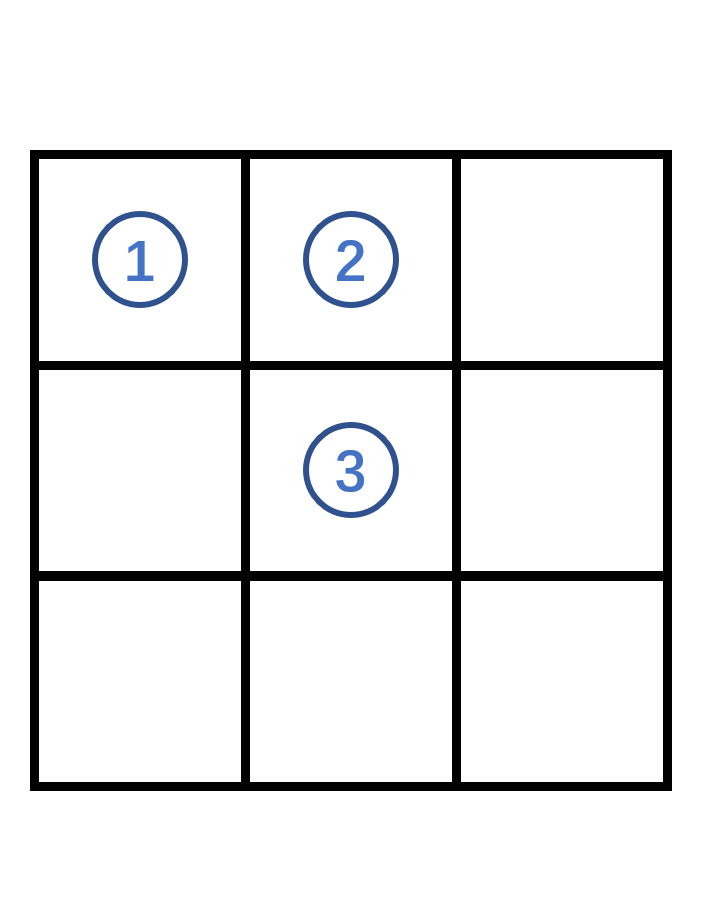}}
\hspace{0.5 cm}
\subfigure[RC-S performance in different scenarios\label{fig:outlet_2}] {\includegraphics[width = 9.5 cm]{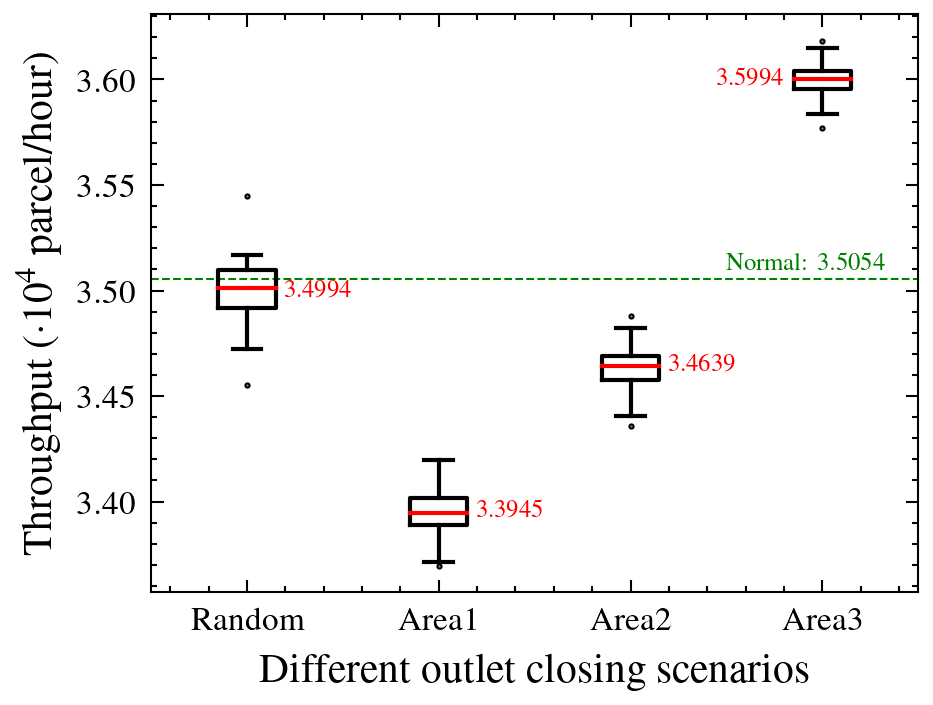}}
\hspace{-0.5 cm}
\caption{The effect of closing outlets in the network}
\label{fig:closing scenarios}  
\end{figure}

We next evaluate the effect of closing a subset of outlets on system throughput. The network size is set to $n_h=n_v=20$, with other experimental parameters kept constant. We consider four scenarios for closing outlets: (1) randomly closing outlets; (2) closing outlets in the corners; (3) closing outlets along the edges; (4) closing outlets in the center. To control for other factors, the number of closed outlets in each scenario is set to one ninth of the total. Scenarios (2), (3), and (4) correspond to Areas~1, 2, and 3 in Figure~\ref{fig:outlet_1}, respectively. Results from 50 independent experiments for each scenario are shown in Figure~\ref{fig:outlet_2}. The green dashed line represents the average throughput when all outlets are open. It can be observed that randomly closing outlets does not have a significant impact on system throughput. Closing outlets in Areas~1 and 2 results in a decrease in average throughput, while closing outlets in Area~3 leads to a slight increase in average throughput. Intuitively, closing outlets in a specific area causes an uneven distribution of demand, leading to imbalanced network traffic and increased traffic load in other areas of the network. However, from another perspective, it reduces robot turns within the closed area, thereby improving the utilization of spatio-temporal resources. Since most of the shortest paths pass through the central area, closing central outlets in scenario 4 results in higher throughput compared to the normal case.

\begin{figure}[ht!]
	\centering
	\includegraphics[width = 10cm]{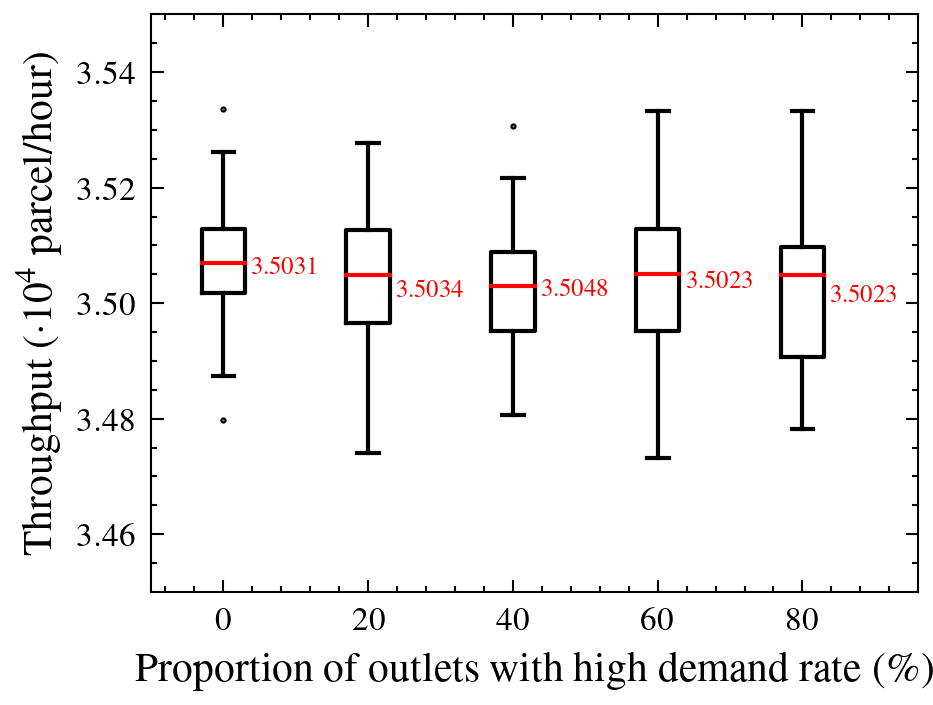}
	\caption{The effect of heterogeneous demand distribution}
	\label{fig:hete_demand}
\end{figure}

In addition, we evaluate the effect of heterogeneous demand distribution on system efficiency. In real sorting centers, sorting destinations with higher demand arrival rates are often assigned multiple outlets to minimize differences in the visit frequency between outlets; however, this cannot completely eliminate demand heterogeneity. We assume there are two types of outlets in the RSS: one type has an average demand arrival rate twice that of the other, with both following a normal distribution. For outlets with extremely low arrival rates, the situation is similar to the previous experiments involving closed outlets, and their impact can be minimized by random assignment. The two types of outlets are spatially uniformly distributed throughout the network. Keeping other experimental parameters constant, we adjust the proportion of the outlet types, and the resulting system performance is shown in Figure~\ref{fig:hete_demand}. The results indicate that, with a fixed total demand and varying proportions of high-demand outlets, RC-S maintains its efficiency, and the impact of demand heterogeneity on system throughput is negligible. This is achieved by leveraging demand arrival rates for different destinations and balancing the visit frequency to each area, ensuring a well-distributed traffic flow throughout the network.

\begin{table*}[!htbp]
	\centering
	\caption{RC-S performance in large-scale systems}
	\begin{tabular}{cccccc}
		\toprule
		\multirow{2}{*}{$n_h,n_v$} & Throughput  & Avg. path  & Avg. time  & Max. time & \multirow{2}{*}{\# Feasible paths} \\
            & ($\cdot {10}^4$/h) & length (m) & per cycle (s) & per cycle (s) & \\
		\midrule
		22 & 3.871 & 27.832 & 0.0359 & 0.0942 & 118,136 \\
            24 & 4.331 & 30.062 & 0.0480 & 0.1255 & 167,084 \\
            26 & 4.699 & 32.517 & 0.0614 & 0.2393 & 229,888 \\
            28 & 5.156 & 34.734 & 0.0836 & 0.2436 & 308,948 \\
            30 & 5.517 & 37.261 & 0.1004 & 0.2820 & 406,856 \\
		\bottomrule
	\end{tabular}
	\label{tab:table12}
\end{table*}

Finally, we evaluate the performance of RC-S in large-scale systems. The experimental scenarios include five large-scale networks with $n_h = n_v = 22,\ 24,\ 26,\ 28$, and $30$, respectively. Among these, the $30+30$ network covers an actual area of over 5{,}000 square meters, significantly larger than most sorting systems in modern warehouses. For example, a real case from Deppon Express mentioned in the study by \cite{zou2021robotic} has a network with a scale of $18+6$, approximately one-quarter of the $30+30$ scenario (assuming the same aisle width and loading station space occupation). The results are shown in Table~\ref{tab:table12}, including system throughput, the average and maximum run-time per cycle, and the number of precomputed feasible paths associated with the network size. As the network scale increases, the total number of feasible paths grows, and the path selection time per cycle gradually increases. However, the maximum run-time remains acceptable and is much shorter than the length of each cycle.
	\section{Layout Design Optimization} \label{sec:layout}
The layout design problem aims to balance facility costs and long-term operating costs. A larger facility will achieve higher maximum throughput; however, it will also increase construction and equipment costs. Considering transportation costs and supply stability, we analyze the following scenario: A warehouse signs a long-term RaaS contract with a robotics company, with each year divided into several operating periods denoted by $\sigma$. The proportion of each period's duration within the year is represented by $\theta^\sigma$, and the planned throughput level during each period is $T^\sigma$. The facility cost $C_f$ and expected operations cost $C_o$ are calculated as follows:

\begin{align}
	&C_f = P_s \cdot [2D \cdot (n_h-1) + W_w + W_l]\cdot [2D \cdot (n_v-1) + W_w + W_l] + P_{l} \cdot n_l \label{eq:C_s}\\
	&C_o = \sum_{\sigma \in \mathcal{S}} \theta^{\sigma} (P_w n_w^{\sigma} + P_r n_r^{\sigma}) \label{eq:C_o}
\end{align}
Parameters $P_s$ and $P_{l}$ in Equation~(\ref{eq:C_s}) denote the discounted site rental cost per square meter and equipment cost per loading station, respectively, while $P_w$ and $P_{r}$ in equation (\ref{eq:C_o}) denote the discounted unit labor cost and the unit rental cost of a robot. The objective of the problem is modeled as follows:
\begin{align}
	&C_{d} = C_f + C_o \label{objective}
\end{align}
The trade-off between these two types of costs in the objective function is the focus of the site planning stage. We formulate the layout design problem (LDP) as the following integer programming model:
\begin{align}
	\textbf{(LDP)}& \nonumber\\
	&\underset{n_h,n_v, n_l, n_w^\sigma,n_r^\sigma}{\min}\ C_{d} \nonumber \\
	\text{s.t.} \quad & \text{Constraints}\ (\ref{eq:alpha}),(\ref{eq:N_VP_final})\text{--}(\ref{eq:T_O_final}),(\ref{eq:C_s})\text{--}(\ref{objective}) \nonumber\\
	& \tilde{T}_{O}(n_h,n_v,n_w^\sigma,n_r^\sigma) \geq T^\sigma, \qquad \forall \sigma \in \mathcal{S}\label{const:THdemand} \\
	& (n_h - 1)(n_v -1 ) \geq N_{o} \label{const:outlet}\\
	& n_l \leq n_h + n_v \label{const:number_of_station}\\
	& n_w^\sigma \leq n_l, \qquad \forall \sigma \in \mathcal{S} \label{const:number_of_worker_s}\\
	& n_h = 2k_h\ , \qquad k_h \in \mathcal{Z}_{+} \label{const:even_1}\\
	& n_v = 2k_v\ , \qquad k_v \in \mathcal{Z}_{+} \label{const:even_2}\\
    & n_w^\sigma, n_r^\sigma \in \mathcal{Z}_+, \qquad \forall \sigma \in \mathcal{S}\label{const:variable_LDP}
\end{align}
Constraint~(\ref{const:THdemand}) ensures that the system throughput meets the sorting requirements in different periods. Constraint (\ref{const:outlet}) imposes a lower bound on the site size to accommodate an sufficient number of outlets $N_o$. Constraint~(\ref{const:number_of_station}) represents the maximum number of loading stations that can be accommodated in the aisle network under RC-S. Constraint~(\ref{const:number_of_worker_s}) limits the maximum number of workers to be stationed. Constraints (\ref{const:even_1})-(\ref{const:even_2}) ensure that the numbers of aisles are even. The remaining constraint~(\ref{const:variable_LDP}) requires the decision variables to be positive integers.

It is difficult to solve LDP because of the non-linearity of the inequality constraints, especially for large-scale problems. In this research, we apply the method of penalty successive linear programming (PSLP) to solve the LP relaxation of LDP, which exhibits good robustness and convergence properties for large-scale problems \citep{bazaraa2013nonlinear}. Specifically, PSLP sequentially solves a linearized feasible direction finding subproblem augmented with the penalty function, and utilizes the concept of a trust region (updated at each iteration) to control the step size. To ensure that the constraints continuous and differentiable, we introduce additional constraints requiring that the number of robots is less than the number of available VPs, given by
\begin{align}
	& n_{r}^\sigma - \beta(n_h,n_v)\cdot (1-{(1-\frac{n_w^\sigma}{n_h+n_v})}^2)\cdot \frac{2\tau_e \cdot [n_h(n_v-1) + n_v(n_h-1)]}{\tau_c} \leq 0 \quad \forall \sigma \in \mathcal{S}. \label{ineq:g1}
\end{align}

At the end of the algorithm, the optimal solution of the relaxed problem is converted into a feasible solution for the original problem by rounding up. The details of the PSLP algorithm are shown in Appendix~C.

	\section{Numerical Examples for Layout Design} \label{sec:numerical}
In this section, we conduct a sensitivity analysis on the site rental cost per square meter $P_s$ and the labor cost per man-month $P_w$, focusing on the trend of total costs under different scenarios as the sorting throughput level changes. Due to increased competition in the robot market, prices remain relatively stable across different companies. Consequently, we do not investigate variations in the robot rental cost in the experiments, nor do we examine the equipment cost at loading stations. We consider a five-year investment plan with a monthly interest rate $\gamma_0 = 0.5\%$, and the discounted costs of different components in LDP are calculated as:
\begin{align}
P_i = \sum_{t=0}^{60} \frac{M_i}{{(1+\gamma_0)}^t},\quad i = s,l,w,r\label{eq:gamma}
\end{align}
where $M_i$ denotes the unit cost of component $i$ per month. We estimate the monthly warehousing rental cost using JD Logistics' financial reports for the first three quarters of 2023. The labor cost is sourced from data related to warehousing job recruitment in Beijing, obtained from a job-posting website (58.com). The data regarding the equipment cost of loading stations and the rental cost of robots are obtained from the Geekplus company. 

We assume that the warehouse only expands its sorting capacity during a few major shopping events, meaning that there are two typical levels of sorting demand, namely the average off-peak season demand ${T}^L$ and the average peak season demand ${T}^H$, with ${T}^L=0.8\cdot {T}^H$. The ratios of the two periods are set as $(\theta^L, \theta^H) = (5/6,1/6)$, which means that each quarter typically includes a two-week shopping event. The default values of the parameters are shown in Table~\ref{tab:table71} and the configurations of RC-S remain the same as in Table~\ref{tab:table51}. Our algorithms are implemented in Python~3.11 with the solver Gurobi Optimizer~9.0, and all experiments are conducted on a personal computer running Windows~11 with an Intel i7-9700F CPU and 16~GB of RAM. All processes are run single-threaded.

\begin{table}[!htbp]
	\centering
	\caption{Default parameter values  in Section~\ref{sec:numerical}}
	\begin{tabular}{cccccc}
		\toprule
		$M_{s}$ (CNY / $m^2\cdot$ mo) & $M_{l}$  (CNY/mo) & $M_w$ (CNY/mo) & $M_{r}$ (CNY/mo) & $W_w$ (m) & $W_l$ (m) \\
		\midrule
		10 & 400 & 5,000 & 200 & 5 & 5\\
		\bottomrule
	\end{tabular}
	\label{tab:table71}
\end{table}

\subsection{Sensitivity analysis on site rental cost}
We begin by investigating the impact of site rental cost under two different scenarios: one with $N_{o}=100$ (scenario 1) and another with $N_{o}=400$ (scenario 2). The monthly rental cost of floor space per square meter varies from 10~CNY to 30~CNY. The results are presented in Figures~\ref{fig:sitecost1} and ~\ref{fig:sitecost2}. The findings reveal that:

\begin{figure}[hbt!]
	\centering
	\subfigure[$N_{o}=100$ (scenario 1)] {\includegraphics[width=0.95\linewidth]{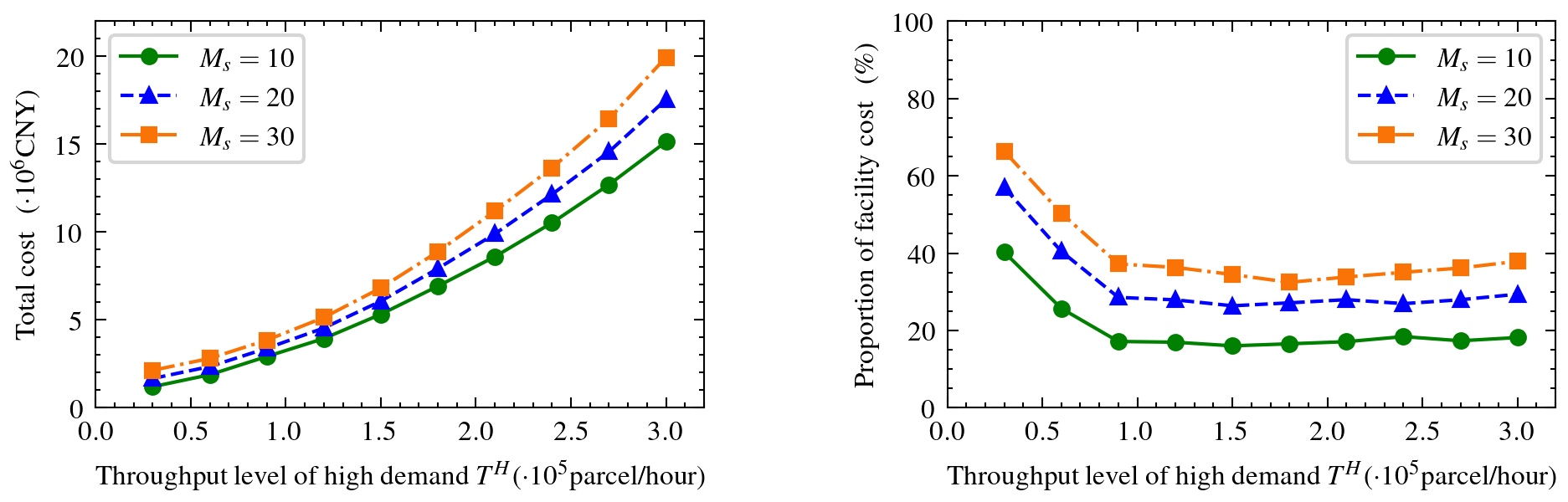}\label{fig:sitecost1}} 
	\qquad
	\subfigure[$N_{o}=400$ (scenario 2)] {\includegraphics[width=0.95\linewidth]{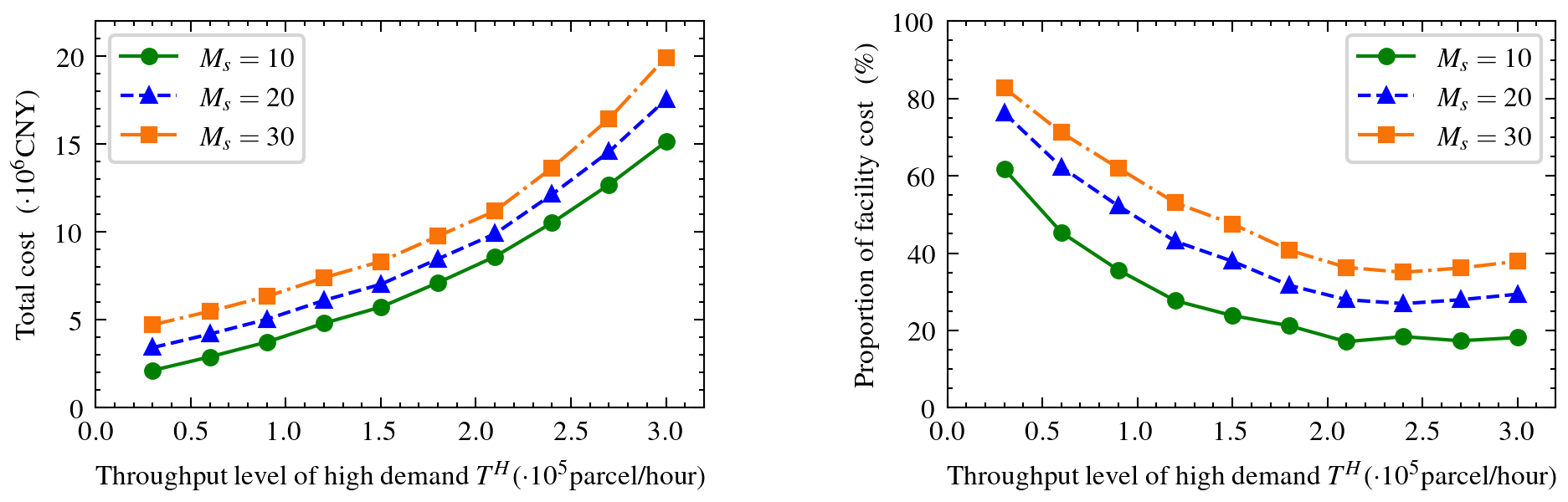}\label{fig:sitecost2}} 
	\caption{Comparison of total costs and facility cost proportions on different site rental costs}
	\label{fig:sitecost}
\end{figure}

\begin{table}[!htbp]
	\centering
	\caption{Results for scenario $M_s=10$, $N_{o} = 100$}
	\begin{tabular}{cccccccc}
		\toprule
		${T}^H$ & $(n_h,n_v)$  & $(n_w^L,n_{r}^L)$ & $(n_w^H,n_{r}^H)$ & $C_f$  & $C_o$ & $C_d$ & Proportion of \\
		&&&&$(\times {10}^6)$&$(\times {10}^6)$&$(\times {10}^6)$ & site rental cost (\%)\\
		\midrule
3000 & (10,12) & (2,10) & (3,13) & 0.53 & 0.67 & 1.20 & 38.80 \\
6000 & (10,12) & (4,21) & (5,26) & 0.57 & 1.31 & 1.88 & 24.78 \\
9000 & (10,12) & (7,33) & (9,42) & 0.65 & 2.26 & 2.92 & 15.96 \\
12000 & (12,14) & (9,51) & (13,67) & 0.87 & 3.07 & 3.94 & 15.20 \\
15000 & (14,16) & (12,75) & (18,98) & 1.12 & 4.20 & 5.32 & 14.07 \\
18000 & (18,18) & (15,109) & (23,142) & 1.48 & 5.44 & 6.92 & 14.54 \\
21000 & (20,22) & (18,147) & (28,194) & 1.88 & 6.72 & 8.60 & 15.09 \\
24000 & (24,26) & (21,199) & (32,261) & 2.41 & 8.11 & 10.52 & 16.60 \\
27000 & (26,26) & (25,236) & (44,315) & 2.79 & 9.91 & 12.70 & 14.74 \\
30000 & (30,30) & (29,302) & (49,401) & 3.42 & 11.72 & 15.14 & 15.88 \\
		\bottomrule
	\end{tabular}
	\label{tab:table72}
\end{table}

\begin{itemize}
	\item When the target throughput level is low, the total cost in scenario 1 is significantly lower than that in scenario 2. This is because, to accommodate more outlets, the latter requires renting a larger initial space. In both scenarios, the curve of the optimal system cost exhibits a turning point $T^*$. When $T^H > T^*$, the cost escalation rate increases significantly. The turning point indicates that the initial space in the current scenario is no longer sufficient to accommodate a larger throughput, necessitating an expansion of the aisle network scale. In scenario~1, this turning point occurs at $T^*=9{,}000$, whereas in scenario 2, the turning point is at $T^*=21{,}000$.  
	\item The proportion of facility costs increases as $M_s$ increases, and this effect is particularly pronounced at lower sorting throughput levels. As the throughput level becomes higher, this proportion gradually decreases and stabilizes after the turning point $T^*$ in each case. By examining the system configurations in each dataset, we find that this phenomenon arises because, at higher throughput demands, efficiency improvements necessitate expanding the space, which leads to an increase in the site rental cost. Table~\ref{tab:table72} shows the optimal system configuration and the corresponding costs in the case of $M_s =10$ and $N_{o}=100$. When $T^H > T^*$, the optimal solutions not only assign more workers and robots to the loading stations, but also expand the size of the network to alleviate traffic pressure.
 
    \item Comparing the data for $T^H=15{,}000$ and $T^H=30{,}000$ in Table~\ref{tab:table72}, the total cost for the latter is nearly three times that for the former. This indicates that due to the simultaneous growth of $C_f$ and $C_o$, the growth rate of the total cost exceeds the growth rate of throughput.
\end{itemize}

From the results of the site rental cost analysis, the derived managerial insights are summarized as follows. 
\begin{insight}
	The advantage of RSS is evident in its low initial investment in scenarios where both throughput demand and the number of sortation categories are low.
\end{insight}
\begin{insight}
	High density leads to decreased worker efficiency. Instead of running an RSS at full capacity, it is more cost-effective to appropriately expand the site size and reduce the proportion of activated loading stations.
\end{insight}

\begin{figure}[hbt!]
	\centering
	\subfigure[$N_{o}=100$ (scenario 1)] {\includegraphics[width=0.95\linewidth]{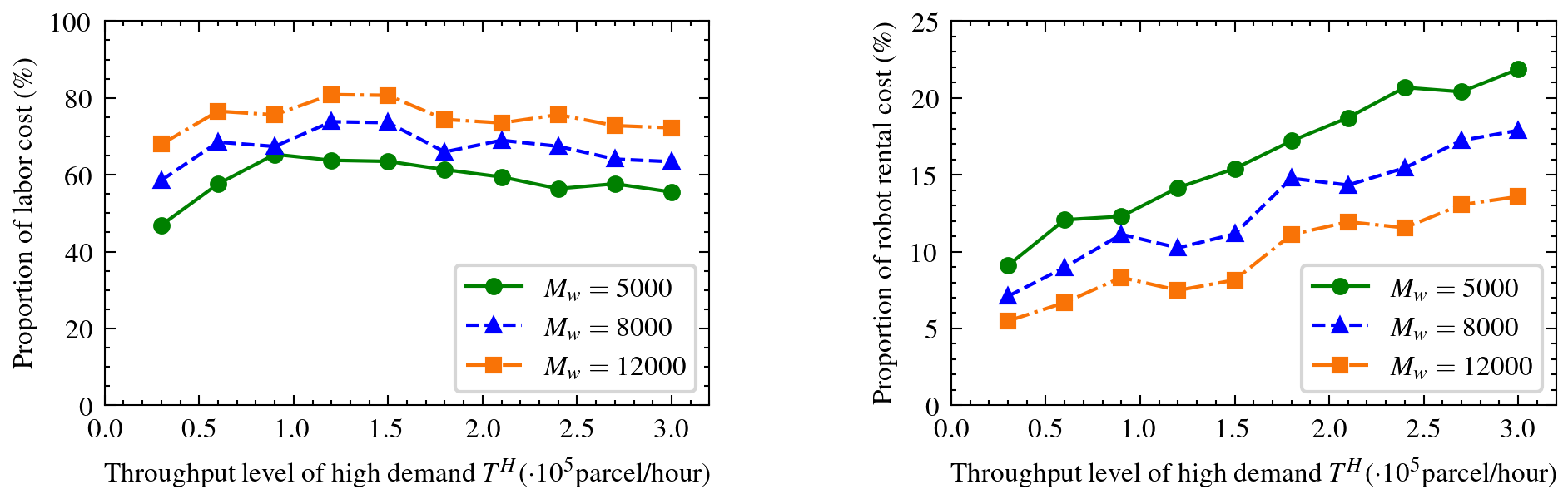}\label{fig:employmentcost1}} 
	\qquad
	\subfigure[$N_{o}=400$ (scenario 2)] {\includegraphics[width=0.95\linewidth]{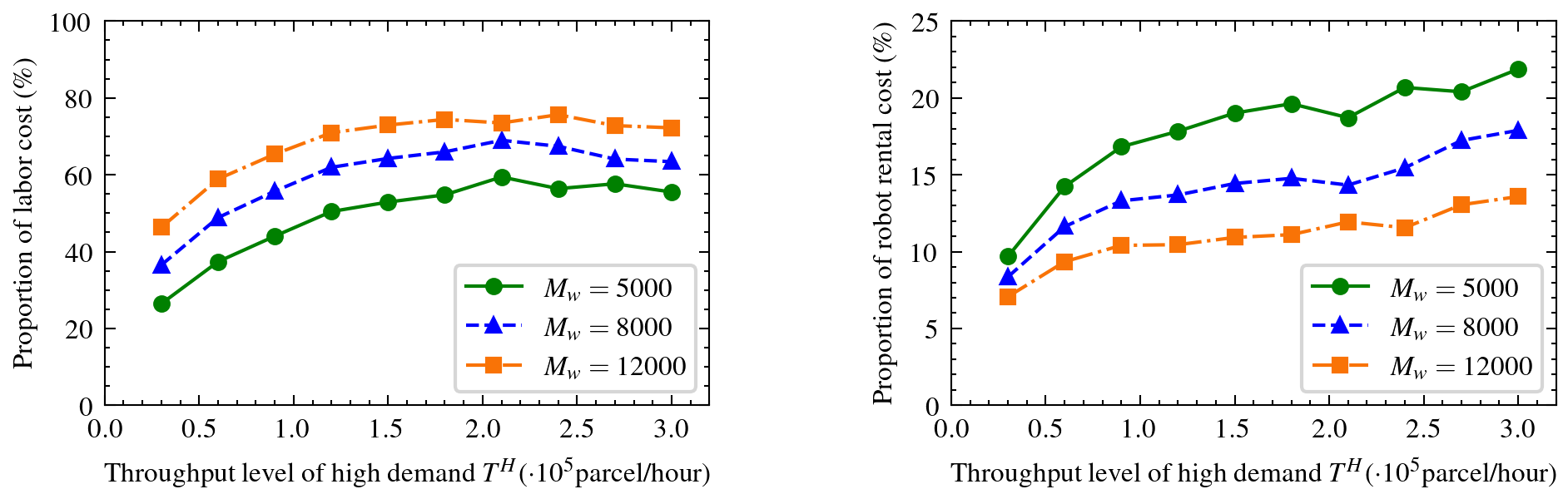}\label{fig:employmentcost2}} 
	\caption{Proportions of labor cost and robot rental cost under different unit labor costs}
	\label{fig:employmentcost}
\end{figure}
\subsection{Sensitivity analysis on labor cost}
Similarly to our previous experiments, we now explore the influence of different labor costs on the total cost in two scenarios with $N_{o}=100$ and $N_{o}=400$. The monthly labor cost per worker varies from 5{,}000~CNY to 12{,}000~CNY. The proportions of labor cost and robot rental cost are both visualized in Figures~\ref{fig:employmentcost1} and \ref{fig:employmentcost2}. Our observations are summarized as follows:

\begin{itemize}
	\item As $T^H$ increases, the proportion of labor cost gradually increases and eventually stabilizes. In every scenario, the proportion of labor costs in the optimal solution for $T^H=30{,}000$ exceeds 50\%.
	\item The rental cost of robots gradually increases with the level of throughput in each scenario. This is because, to maintain traffic efficiency, the number of robots needs to increase with the expansion of the network. Meanwhile, as discussed in the performance analysis in Section~\ref{sec:analysis}, the rate of traffic capacity loss also increases, resulting in a considerable robot rental cost.
    \item Comparing the two scenarios, it is evident that in larger-scale sorting sites (scenario~2), the initial increase in the proportion of robot rental cost is faster, making it challenging for RSS systems to balance the advantages of high scalability and low cost.
\end{itemize}

From the labor cost analysis, we summarize the following managerial insight: 
\begin{insight}
	At lower levels of throughput, facility cost is the primary expense of RSS, whereas at higher levels of throughput, labor cost becomes the primary expense.
\end{insight}

Additionally, the cost structure undergoes notable changes with varying sorting demand. In scenarios with high sorting demand, the proportion and growth rate of robot rental cost become more significant. This reveals a distinctive feature of RSS systems: Robots, as carriers of parcels, become less efficient in larger sites because they have to cover increased distances. Merely enhancing management methods is insufficient to address the challenge of diminishing traffic efficiency. Today's commonly used conveyor sorters, equipped with high-power motors, require substantial initial investments but experience only slight increases in total cost as the sorting demand grows \citep{zou2021robotic}. Furthermore, the optimal system configurations of a conveyor sorter are not easily affected by fluctuations in the unit prices of cost components \citep{russell2003cost}, which makes it suitable for long-term operation in large-scale sorting scenarios. Therefore, we conclude the following insight regarding the application of RSS:
\begin{insight}
	Due to the traffic limitations of robots, RSS is suitable for application in small-scale scenarios, such as distribution centers at the end of the supply chain.
\end{insight}

	\section{Conclusion} \label{sec:conclusion}
This study conducts a comprehensive analysis of RSS in modern warehouses, encompassing robot traffic management, system efficiency, and the cost composition and layout design of the system. By incorporating throughput estimation, this research provides guidance for performance prediction and resource allocation in modern warehouses or factories that utilize numerous robots, promoting the ongoing development of automation in the logistics industry. 

In the operations stage, efficient robot scheduling emerges as a crucial element for unlocking the business value of RSS. We propose an innovative RC-S scheme, which serves as a framework for managing a large number of robots simultaneously and efficiently. We provide a detailed description of the phase structure and the parameter-setting method in RC-S. Additionally, we present a mathematical programming model, FPA, to minimize the travel cost at each cycle, and a heuristic algorithm as an online solver. We then conduct a theoretical analysis of the efficiency of RC-S, exploring the impact of various system configurations on throughput. In the validation section, we first compare our proposed traffic management framework with the classical cooperative A* algorithm. The simulation results indicate that our control method achieves a higher level of performance and computational efficiency. Furthermore, we validate the throughput estimation formula. The results show that queueing models can produce significant errors in some scenarios, whereas our model fits the real traffic flows of robots more precisely.

In the site planning stage, warehouse managers need to make decisions that involve a trade-off between initial facility investment $C_f$ and discounted operations costs $C_o$. We analyze the cost composition of RSS and propose a layout design optimization model, LDP, that minimizes total system costs. Through a detailed sensitivity analysis examining the cost dynamics of site rental and labor, we investigate cost proportions at different throughput levels. Specifically, at lower throughput, the rental cost for the site represents a considerable proportion of the total cost. However, at higher throughput, labor expenditure becomes the predominant cost component. The results also support our analysis of the model properties. Key insights can guide managers in understanding the investment and returns associated with the application of RSS, thus reducing the total cost of the warehouse.

For future work, there are two directions to explore. First, optimizing layout design by considering non-uniform demand and various outlet distributions would further improve the model's effectiveness. Second, extending the analysis to cover the entire process from order fulfillment to sorting in modern warehouses would provide valuable insights for enhancing system performance. 

        \section*{Declaration of generative AI and AI-assisted technologies in the writing process}
During the preparation of this work the authors used ChatGPT in order to check for grammatical errors and improve readability. After using this tool, the authors reviewed and edited the content as needed and take full responsibility for the content of the publication.

	\bibliographystyle{apalike}
	\bibliography{ref}

@article{fang2025dynamic,
  title={\href{https://doi.org/10.1287/trsc.2023.0458}{Dynamic Robot Routing and Destination Assignment Policies for Robotic Sorting Systems}},
  author={Fang, Yuan and De Koster, Ren{\'e} and Roy, Debjit and Yu, Yugang},
  journal={Transportation Science},
  year={2025},
  publisher={INFORMS}
}

@book{buzacott1993stochastic,
  title={Stochastic models of manufacturing systems},
  author={John A. Buzacott, J. George Shanthikumar},
  publisher={Prentice Hall},
  year={1993}
}

@article{russell2003cost,
  title={\href{https://doi.org/10.1080/07408170304358}{Cost and throughput modeling of manual and automated order fulfillment systems}},
  author={Russell, Mardi L and Meller, Russell D},
  journal={IIE Transactions},
  volume={35},
  number={7},
  pages={589--603},
  year={2003}
}

@inproceedings{ma2019searching,
  title={\href{doi.org/10.1609/aaai.v33i01.33017643}{Searching with consistent prioritization for multi-agent path finding}},
  author={Ma, Hang and Harabor, Daniel and Stuckey, Peter J and Li, Jiaoyang and Koenig, Sven},
  booktitle={Proceedings of the AAAI conference on artificial intelligence},
  volume={33},
  pages={7643--7650},
  year={2019}
}

@inproceedings{phillips2011sipp,
  title={\href{10.1109/ICRA.2011.5980306}{Sipp: Safe interval path planning for dynamic environments}},
  author={Phillips, Mike and Likhachev, Maxim},
  booktitle={2011 IEEE international conference on robotics and automation},
  pages={5628--5635},
  year={2011},
  organization={IEEE}
}

@inproceedings{Silver_2021,
  title={\href{10.1609/aiide.v1i1.18726}{Cooperative pathfinding}},
  author={Silver, David},
  booktitle={Proceedings of the aaai conference on artificial intelligence and interactive digital entertainment},
  volume={1},
  number={1},
  pages={117--122},
  year={2005}
}

@article{de2007design,
  title={\href{https://doi.org/10.1016/j.ejor.2006.07.009}{Design and control of warehouse order picking: A literature review}},
  author={De Koster, Ren{\'e} and Le-Duc, Tho and Roodbergen, Kees Jan},
  journal={European journal of operational research},
  volume={182},
  number={2},
  pages={481--501},
  year={2007},
  publisher={Elsevier}
}

@inproceedings{stern2019multi,
  title={\href{doi.org/10.1609/socs.v10i1.18510}{Multi-agent pathfinding: Definitions, variants, and benchmarks}},
  author={Stern, Roni and Sturtevant, Nathan and Felner, Ariel and Koenig, Sven and Ma, Hang and Walker, Thayne and Li, Jiaoyang and Atzmon, Dor and Cohen, Liron and Kumar, TK and others},
  booktitle={Proceedings of the International Symposium on Combinatorial Search},
  volume={10},
  pages={151--158},
  year={2019}
}

@Online{ReginaLoBiondo,
 author = {BusinessWire},
 year = {2023},
 title = {New Report Reveals Consumers Embracing Return Fees In Exchange For Convenient, Premium Offerings},
 url = {https://www.businesswire.com/news/home/20230110005398/en/New-Report-Reveals-Consumers-Embracing-Return-Fees-In-Exchange-For-Convenient-Premium-Offerings},
 urldate = {2023-01-10}
}

@article{wurman2008coordinating,
  title={\href{https://doi.org/10.1609/aimag.v29i1.2082}{Coordinating hundreds of cooperative, autonomous vehicles in warehouses}},
  author={Wurman, Peter R and D'Andrea, Raffaello and Mountz, Mick},
  journal={AI magazine},
  volume={29},
  number={1},
  pages={9--9},
  year={2008}
}

@article{BOYSEN2019396,
title = {\href{https://doi.org/10.1016/j.ejor.2018.08.023}{Warehousing in the e-commerce era: A survey}},
journal = {European Journal of Operational Research},
volume = {277},
number = {2},
pages = {396-411},
year = {2019},
issn = {0377-2217},
author = {Nils Boysen and René {de Koster} and Felix Weidinger}
}

@article{BOYSEN2023582,
title = {\href{https://doi.org/10.1016/j.ejor.2023.03.037}{Robotized sorting systems: Large-scale scheduling under real-time conditions with limited lookahead}},
journal = {European Journal of Operational Research},
volume = {310},
number = {2},
pages = {582-596},
year = {2023},
issn = {0377-2217},
author = {Nils Boysen and Stefan Schwerdfeger and Marlin {W. Ulmer}}
}

@article{pub99983,
  title={\href{https://doi.org/10.1287/trsc.2018.0873}{Robotized and automated warehouse systems: Review and recent developments}},
  author={Azadeh, Kaveh and De Koster, Ren{\'e} and Roy, Debjit},
  journal={Transportation Science},
  volume={53},
  number={4},
  pages={917--945},
  year={2019},
  publisher={INFORMS}
}

@article{XU2022102808,
title = {\href{https://doi.org/10.1016/j.tre.2022.102808}{Assignment of parcels to loading stations in robotic sorting systems}},
journal = {Transportation Research Part E: Logistics and Transportation Review},
volume = {164},
pages = {102808},
year = {2022},
issn = {1366-5545},
author = {Xianhao Xu and Yuerong Chen and Bipan Zou and Yeming Gong}
}

@article{Yu_LaValle_2013, 
title={\href{10.1609/aaai.v27i1.8541}{Structure and Intractability of Optimal Multi-Robot Path Planning on Graphs}}, 
volume={27}, 
number={1}, 
journal={Proceedings of the AAAI Conference on Artificial Intelligence}, 
author={Yu, Jingjin and LaValle, Steven}, 
year={2013}, 
month={Jun.}, 
pages={1443-1449} 
}

@article{poms13626,
author = {Shi, Ye and Yu, Hu and Yu, Yugang and Yue, Xiaohang},
title = {\href{doi.org/10.1111/poms.13626}{Analytics for IoT-Enabled Human–Robot Hybrid Sortation: An Online Optimization Approach}},
journal = {Production and Operations Management},
year = {2021},
volume = {},
number = {},
pages = {},
}

@article{TAN2021101279,
title = {\href{doi.org/10.1016/j.aei.2021.101279}{Optimizing parcel sorting process of vertical sorting system in e-commerce warehouse}},
journal = {Advanced Engineering Informatics},
volume = {48},
pages = {101279},
year = {2021},
issn = {1474-0346},
author = {Zheyi Tan and Haolin Li and Xueting He}
}

@article{Surynek_2010, 
title={\href{10.1609/aaai.v24i1.7767}{An Optimization Variant of Multi-Robot Path Planning Is Intractable}}, 
volume={24}, 
number={1}, 
journal={Proceedings of the AAAI Conference on Artificial Intelligence}, 
author={Surynek, Pavel}, 
year={2010}, month={Jul.}, 
pages={1261-1263} }

@article{dekhne2019automation,
  title={Automation in logistics: Big opportunity, bigger uncertainty},
  author={Dekhne, Ashutosh and Hastings, Greg and Murnane, John and Neuhaus, Florian},
  journal={McKinsey Q},
  year={2019},
pages = {1-12}
}

@article{BOYSEN2019796,
title = {\href{https://doi.org/10.1016/j.ejor.2018.08.014}{Automated sortation conveyors: A survey from an operational research perspective}},
author = {Nils Boysen and Dirk Briskorn and Stefan Fedtke and Marcel Schmickerath},
journal = {European Journal of Operational Research},
volume = {276},
number = {3},
pages = {796-815},
year = {2019},
issn = {0377-2217},
}

@INPROCEEDINGS{wagner6095022,
  author={Wagner, Glenn and Choset, Howie},
  booktitle={2011 IEEE/RSJ International Conference on Intelligent Robots and Systems}, 
  title={\href{10.1109/IROS.2011.6095022}{M*: A complete multirobot path planning algorithm with performance bounds}}, 
  year={2011},
  volume={},
  number={},
  pages={3260-3267}
}

@inproceedings{nguyen2019generalized,
  title={\href{doi.org/10.1609/socs.v10i1.18487}{Generalized target assignment and path finding using answer set programming}},
  author={Nguyen, Van and Obermeier, Philipp and Son, Tran and Schaub, Torsten and Yeoh, William},
  booktitle={Proceedings of the International Symposium on Combinatorial Search},
  volume={10},
  number={1},
  pages={194--195},
  year={2019}
}

@inproceedings{Liu2019path,
author = {Liu, Minghua and Ma, Hang and Li, Jiaoyang and Koenig, Sven},
title = {\href{10.5555/3306127.3331816}{Task and Path Planning for Multi-Agent Pickup and Delivery}},
year = {2019},
isbn = {9781450363099},
publisher = {International Foundation for Autonomous Agents and Multiagent Systems},
address = {Richland, SC},
booktitle = {Proceedings of the 18th International Conference on Autonomous Agents and MultiAgent Systems},
pages = {1152–1160},
numpages = {9},
location = {Montreal QC, Canada},
series = {AAMAS '19}
}

@article{Li_Tinka2021, 
title={\href{10.1609/aaai.v35i13.17344}{Lifelong Multi-Agent Path Finding in Large-Scale Warehouses}}, 
volume={35},
number={13}, 
journal={Proceedings of the AAAI Conference on Artificial Intelligence}, 
author={Li, Jiaoyang and Tinka, Andrew and Kiesel, Scott and Durham, Joseph W. and Kumar, T. K. Satish and Koenig, Sven}, 
year={2021}, 
month={May}, 
pages={11272-11281} 
}

@article{RCII,
 title={\href{https://doi.org/10.1287/trsc.2021.1061}{Rhythmic Control of Automated Traffic—Part II: Grid Network Rhythm and Online Routing}},
 author={Lin, Xi and Li, Meng and Max, Shen Zuo-Jun and Yin, Yafeng and He, Fang},
 journal={Transportation Science},
 volume={55},
 number={5},
 pages={988--1009},
 year={2021},
 publisher={INFORMS}
}

@InProceedings{10.1007/978-981-16-8174-5_10,
author="Wang, Ke
and Liang, Wei
and Shi, Huaguang
and Zhang, Jialin
and Wang, Qi",
editor="Cui, Li
and Xie, Xiaolan",
title={\href{https://doi.org/10.1007/978-981-16-8174-5_10}{A Calculation Time Prediction-Based Multiflow Network Path Planning Method for the AGV Sorting System}},
booktitle="Wireless Sensor Networks",
year="2021",
publisher="Springer Singapore",
address="Singapore",
pages="123--135",
isbn="978-981-16-8174-5"
}

@article{zou2021robotic,
  title={\href{https://doi.org/10.1287/trsc.2021.1053}{Robotic Sorting Systems: Performance Estimation and Operating Policies Analysis}},
  author={Zou, Bipan and De Koster, Ren{\'e} and Gong, Yeming and Xu, Xianhao and Shen, Guwen},
  journal={Transportation Science},
  volume={55},
  number={6},
  pages={1430--1455},
  year={2021},
  publisher={INFORMS}
}

@INPROCEEDINGS{9532245,  
author={Zou, Bipan and Chen, Yuerong},  
booktitle={2020 7th International Conference on Information Science and Control Engineering (ICISCE)},   
title={\href{10.1109/ICISCE50968.2020.00079}{Assign orders to workstations in robotic sorting systems considering destination information}},   
year={2020},  
volume={},  
number={},  
pages={343-348}}

@article{liu2019multi,
  title={\href{https://doi.org/10.1371/journal.pone.0226161}{Multi-objective AGV scheduling in an automatic sorting system of an unmanned (intelligent) warehouse by using two adaptive genetic algorithms and a multi-adaptive genetic algorithm}},
  author={Liu, Yubang and Ji, Shouwen and Su, Zengrong and Guo, Dong},
  journal={PloS one},
  volume={14},
  number={12},
  pages={e0226161},
  year={2019},
  publisher={Public Library of Science San Francisco, CA USA}
}

@INPROCEEDINGS{9102081,  
author={Zi, Lanfang and Gao, Benhe},  
booktitle={2020 IEEE 7th International Conference on Industrial Engineering and Applications (ICIEA)},  
title={\href{10.1109/ICIEA49774.2020.9102081}{Layout Planning Problem in an Innovative Parcel Sorting System Based on Automated Guided Vehicles}},   
year={2020},  
volume={},  
number={},  
pages={630-637}}

@INPROCEEDINGS{9216812,  
author={Zi, Lanfang and Gao, Benhe},  
booktitle={2020 IEEE 16th International Conference on Automation Science and Engineering (CASE)}, 
title={\href{10.1109/CASE48305.2020.9216812}{Performance Estimating in an Innovative AGVs-based Parcel Sorting System Considering the Distribution of Destinations}},   
year={2020}, 
volume={},  
number={},  
pages={1129-1134}}

@book{bazaraa2013nonlinear,
  title={\href{https://doi.org/10.1016/j.trb.2009.11.001}{Nonlinear programming: theory and algorithms}},
  author={Bazaraa, Mokhtar S and Sherali, Hanif D and Shetty, Chitharanjan M},
  year={2013},
  publisher={John Wiley \& Sons}
}

@article{DAGANZO2010434,
title = {\href{https://doi.org/10.1016/j.trb.2009.11.001}{Structure of competitive transit networks}},
journal = {Transportation Research Part B: Methodological},
volume = {44},
number = {4},
pages = {434-446},
year = {2010},
issn = {0191-2615},
author = {Carlos F. Daganzo}
}
	
	\newpage
	\section*{Appendix A: Main notations}
\setcounter{table}{0}
\begin{table}[!ht]
	\renewcommand{\thetable}{A-1}
	\centering
	\caption{Summary of main abbreviations}
	\begin{tabular} {c l}
		\hline
		\textbf{Abbreviations} & \textbf{Explanation} \\
		\hline
		RSS & robotic sorting system\\
            RMFS & robotic mobile fulfill system\\
            RaaS & Robot-as-a-Service\\
            MAPF & multi-agent path finding problem\\
            SOQN & semi-open queueing network\\
            OQN & open queueing network\\
            RC-S & rhythmic control for sorting scenario\\
            VP & virtual platoon\\
            FPA & feasible path assignment problem\\
            LDP & layout design problem\\
            PSLP & penalty successive linear programming\\
		\hline
	\end{tabular}
	\label{Abbreviation_tab}
\end{table}

\renewcommand{\thetable}{A-2}
\begin{longtable}{p{.20\textwidth} p{.80\textwidth}}
	\caption{Notations and explanation} \\
	\hline
	\textbf{Notations} & \textbf{Explanation} \\ \hline
	\multicolumn{2}{l}{\textit{Sets}} \\
		$\mathcal{G}$ & graph constituted of all nodes and links\\
            $\mathcal{V}$ & set of all nodes\\
		$\mathcal{V}_c$ & set of conflict nodes\\
		$\mathcal{V}_u$ & set of unloading nodes\\
            $\mathcal{V}_e$ & set of entrance/exit nodes\\
		$\mathcal{E}$ & set of all edges\\
		$\mathcal{L}$ & set of all loading stations\\
		$\mathcal{O}$ & set of all outlets\\
            $\mathcal{C}$ & set of all cycles starting from the current cycle\\
		$\mathcal{R}$ & set of all feasible paths\\
		$\mathcal{R}(i,j,k)$ & set of feasible paths connecting loading station $i\in \mathcal{L}$ and $j\in \mathcal{L}$, passing outlet $k\in \mathcal{O}$\\
            $\mathcal{R}_{i,k}$ & set of feasible path patterns starting from loading station $i\in \mathcal{L}$ and passing by outlet $k\in \mathcal{O}$, independent of starting cycle\\
            $\mathcal{R}_i^{(0)}$ & Set of feasible paths available to loading station $i$ at the beginning of current cycle, based on the current destination outlet of the pending sorting task. If no task is pending, $\mathcal{R}_i^{(0)} = \emptyset$ \\
            $\hat{\mathcal{R}}_i$ & set of feasible paths ending at loading station $i\in \mathcal{L}$\\
		$\mathcal{I}_x$ & set of locations of loading stations on the top and bottom sides\\
		$\mathcal{I}_y$ & set of locations of loading stations on the left and right sides\\
            $\mathcal{S}$ & set of operating periods\\
		\hline	
		\multicolumn{2}{l}{\textit{Parameters}} \\
		$D$ & length of each grid in the sorting zone\\
		$W_w$ & width of waiting zone\\
		$W_l$ & width of loading zone\\
		$C_f$ & the facility costs of RSS\\
		$C_o$ & the operations cost of RSS\\
		$C_d$ & the total costs of RSS\\
		$P_f$ & cost of facility per square floor space\\
		$P_{l}$ & cost of equipment in one loading station\\
		$P_{w}$ & cost of labor per person per month\\
		$P_{r}$ & cost of robot per vehicle per month\\
		${T}^\sigma$ & target throughput level in operating period $\sigma \in \mathcal{S}$\\
		$\theta^\sigma$ & Ratio of operating period $\sigma\in \mathcal{S}$\\
		$N_{o}$ & minimum number of outlets required to meet the sorting category demands\\
		$\gamma$ & weight of operations costs in optimization model\\
		$\tau_c$ & time interval of VPs in RC-S\\
		$\tau_e$ & travel time of VPs on each link in RC-S\\
            $v_{VP}$ & speed of VPs \\
            $v_{max}$ & maximum speed of robots\\
            $c_{max}$ & maximum acceleration/deceleration of robots\\
            $\omega_r$ & rotation speed of robots\\
            $r_{l}$ & maximum loading rate of a loading station\\
            $c_{i}^r$ & completion time of feasible path $r \in \mathcal{R}_i^{(0)}$, measured in number of cycles\\
            $\hat{c}_{i}$ & penalty for delaying the assignment of a feasible path to a sorting task from loading station $i\in \mathcal{L}$ until the next cycle\\
            $\delta_{i}^{r,t,\nu,l}$ & incidence between feasible path $r \in \mathcal{R}_i^{(0)}$ starting at cycle $t\in \mathcal{C}$, node $\nu \in \mathcal{V}$ and cycle $l\in \mathcal{C}$\\
            $N_\nu^l$ & remaining capacity in node $\nu\in \mathcal{V}$ in cycle $l\in \mathcal{C}$\\
            $d_{i}$ & sorting demand from loading station $i\in \mathcal{L}$, $d_{i} \in \{0,1\}$\\
		\hline	
		\multicolumn{2}{l}{\textit{Variables}} \\
            $x_{i}^{r,t}$ & decision on whether feasible path $r \in \mathcal{R}_i^{(0)}$ starting at cycle $t\in \mathcal{C}$ is reserved\\
            $\hat{x}_{i}$ & decision on whether the sorting task from loading station $i\in \mathcal{L}$ should be postponed to the next cycle\\
		$n_h,n_v$ & number of horizontal and vertical aisles in sorting zone\\
		$n_{l}$ & number of loading stations in loading zone\\
		$n_w^\sigma$ & number of workers in operating period $\sigma\in \mathcal{S}$\\
		$n_r^\sigma$ & number of robots in operating period $\sigma\in \mathcal{S}$\\
		
		\hline
	\label{Notataion_tab}
\end{longtable}
	\setcounter{figure}{0}
\section*{Appendix B: Proofs}
\subsection*{B.1 Proof of Proposition 1} 
Assume that the number of vertical and horizontal aisles in the network is $(n_v, n_h)$. We illustrate how each VP–cycle reservation in a feasible path can be translated into a set of node–cycle pairs, using the VPs moving along a horizontal aisle as an example.

For a horizontal aisle $a$, we denote the sequence of nodes along the aisle (from entrance node to exit node) as $\mathcal{V}_{aisle}^{(a)} = \{\nu_1^{(a)}, \nu_2^{(a)}, \nu_3^{(a)}, \ldots\}$. The total number of nodes in this aisle is:
\begin{eqnarray*}
    \lvert \mathcal{V}_{ailse}^{(a)} \rvert = 2n_v + 1
\end{eqnarray*}
Under the RC-S control scheme, each cycle is divided into four phases, and each phase has a duration of $\tau_e$. A VP advances to the next node along its path after each phase. As illustrated in Figure 7 (Section 4), a newly generated VP enters the entrance node of its aisle during the second phase of the cycle in which it is released. In the cycle when it exits the network, it reaches the exit node during the second phase.

Thus, for a VP $p_0$ generated in cycle $l_0$ along aisle $a$, the set of nodes it occupies in each cycle $l_i \in [l_0, l_0 + \frac{n_v}{2} - 1]$ is given by:
\begin{eqnarray*}
    (p_0, l_i) \iff \begin{cases}
        \{(\nu_1^{(a)},l_i), (\nu_2^{(a)},l_i), (\nu_3^{(a)},l_i)\}, & \text{if }\ l_i = l_0 \\
        \{(\nu_{4(l_i-l_0)}^{(a)},l_i), (\nu_{4(l_i-l_0)+1}^{(a)},l_i), (\nu_{4(l_i-l_0)+2}^{(a)},l_i), (\nu_{4(l_i-l_0)+3}^{(a)},l_i)\}, & \text{if }\ l_0 < l_i < l_0 + \frac{n_v}{2} \\
        \{(\nu_{2n_v}^{(a)},l_i), (\nu_{2n_v + 1}^{(a)},l_i)\}, & \text{if }\ l_i = l_0 + \frac{n_v}{2} \\
        \emptyset, & \text{otherwise.}
    \end{cases}
\end{eqnarray*}
Similarly, for a VP $p_0$ moving along a vertical aisle $a$, we have:
\begin{eqnarray*}
    (p_0, l_i) \iff \begin{cases}
        \{(\nu_1^{(a)},l_i), (\nu_2^{(a)},l_i), (\nu_3^{(a)},l_i)\}, & \text{if }\ l_i = l_0 \\
        \{(\nu_{4(l_i-l_0)}^{(a)},l_i), (\nu_{4(l_i-l_0)+1}^{(a)},l_i), (\nu_{4(l_i-l_0)+2}^{(a)},l_i), (\nu_{4(l_i-l_0)+3}^{(a)},l_i)\}, & \text{if }\ l_0 < l_i < l_0 + \frac{n_h}{2} \\
        \{(\nu_{2n_h}^{(a)},l_i), (\nu_{2n_h + 1}^{(a)},l_i)\}, & \text{if }\ l_i = l_0 + \frac{n_h}{2} \\
        \emptyset, & \text{otherwise.}
    \end{cases}
\end{eqnarray*}
\begin{flushright} 
	$\square$ 
\end{flushright}

\subsection*{B.2 Proof of Proposition 2} 
\noindent (\romannumeral1):

By definition 1, a feasible path of RC-S connects two active loading stations. If a feasible path has no more than two turns, all its turning points are within the blue region; otherwise, its endpoint will fall into a non-active loading station. The movement trajectories of occupied VPs will not cross the boundaries of the blue region; therefore, their upper limit is $\kappa$. Condition (\romannumeral1) is sufficient.

\noindent (\romannumeral2):

We derive the average road length that each loading station can be allocated, noted as $\tilde{d}_{ls}$:
\begin{eqnarray*}
    \tilde{d}_{ls}(\alpha) &=& \frac{\big(1-{(1-\alpha)}^2\big)\cdot n_h n_v}{\frac{\alpha(n_h + n_v)}{2}} \\
    &=& \frac{(2-\alpha)\cdot 2n_h n_v }{n_h + n_v} \geq (2-\alpha) \cdot min\{n_h, n_v\}
\end{eqnarray*}
When the average travel distance of a task is less than $\tilde{d}_{ls}$, the proportion of the union of segments passed by the occupied VPs in the total set of segments is less than $\kappa$. Condition (\romannumeral2) is sufficient.
\begin{flushright} 
	$\square$ 
\end{flushright}

\subsection*{B.3 Derivation of Maximum Travel Distance in Constraint (7)} 
\begin{figure}[hbt!]
    \renewcommand{\thefigure}{B-1}
	\centering
	\subfigure[] {\includegraphics[width=0.45\linewidth]{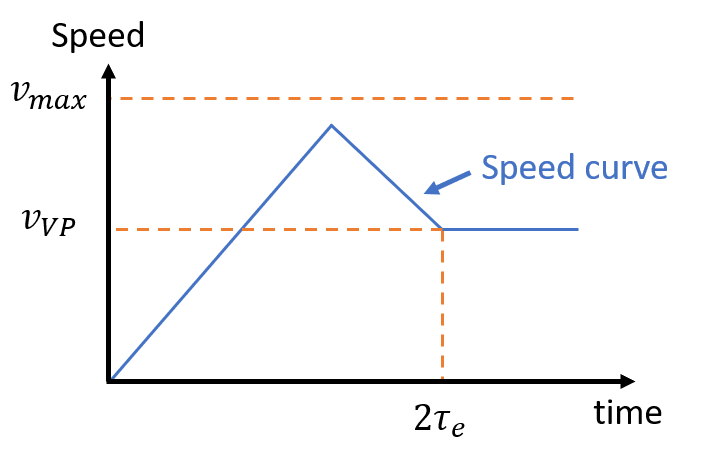}\label{fig:curve_1}} 
	\quad
	\subfigure[] {\includegraphics[width=0.46\linewidth]{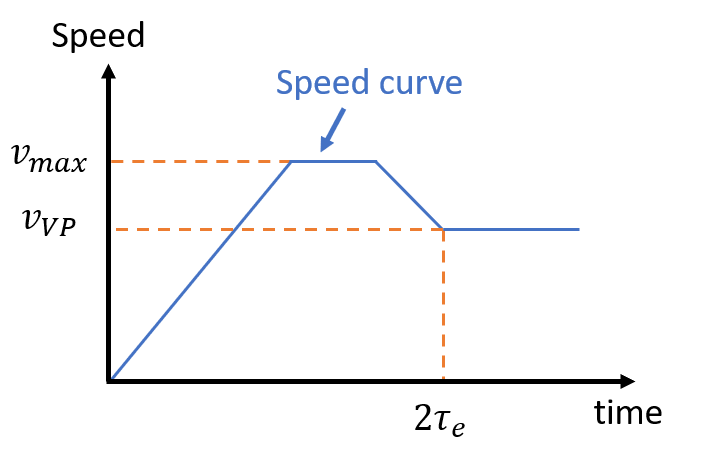}\label{fig:curve_2}} 
	\caption{Speed curves in the acceleration process}
	\label{fig:curve_b1}
\end{figure}

We first examine the acceleration process of the robot after a turn. It needs to cover a distance of $2D$ within $2\tau_e$, reaching a final velocity of $v_{VP}$. We discuss two cases regarding the maximum travel distance the robot can cover while satisfying the final velocity $v_{VP}$, which are shown in Figure B-\ref{fig:curve_1} and B-\ref{fig:curve_2}, respectively.

\noindent CASE 1:  $2v_{max}-v_{VP} \geq 2\tau_e c_{max}$

The speed curve $v_1(t)$ and the associated travel distance $d_1$ are formulated as follows:
\begin{eqnarray*}
    v_1(t) = \left\{
    \begin{aligned} 
    &c_{max}t,\quad t \leq \tau_e + \frac{v_{VP}}{2c_{max}}\\
    &v_{VP} + c_{max}(2\tau_e-t),\quad  \tau_e + \frac{v_{VP}}{2c_{max}} < t \leq 2\tau_e\\
    \end{aligned}
    \right.
\end{eqnarray*}
\begin{eqnarray*}
    d_1 = \int_0^{2\tau_e} v_1(t) = \frac{v_{VP}^2}{4c_{max}}  + \tau_e v_{VP} + c_{max}\tau_e^2
\end{eqnarray*}
\noindent CASE 2:  $2v_{max}-v_{VP} < 2\tau_e c_{max}$

The speed curve $v_1(t)$ and the associated travel distance $d_1$ are formulated as follows:
\begin{eqnarray*}
    v_2(t) = \left\{
    \begin{aligned} 
    &c_{max}t,\quad t \leq \frac{v_{max}}{c_{max}}\\
    &v_{max},\quad \frac{v_{max}}{c_{max}} < t \leq 2\tau_e - \frac{v_{max}-v_{VP}}{c_{max}} \\
    &v_{VP} + c_{max}(2\tau_e-t),\quad  2\tau_e - \frac{v_{max}-v_{VP}}{c_{max}} < t \leq 2\tau_e\\
    \end{aligned}
    \right.
\end{eqnarray*}
\begin{eqnarray*}
    d_2 = \int_0^{2\tau_e} v_2(t) = 2\tau_e v_{max} -\frac{2v_{max}^2-2v_{max}v_{VP}+v_{VP}^2}{2c_{max}}
\end{eqnarray*}
Above all, the maximum travel distance is expressed by:
\begin{eqnarray*}
    d(\tau_e,v_{VP},v_{max},c_{max}) = \left\{
    \begin{aligned} 
    &d_1,\quad 2v_{max}-v_{VP} \geq 2\tau_e c_{max}\\
    &d_2,\quad 2v_{max}-v_{VP} < 2\tau_e c_{max}\\
    \end{aligned}
    \right.
\end{eqnarray*}
Similarly, the analysis for the deceleration process before a turn is consistent with that of the acceleration process, yielding the same result.
\begin{flushright} 
	$\square$ 
\end{flushright}

\subsection*{B.4 Derivation of Average Travel Distance in Section 5}
Let $I_x$ represents the set of x-coordinates of loading stations on the top and bottom sides, while $I_y$ represents the set of y-coordinates of loading stations on the left and right sides. Similarly, $W_x$ and $W_y$ denote the set of all coordinates of potential locations on the top and bottom sides, and on the left and right sides, respectively. 

\noindent Area 1:
Select turning points with equal probability in the central area and perform a weighted sum of all path lengths, with weights proportional to the route lengths.
\begin{eqnarray*}
	E[l_1] &=& \frac{2D}{\lvert I_x \rvert\cdot \lvert I_y \rvert}\cdot \sum_{x\in I_x}\sum_{y\in I_y}\{\frac{{(x+y)}^2}{2(n_h+n_v)}+\frac{{[x+(n_h-y)]}^2}{2(n_h+n_v)}+\frac{{[(n_v-x)+y]}^2}{2(n_h+n_v)}\\
    &&+\frac{{[(n_v-x)+(n_h-y)]}^2}{2(n_h+n_v)}\} \\
	&=&2D\cdot [\frac{\frac{9+a^2}{6}(n_h^2+n_v^2)-n_h n_v-\frac{1}{3}}{n_h+n_v}+1]
\end{eqnarray*}

\noindent Area 2:
Calculate the average path lengths for starting points on the top and bottom sides and on the left and right sides, separately. According to the number of loading stations, perform a weighted sum. Based on the distance between the first turning point and the starting point, paths are divided into those that return to the original side and those that reach the opposite side, with each scenario having a probability of 1/2. For the former, the vertical movement distance is $n_v/2$ or $n_h/2$, while for the latter, it is $n_v$ or $n_h$, depending on the start point. 
\begin{eqnarray*}
	E[l_2] &=& 2D\cdot \big\{ \frac{n_v}{n_h+n_v}\sum_{x\in I_x}\sum_{z\in I_x} \frac{\frac{1}{2}(\frac{n_h}{2}+\lvert z-x \rvert) + \frac{1}{2}(n_h+\lvert z-x \rvert)}{{\lvert I_x \rvert}^2} \\
	&&+ \frac{n_h}{n_h+n_v}\sum_{y\in I_y}\sum_{z\in I_y} \frac{\frac{1}{2}(\frac{n_v}{2}+\lvert z-y \rvert) + \frac{1}{2}(n_v+\lvert z-y \rvert)}{{\lvert I_y \rvert}^2} \big\} \\
	&=& 2D\cdot [\frac{\alpha(n_h^2+n_v^2)}{3(n_h+n_v)}-\frac{2}{3\alpha(n_h+n_v)}+\frac{3n_h n_v}{2(n_h + n_v)}]\\
 \end{eqnarray*}

 \noindent Area 3:
We first calculate the basic average path length. In accordance with workload balance, the starting and ending points of the path should be equally likely to fall on each loading station. Then add the detour lengths for the two different scenarios multiplied by their respective probabilities, according to the path allocation rules.
 \begin{eqnarray*}
	E[l_3] &=& 2D\cdot \sum_{x\in I_x}\sum_{y\in I_y} [\frac{x+y}{\lvert I_x \rvert\cdot \lvert I_y \rvert} + \frac{n_h}{2(n_h+n_v)}\cdot \frac{\sum_{x\in W_x\backslash I_x} 2(x-\frac{\alpha n_v}{2})}{\lvert W_x\backslash I_x \rvert}\\
 &&+ \frac{n_v}{2(n_h+n_v)}\cdot \frac{\sum_{y\in W_y\backslash I_y} 2(y-\frac{\alpha n_h}{2})}{\lvert W_y\backslash I_y \rvert}] \\
	&=& 2D\cdot (\frac{n_h+n_v}{2} + \frac{(1+\alpha)}{4}\cdot \frac{n_h n_v}{n_h + n_v})
\end{eqnarray*}

\subsection*{B.5 Proof of Lemma 2}
We first prove (\romannumeral1). By equation (\ref{eq:travel distance}), the expression of average travel distance of a sorting task is as follows:
\begin{align}
	\bar{l}(n_h,n_v,\alpha) =& 2D\cdot[\frac{n_h^2+n_v^2}{n_h+n_v}\cdot\frac{-\alpha^4-5\alpha^2+18\alpha}{12} + \frac{n_h n_v}{n_h+n_v}\cdot\frac{-\alpha^2+2\alpha}{4}\nonumber\\
	&+ \frac{1}{n_h+n_v}\cdot\frac{3\alpha^3-\alpha^2-3\alpha-5}{12} + (n_h+n_v)\cdot \frac{\alpha^2-2\alpha+1}{2} +\frac{-\alpha^2 + 2\alpha}{2}] \nonumber
\end{align}
where $\alpha$ is the ratio of the number of workers to the number of aisle entrances, $\alpha\in (0,1]$. We consider the partial derivative:
\begin{align}
	\frac{\partial \bar{l}(n_h,n_v,\alpha)}{\partial \alpha} =& 2D\cdot[\frac{n_h^2+n_v^2}{n_h+n_v}\cdot \frac{-4\alpha^3-10\alpha+18}{12} + \frac{n_h n_v}{n_h+n_v}\cdot \frac{-\alpha+1}{2} \nonumber \\
    & + \frac{1}{n_h+n_v}\cdot \frac{9\alpha^2-2\alpha-3}{12}  + (n_h+n_v)(\alpha-1)-\alpha+1] \nonumber
\end{align}
At the two endpoints of the range of $\alpha$, the derivative has the value:\begin{align}
	\lim_{\alpha\xrightarrow{}0} \frac{\partial \bar{l}(n_h,n_v,\alpha)}{\partial \alpha} &= 2D\cdot[\frac{3(n_h^2+n_v^2)}{2(n_h+n_v)} +\frac{n_h n_v}{2(n_h+n_v)} - \frac{1}{4(n_h+n_v)}-(n_h+n_v)+1]\nonumber\\
	&\leq 2D\cdot[\frac{n_h^2+n_v^2-3n_h n_v}{2(n_h+n_v)} + 1] < 0\nonumber\\
        &\leq 2D\cdot[-\frac{1}{4}min\{n_h,n_v\}+1] \leq 0 \nonumber
\end{align}
\begin{align}
	\lim_{\alpha\xrightarrow{}1} \frac{\partial \bar{l}(n_h,n_v,\alpha)}{\partial \alpha} &= 2D\cdot[\frac{n_h^2+n_v^2}{3(n_h+n_v)} +\frac{1}{3(n_h+n_v)}-1] \nonumber\\
 &\geq 2D\cdot[\frac{1}{3}min\{n_h,n_v\}-1] > 0\nonumber
\end{align}
Similarly, we could calculate the second-order derivative and obtain that $\frac{\partial^2 \bar{l}(n_h,n_v,\alpha)}{\partial \alpha^2} \geq 0$. Above all, we can prove that $\bar{l}(n_h,n_v,\alpha)$ decreases initially and then increases.

Next, we prove the advantage of square-shape site in (\romannumeral2), namely the network with $n_h=n_v$. Without loss of generality, we assume $n_h=k\cdot n_v$, $k \geq 1$. Our objective is to show that the average travel distance is minimum when $k=1$. Let the area of the sorting zone be $S\cdot 4D^2$, then $n_h = \sqrt{kS}$, $n_v = \sqrt{S/k}$. Consider the derivative:
\begin{align}
\frac{\partial \bar{l}(S,k,\alpha)}{\partial k} =& 2D\cdot[\frac{k^2-1}{2k^{\frac{3}{2}}}\cdot \sqrt{S}\cdot\frac{-\alpha^4-5\alpha^2+18\alpha}{6} + \frac{k-1}{2k^{\frac{1}{2}}{(k+1)}^2}\cdot \sqrt{S}\cdot\frac{-\alpha^2+2\alpha}{4}\nonumber\\
	&+ \frac{k-1}{2k^{\frac{1}{2}}{(k+1)}^2}\cdot\frac{1}{\sqrt{S}}\cdot\frac{3\alpha^3-\alpha^2-3\alpha-5}{12} + \frac{k-1}{2k^{\frac{3}{2}}}\cdot\sqrt{S}\cdot \frac{\alpha^2-2\alpha+1}{2}] \nonumber
\end{align}
It is easy to prove that the right side is consistently non-negative. As a result, the derivative is non-negative when $k\geq 1$, thus $\bar{l}(S,k,\alpha)\geq \bar{l}(S,1,\alpha)$.

To prove (\romannumeral3), we first obtain the upper bound of $\bar{l}(n_h,n_v,\alpha)$. By (\romannumeral1), the upper bound can only be attained at the two endpoints:
\begin{align}
	\bar{l}(n_h,n_v,0) =& n_h \nonumber\\
	\bar{l}(n_h,n_v,1) =& \frac{9}{8} n_h - \frac{1}{4n_h} +\frac{1}{2} < \frac{9}{8} n_h \nonumber
\end{align}
While $\alpha=0$ is not feasible for operation, we conclude $\bar{l}(n_h,n_v,\alpha) < \frac{9}{8} n_h$. According to the sorting demands with an average distribution, we can easily derive that the average distance for all three paths are greater than $n_h$, then we have: $\bar{l}(n_h,n_v,\alpha) < n_h$. The proof is completed. 
\begin{flushright} 
	$\square$ 
\end{flushright}

\subsection*{B.6 Proof of Proposition 3}
We first obtain the expression of $\tilde{T}_M(n_h,n_v,n_l)$: 
\begin{align}
	\tilde{T}_M(n_h,n_v,n_l) = \frac{D}{\tau_e}\cdot\frac{n_h n_v}{a+b(n_h+n_v)}\cdot \frac{\frac{2n_l}{n_h+n_v} - {(\frac{n_l}{n_h+n_v})}^2}{\bar{l}(n_h,n_v,\frac{n_l}{n_h+n_v})} \nonumber
\end{align}
We further denote:
\begin{align}
	G_{n_h,n_v}(\alpha) = \frac{2\alpha - \alpha^2}{\bar{l}(n_h,n_v,\alpha)} \nonumber
\end{align}
\begin{align}
	\frac{\partial G_{n_h,n_v}(\alpha)}{\partial \alpha} = \frac{2(1-\alpha)\cdot \bar{l}(\cdot) - \frac{\partial \bar{l}(\cdot)}{\partial \alpha}\cdot (2\alpha - \alpha^2)}{\bar{l}^2(\cdot)} \nonumber
\end{align}
From Lemma 1 and Lemma 2, we have:
\begin{align}
	\lim_{\alpha \xrightarrow{} 0} \frac{\partial G_{n_h,n_v}(\alpha)}{\partial \alpha} = \frac{2}{\lim_{\alpha \xrightarrow{} 0} \bar{l}(\cdot)} \geq \frac{2}{\sqrt{n_h n_v}} > 0 \nonumber \\
	\lim_{\alpha \xrightarrow{} 1} \frac{\partial G_{n_h,n_v}(\alpha)}{\partial \alpha} = \frac{-\frac{\partial \bar{l}(\cdot)}{\partial \alpha}}{\lim_{\alpha \xrightarrow{} 0} \bar{l}^2(\cdot)} <0  \nonumber
\end{align}
The non-negativity of the second-order derivative of $G_{n_v,n_h}(\alpha)$ can be easily proved by obtaining the expression and checking the bounds of each term in the numerator. Hence, $G_{n_v,n_h}(\alpha)$ has one zero point. The expression of the derivative of $\tilde{T}_M(n_h,n_v,n_l)$ can be written as:
\begin{align}
	\frac{\partial \tilde{T}_M(n_h,n_v,n_l)}{\partial n_l} = \frac{D}{\tau_e}\cdot\frac{n_h n_v}{a+b(n_h+n_v)}\cdot \frac{1}{n_h+n_v}\cdot \frac{\partial G_{n_h,n_v}(\alpha)}{\partial \alpha} \nonumber
\end{align}
It is immediate to prove that $\tilde{T}_M(n_h,n_v,n_l)$ has one zero point within the range $n_l\in (0,n_h+n_v]$ and it initially increases then decreases.
\begin{flushright} 
	$\square$ 
\end{flushright}
	\section*{Appendix C: Penalty Successive Linear Programming}

\setcounter{equation}{0}
\renewcommand{\theequation}{C-\arabic{equation}}
\setcounter{algorithm}{0}
\renewcommand{\thealgorithm}{C-\arabic{algorithm}}
The classic penalty successive linear programming (PSLP) enjoys good robustness and convergence properties for large-scale problems \citep{bazaraa2013nonlinear}. Specifically, PSLP sequentially solves a linearized feasible direction finding subproblem along with the penalty function, and utilizes the concept of trust region (updated at each iteration) to control the step size. In each iteration $k$, a direction-finding linear program is formulated based on first-order Taylor series approximations to the objective and constraint functions, in addition to appropriate trust region restrictions on the direction components. The subproblem in iteration $k$ are formed as follows:
\begin{align}
	\textbf{(LP-S)}(\omega_k, \Delta_k) &\nonumber\\
	&\underset{d}{\min}\ \nabla {C_{d}(\omega_k)}^T\cdot d + \mu \cdot (\sum_{i=1}^{2\left|S\right|+2} y_i) \nonumber \\
	s.t.\quad & y_i \geq g_i(\omega_k) + \nabla {g_i(\omega_k)}^T d, \quad i = 1,2,\cdots,2\left|S\right|+2 \\
	& -\Delta_k \leq d \leq \Delta_k\\
	& y_i \geq 0, \quad i = 1,2,\cdots,2\left|S\right|+2
\end{align}
where $\omega = \{n_h,n_v, n_w^\sigma,n_r^\sigma\}$ is the set of all decision variables. $\Delta_k$ is the bound of $d$ in iteration $k$. $g_i(\omega_k)$ represents the left-hand side of standard form inequality constraints (\ref{gi_1})-(\ref{gi_4}) derived from the original \textbf{LDP}, assuming that the number of robots is less than the available VPs (considering the parking demand of robots during downtime, this setting is reasonable):
\begin{align}
	& - (n_h-1)(n_v-1) + N_{o} \leq 0 \label{gi_1}\\
	& \underset{\sigma}{max}\{n_{w}^\sigma\} - n_h - n_v \leq 0 \label{gi_2}\\
	& n_{r}^\sigma - \beta(n_h,n_v)\cdot (1-{(1-\frac{n_w^\sigma}{n_h+n_v})}^2)\cdot \frac{2\tau_e \cdot [n_h(n_v-1) + n_v(n_h-1)]}{\tau_c} \leq 0 \qquad \forall \sigma \in \mathcal{S} \label{gi_3}\\
	& -\frac{D \cdot n_{r}^\sigma}{\tau_e \cdot E[l(n_h,n_v,\frac{n_w^\sigma}{n_h+n_v})]} + {T}^\sigma \leq 0 \qquad \forall \sigma \in \mathcal{S} \label{gi_4}
\end{align}  

To avoid waste of resources, the number of loading station is always equal to the peak number of workers, thus $n_{l}$ is replaced by $\underset{\sigma}{max}\{n_w^\sigma\}$ in the constraints. The details of solving PLSP is shown in Algorithm~\ref{alg6.1}. It could be roughly divided into 2 steps in each iteration: (1) Generate \textbf{LP-S} and obtain the step size $d_k$, determine whether to stop; (2) Adjust step bounds for the next iteration. 

\begin{algorithm}[!htb]
	\caption{The PLSP algorithm to solve model \textbf{LDP}}
	\label{alg6.1}
	\small
	\begin{algorithmic}[1]
		\REQUIRE ~~\\ 
		An initial feasible solution $\omega_0$\\
		Confidence intervals $0 \leq \Delta_{LB} \leq \Delta_1$\\
		Parameters $0\leq \rho_0 \leq \rho_1 \leq \rho_2 \leq 1$, $\phi \leq 1$\\
		Large enough constant $\mu$\\
		Maximum iterations $N_{m}$\\
		\ENSURE The optimal value of layout decision variables, $\omega$\\
		\STATE Initialization: $k=1$, $R_k = 0$, $iteration = 1$
		\WHILE{$iteration < N_m$}
		\REPEAT
		\STATE Solve $\textbf{LP-S}(\omega_k, \Delta_k)$ to obtain $d_k$
		\IF {$C_{E L_k}(\omega_k) - C_{E L_k}(\omega_k + d_k) = 0$ \OR $\lceil\omega_k + d_k \rceil = \omega_k$}
		\STATE \textbf{Return} $\omega_k$
		\STATE \textbf{Stop}
		\ELSE 
		\STATE Calculate $R_k(\omega_k, d_k)$
		\ENDIF
		\IF{$R_k < \rho_0$}
		\STATE $\Delta_k = \alpha \Delta_k$
		\ENDIF
		\UNTIL{$R_k \geq \rho_0$}
		\STATE $\omega_{k+1} = \lceil\omega_k + d_k \rceil$
		\IF{$\rho_0 \leq R_k < \rho_1$}
		\STATE $\Delta_{k+1} = \phi \Delta_k$
		\ELSIF{$\rho_1 \leq R_k < \rho_2$}
		\STATE $\Delta_{k+1} = \Delta_k$
		\ELSIF{$R_k \geq \rho_2$}
		\STATE $\Delta_{k+1} = \Delta_k / \phi$
		\ENDIF
		\STATE $\Delta_{k+1} = max\{\Delta_{k+1},\Delta_{LB}\}$
		\STATE $k = k+1$
		\ENDWHILE
	\end{algorithmic}
\end{algorithm} 
In iteration $k$, the termination of algorithm could be determined by calculating the intermediate variables as follows:
\begin{align}
	&C_E(\omega_k) = C_{d}(\omega_k) + \mu\sum_{i} max\{0,g_i(\omega_k)\}\\
	&C_{E L_k}(\omega_k) =  C_{d}(\omega_k) + \nabla {C_{d}(\omega_k)}^T d_k + \mu \sum_{i} max\{0,g_i(\omega_k) + {\nabla g_i(\omega_k)}^T \cdot d_k\}\\
	&R_k = \frac{C_E(\omega_k) - C_E(\omega_k + d_k)}{C_{E L_k}(\omega_k) - C_{E L_k}(\omega_k + d_k)}
\end{align}
where $\mu$ is a large enough constant. Given that the variables in the original problem are integers, we introduced a rounding-up step for $\omega_k$ and the bounds in the search process. Additionally, due to the relatively limited solution space, the number of iterations in experiments remained below 50, and the total solving time was in the order of seconds. To ensure the feasibility of the solution for the original problem, the numerical update process includes a rounding-up step.

\end{document}